\documentclass[11pt]{article}

\usepackage[margin=.7in]{geometry}%modify the margins
\usepackage{amsfonts}%for mathcal, mathbb etc..
\usepackage[mathcal]{euscript}%for mathcal, mathbb etc..
\usepackage{amsmath}%for eqref
\usepackage{amsthm}%for proofs
\usepackage{graphicx}
\usepackage{xcolor}% to use colors

\usepackage{cleveref}

\usepackage{enumitem}%change the enumerate style

\usepackage{pifont}%ding symbols

\usepackage{pgf,tikz}
\usetikzlibrary{matrix, calc, arrows,shapes,decorations,plotmarks}

\newcommand{\varx}{x}% \varx\in\mathbb R^3, so maybe bold?
\newcommand{\varb}{\beta}% \varb\in\mathbb N_0^3, for index of Taylor expansion coefficients
\newcommand{\pG}{x_C}% center of the Taylor expansion
\newcommand{\Lc}{\mathcal L_{c}}
\newcommand{\LcAm}{\mathcal L_{c}^{Am,\Lambda}}
\newcommand{\LcPh}{\mathcal L_{c}^{Ph}}
\newcommand{\TC}{T_}
\newcommand{\polP}{P}% polynomial function P(\varx)
\newcommand{\polQ}{Q} % Q(\varx)
\newcommand{\polR}{R} % Q(\varx)
\newcommand{\funcci}{c_i}% the function (coef of PDE)  c_i\varx)
\newcommand{\funcc}{c_}% the function (coef of PDE)  c_\alpha(\varx)
\newcommand{\coefc}[1]{T_{c_{#1}}}% the function (coef of PDE)  c_\alpha(\varb)
\newcommand{\bfz}{\mathbf 0}
\newcommand{\funcczr}{c_{\bfz}}% the function (coef of PDE)   c_{\bfz}(\varx)
\newcommand{\coefczr}{T_{c_{\bfz}}}% the function (coef of PDE)  c_\alpha(\varb)

\newcommand{\dxk}{\partial_{x_k}}
\newcommand{\dxl}{\partial_{x_l}}

\newcommand{\mi}{\textrm i}
\newcommand{\bfm}{\mathbf M}
\newcommand{\CoPo}{\mathbb A}% space of complex polynomials in two variables
\newcommand{\HoPo}[1]{\mathbb A_{#1}}% space of homogeneous complex polynomials in two variables
\newcommand{\indn}{{\mathfrak n}}%a generic index: neither $i$ (for unknowns) nor $\varb$ (for equations)
\newcommand{\cst}{\mathfrak s}%initialization constant
\newcommand{\dir}{\mathbf d}%initialization direction
\newcommand{\Cmat}{\mathsf C}%initialization - PDop matrix
\newcommand{\Pmat}{\mathcal P}%initialization - PDop matrix orth. change of basis
\newcommand{\Dmat}{\mathcal D}%initialization - PDop diag of eig.vecs

\newcommand{\mata}{\mathsf A}%{M\hspace{-.75em}A}
\newcommand{\lata}{\mathsf L_A}%{L\!\!A}
\newcommand{\matp}{\mathsf P}%{M\!\!\!P}
\newcommand{\latp}{\mathsf L_P}%{L\!\!\!\!P}
\newcommand{\matq}{\mathsf Q}%}{M\!\!\!Q}
\newcommand{\matr}{\mathsf R}%{M\!\!\!R}
\newcommand{\mats}{\mathsf S}
\newcommand{\mate}{\mathsf E}%{M\!\!\!E}

\newcommand{\Fspace}{\mathbb F_{n,\Lc}}
\newcommand{\Kspace}{\mathbb K_{n,\Lc}}

\newtheorem{theorem}{Theorem}
\newtheorem{cor}{Corollary}
\newtheorem{prop}{Proposition}
\newtheorem{dfn}{Definition}
\newtheorem{lmm}{Lemma}
\newtheorem{hyp}{Hypothesis}
\newtheorem{rmk}{Remark}

\crefname{cor}{corollary}{corollaries}
\Crefname{cor}{Corollary}{Corollaries}

\crefname{prop}{proposition}{propositions}
\Crefname{prop}{Proposition}{Propositions}

\crefname{lmm}{lemma}{lemmas}
\Crefname{lmm}{Lemma}{Lemmas}

\usepackage{algorithmicx}
\usepackage{algorithm}% http://ctan.org/pkg/algorithms
\usepackage[noend]{algpseudocode}% http://ctan.org/pkg/algorithmicx

\usepackage{xcolor} 
\usepackage{tikz,pgfplots}
\pgfplotsset{compat=newest,
legend style={
font=\footnotesize,
rounded corners=2pt,
}}
\pgfplotscreateplotcyclelist{mycolorlist}{%
blue,every mark/.append style={fill=blue!80!black},mark=*\\%
red,every mark/.append style={fill=red!80!black},mark=square*\\%
brown!60!black,every mark/.append style={fill=brown!80!black},mark=otimes*\\%
black,mark=star\\%
blue,every mark/.append style={fill=blue!80!black},mark=diamond*\\%
red,densely dashed,every mark/.append style={solid,fill=red!80!black},mark=*\\%
brown!60!black,densely dashed,every mark/.append style={
solid,fill=brown!80!black},mark=square*\\%
black,densely dashed,every mark/.append style={solid,fill=gray},mark=otimes*\\%
blue,densely dashed,mark=star,every mark/.append style=solid\\%
red,densely dashed,every mark/.append style={solid,fill=red!80!black},mark=diamond*\\%
}

\pgfplotscreateplotcyclelist{nequalstwenty}{%
blue,every mark/.append style={fill=blue!80!black},mark=*\\%
red,every mark/.append style={fill=red!80!black},mark=square*\\%
brown!60!black,every mark/.append style={fill=brown!80!black},mark=otimes*\\%
black,mark=star\\%
green,every mark/.append style={fill=green!80!green},mark=diamond*\\%
pink,every mark/.append style={solid,fill=pink!80!pink},mark=*\\%
yellow!60!yellow,every mark/.append style={solid,fill=yellow!80!yellow},mark=square*\\%
black,every mark/.append style={solid,fill=gray},mark=otimes*\\%
blue,mark=star,every mark/.append style=solid\\%
red,every mark/.append style={solid,fill=red!80!black},mark=diamond*\\%
blue,every mark/.append style={fill=red!80!red},mark=*\\%
red,every mark/.append style={fill=blue!80!blue},mark=square*\\%
brown!60!black,every mark/.append style={fill=brown!80!black},mark=otimes*\\%
green,every mark/.append style={fill=blue!80!green},mark=diamond*\\%
pink,every mark/.append style={solid,fill=red!80!pink},mark=*\\%
yellow!60!yellow,every mark/.append style={solid,fill=brown!80!black},mark=square*\\%
black,every mark/.append style={solid,fill=red},mark=*\\%
green,every mark/.append style={fill=green!80!green},mark=star\\%
pink,every mark/.append style={solid,fill=pink!80!pink},mark=star\\%
yellow!60!yellow,every mark/.append style={solid,fill=yellow!80!yellow},mark=star\\%
}

\pgfplotscreateplotcyclelist{twotimeseight}{%
blue\\%blue,every mark/.append style={fill=blue!80!black},mark=*\\%
red\\%red,every mark/.append style={fill=red!80!black},mark=square*\\%
brown!60!black\\%brown!60!black,every mark/.append style={fill=brown!80!black},mark=otimes*\\%
black\\%black,mark=star\\%
green\\%green,every mark/.append style={fill=green!80!green},mark=diamond*\\%
pink\\%pink,every mark/.append style={solid,fill=pink!80!pink},mark=*\\%
yellow!60!yellow\\%yellow!60!yellow,every mark/.append style={solid,fill=yellow!80!yellow},mark=square*\\%
teal\\%black,every mark/.append style={solid,fill=gray},mark=otimes*\\%
blue,only marks,mark=star\\%blue,only marks,every mark/.append style={fill=red!80!red},mark=*\\%
red,only marks,mark=star\\%red,only marks,every mark/.append style={fill=blue!80!blue},mark=square*\\%
brown,only marks,mark=star\\%brown!60!black,only marks,every mark/.append style={fill=brown!80!black},mark=otimes*\\%
black,only marks,mark=star\\%green,only marks,every mark/.append style={fill=blue!80!green},mark=diamond*\\%
green,only marks,mark=star\\%pink,only marks,every mark/.append style={solid,fill=red!80!pink},mark=*\\%
pink,only marks,mark=star\\%yellow!60!yellow,only marks,every mark/.append style={solid,fill=brown!80!black},mark=square*\\%
yellow,only marks,mark=star\\%black,only marks,every mark/.append style={solid,fill=red},mark=*\\%
teal,only marks,mark=star\\%green,only marks,every mark/.append style={fill=green!80!green},mark=star\\%
}

\pgfplotscreateplotcyclelist{threetimeseight}{%
blue, mark size=4pt,every mark/.append style={fill opacity = 0},only marks,mark=*\\
red, mark size=4pt,every mark/.append style={fill opacity = 0},only marks,mark=*\\
brown, mark size=4pt,every mark/.append style={fill opacity = 0},only marks,mark=*\\
black, mark size=4pt,every mark/.append style={fill opacity = 0},only marks,mark=*\\
green, mark size=4pt,every mark/.append style={fill opacity = 0},only marks,mark=*\\
pink, mark size=4pt,every mark/.append style={fill opacity = 0},only marks,mark=*\\
yellow, mark size=4pt,every mark/.append style={fill opacity = 0},only marks,mark=*\\
teal, mark size=4pt,every mark/.append style={fill opacity = 0},only marks,mark=*\\
blue,only marks,mark=star\\
red,only marks,mark=star\\
brown,only marks,mark=star\\
black,only marks,mark=star\\
green,only marks,mark=star\\
pink,only marks,mark=star\\
yellow,only marks,mark=star\\
teal,only marks,mark=star\\
blue\\
red\\
brown\\
black\\
green\\
pink\\
yellow\\
teal\\
}

\title{%Computational aspects of GPWs for the 3D convected Helmholtz equation\\
Three types of quasi-Trefftz functions for the 3D convected Helmholtz equation: construction and approximation properties}

\author{Lise-Marie Imbert-G\'erard\footnote{Email: lmig@math.arizona.edu Affiliation: University of Arizona, Department of Mathematics, Tucson AZ, USA. ORCID number 0000000177548582. L.-M. Imbert-G\'erard acknowledges support from the US National Science Foundation: this material is based upon work supported by the NSF under Grants No. DMS-2105487 and DMS-2110407.}, 
Guillaume Sylvand\footnote{Affilitation: Airbus Central R\&T, France.}}

\begin{document}
\maketitle
%\tableofcontents

\begin{abstract}
%Trefftz and quasi-Trefftz
Trefftz methods are numerical methods for the approximation of solutions to boundary and/or initial value problems. They are Galerkin methods with particular test and trial functions, which solve locally the governing partial differential equation (PDE). This property is called the Trefftz property. Quasi-Trefftz methods were introduced to leverage the advantages of Trefftz methods for problems governed by variable coefficient PDEs, by relaxing the Trefftz property into a so-called quasi-Trefftz property: test and trial functions are not exact solutions but rather local  approximate solutions to the governing PDE.
%3D convected Helmholtz eq
%quasi-Trefftz bases with good approximation properties
In order to develop quasi-Trefftz methods for aero-acoustics problems governed by the convected Helmholtz equation,
the present work tackles the question of the definition, construction and approximation properties of three families of quasi-Trefftz functions: two based on generalizations on plane wave solutions, and one polynomial.
{The polynomial basis shows significant promise as it does not suffer from the ill-conditioning issue inherent to wave-like bases.}
\end{abstract}

\section*{Acknowledgments}
L.-M. Imbert-G\'erard acknowledges support from the US National Science Foundation: this material is based upon work supported by the United States National Science Foundation under Grants No. DMS-2105487 and DMS-2110407.

\section*{Statements and declarations}
 L.-M. Imbert-G\'erard has disclosed an outside interest in Airbus Central R\&T to the University of Arizona.  Conflicts of interest resulting from this interest are being managed by The University of Arizona in accordance with its policies.

\section{Introduction}
Several time-harmonic wave propagation phenomena can be modeled by variations of the Helmholtz equation.
We are interested in developing tools for the numerical simulation of  linear acoustic propagation at a fixed frequency in a subsonic flow in three dimensions. % around a scattering object.
Such phenomena can be modeled by the convected Helmholtz equation for the acoustic potential:
\begin{equation}\label{CHinit}
-\nabla\cdot(\rho(\nabla\phi-(\bfm\cdot\nabla\phi)\bfm+\mi \kappa \phi\bfm))
-\rho(\kappa^2\phi+\mi \kappa \bfm\cdot\nabla\phi)=0,
\end{equation}
where $\rho$ is the real-valued fluid density and $\bfm$ is the vector-valued fluid velocity, both assumed to be depending on the space variable, while $\kappa$ is the wavenumber. 
 The density $\rho$ is naturally assumed to be positive. We are interested in the subsonic regime, so the Mach number $|\bfm|$ is assumed to be no more than $1$.
For comparison the standard Helmholtz equation then corresponds to the case of a fluid with velocity equal to zero and a constant density.
The methods of interest here are the so-called quasi-Trefftz methods that were introduced as an extension of Trefftz methods for problems of wave propagation in inhomogeneous media. %governed by variable-coefficient equations.

In a search for bounds on the solutions to boundary value problems, Trefftz introduced in the 1920s the idea to leverage trial functions satisfying the governing equation \cite{Trefftz,transl}.
 Since then, this idea has lead to the development of several numerical methods, generally referred to as Trefftz methods. 
 %{\color{red} cite a bunch, maybe early works Jirousek, Herrera, the more recent like the whole wave stuff and maybe Igor's stuff too - check out the survey} 
 In the present work, Trefftz methods refer to a class of numerical methods falling within the category of Galerkin methods, and specifically relying on functions that satisfy the governing equation, both for the derivation of their weak formulation and for the discretization of this formulation.
 Implementing Trefftz methods requires bases of %the knowledge of 
 exact solutions to the governing equation - called Trefftz functions - in order to discretize the weak formulation. 
This then limits their application since exact solutions are not known for many equations. Nevertheless these methods are particularly popular in the field of wave propagation, as circular, spherical, plane or even evanescent waves are some natural and common candidates Trefftz functions, see for instance \cite{moiolaTH,LuoT,survey} and references therein. 
Early works on such methods include the introduction \cite{uwvf,uwvfart} and study \cite{Huttuwvf,buffam,gabg,Hutt1,Hutt2,Hutt3} of the ultra-weak variational formulation,  as well as more general wave-based methods \cite{PWDGh,PWDGp,PWDGhp}.
More recent works include \cite{TVEM,TVEMnum} focusing on Trefftz Virtual Element Methods, and \cite{barucond,cong} focusing on conditioning aspects.
Existing extensions to space-time problems include work on acoustics and elasto-acoustics \cite{STTDGac,STTDGea}, as well as Friedrichs systems \cite{Morel1,Morel2}.

%OR MAYBE DON:T SSAY THAT JUST YET  Despite some drawbacks in particular in terms of conditioning, 
Intuitively, the advantage of Trefftz methods relies on their PDE-dependent function spaces: since Trefftz functions solve exactly the PDE, they can be expected to have better approximation properties than non-PDE dependent function spaces. In practice, local approximation properties of the discrete spaces are a corner stone of convergence properties of Galerkin methods. Typically these would be stated for a space $\mathbb V$ of functions to be approximated by functions in a discrete space $\mathbb V_h$ locally on some region $\Omega$ containing a point $\pG$ as:
$$
\forall u \in \mathbb V,
\exists u_a\in\mathbb V_h \text{ s. t. }
\forall \varx\in\Omega,
|u(\varx)-u_a(\varx)|\leq C | \varx-\pG|^{n}
$$
for some order of approximation $n$.
This is precisely where the advantage of Trefftz functions over standard polynomial spaces can be emphasized, as spaces of Trefftz functions require less degrees of freedom to achieve a given order of accuracy than standard polynomial spaces.
It is however important to keep in mind that these approximation properties of spaces of Trefftz functions hold only for the approximation of exact solutions to the governing PDE (as will be illustrated in \Cref{thm:approx}),
whereas approximation properties of standard polynomial spaces hold for the approximation of smooth enough function that do not necessarily solve the governing PDE. 
Moreover, bases of Trefftz functions are also known to suffer from ill-conditioning issues in certain regimes. {One of the type of basis functions proposed in this work represents an avenue to overcome these issues.}
% these aspects are beyond the scope of this work. \textcolor{we need to update this statement and announce here our great findings}

In the general context of wave propagation, the application of Trefftz methods in their standard form to problems of propagation through inhomogeneous media is similarly limited since one more time global exact solutions are not known for most variable-coefficient equations. 
Quasi-Trefftz methods, relying on approximate solutions - called quasi-Trefftz functions - rather than exact solutions to the governing equation, were introduced to extend Trefftz methods to problems governed by variable-coefficient equations. They were first introduced in \cite{LMD} under the name of Generalized Plane Wave (GPW) methods for 2D problems governed by the Helmholtz equation. 
%The name was derived from the original idea to 
The original idea behind the GPW concept was to retain the oscillating behavior of a plane wave (PW) while allowing for some extra degrees of freedom to be adapted to the varying PDE coefficient, and this is where their name came from. Initially this was performed via the introduction of Higher Order Terms (HOT) in the phase of a PW as follows:
\begin{equation}\label{eq:c/GPW}
\left\{\begin{array}{l}
\varphi(\varx)= \exp( i \widetilde{\kappa} %\kappa\sqrt{\epsilon(\pG)}
\mathbf k\cdot\varx+\text{HOT})\\
\left[-\Delta-\kappa^2\epsilon(\varx)\right]\varphi(\varx)\approx0,
\end{array}\right.
\text{ instead of }
\left\{\begin{array}{l}
\phi(\varx)=\exp i\kappa\mathbf k\cdot\varx\\
\left[-\Delta-\kappa^2\right]\phi=0,
\end{array}\right.
\text{ for any unit vector } \mathbf k\in\mathbb R^3,
\end{equation}
where Helmholtz equation has either a constant coefficient $\kappa^2$ or a variable coefficient $\kappa^2\epsilon$; so $\kappa$ is the wavenumber of the PW while $\widetilde{\kappa} %\kappa\sqrt{\epsilon(\pG)}
$ can be interpreted as the local wavenumber of the GPW.
A procedure to construct a basis of such GPWs was proposed in \cite{LMinterp}, together with a study of the approximation properties of the basis.
A systematic procedure to study these approximation properties was introduced in \cite{IGS}, and used in \cite{AbGPW} on a new type of GPWs.
The idea of the associated GPW-based Galerkin method was presented in \cite{LMD}, a proof of convergent for a variant was studied in \cite{LMM}. The method was applied a problem of mode conversion for wave propagation in plasmas in \cite{LMmc}.

Other works on quasi-Trefftz methods include  
\cite{Yuan} on the convergence of a GPW Discontinuous Galerkin method for anisotropic Helmholtz problems,
\cite{Morel2}, on linear transport problems, %Trefftz basis but also complicated BF construction
as well as \cite{IGMS} on time-dependent wave propagation problems. 
 
 Quasi-Trefftz methods rely on function spaces of approximate solutions to the governing equation, as opposed to exact solutions, and this is their fundamental property. In our work, we define this approximation as a local property in the sense of a Taylor expansion. 
Given the partial differential operator $\mathcal L$ of the governing equation, and any point $\pG$ in the domain of interest, we consider functions $\varphi$  with the following property:
\begin{equation}
\label{eq:TaylExp}
%%%\text{Find }\varphi \\%(\mathbf x)= \exp( i\kappa\mathbf k\cdot\mathbf x+\text{HOT})\\
%\left[-\Delta-\kappa^2\epsilon(\varx)\right]\varphi(\varx)=O(|\varx - \pG|^q), 
%%\text{ for some } q\in\mathbb N.
\mathcal L \varphi(\varx)=O(|\varx - \pG|^q), 
\end{equation}
for some parameter $q$ providing some flexibility in the desired order of approximation with respect to the distance $|\varx - \pG|$. In other words the degree $q-1$ Taylor polynomial of the image of $\varphi$ through the operator $\mathcal L$ is zero.
In the context of Discontinuous Galerkin methods, with function spaces of local functions defined element-wise on a computational mesh, then if \eqref{eq:TaylExp} holds within each element, with $\pG$ in the element, then the remainder can be described as $O(h^q)$ where $h$ denotes the mesh size, as $|\varx - \pG|\leq h$ for all $\varx$ in the element.
There the goal is to establish so-called $h$-convergence properties, that is in the regime $h\to 0$. 
Note that, in order to prove their convergence, quasi-Trefftz methods so far include in their weak formulations a stabilization term to handle  the non-zero remainder in the quasi-Trefftz property  \eqref{eq:TaylExp}, see \cite{LMM,Yuan,IGMS}.

\subsection{Central results}
Our goal is to address here the fundamental question of basis functions at the centre of quasi-Trefftz methods for the three-dimensional convected Helmholtz equation: the actual construction of basis functions is fundamental to the discretization stage, and therefore to the implementation of the methods, while the approximation properties of the discrete space are a fundamental element in the proof of convergence of the methods.
Since quasi-Trefftz functions satisfy a local quasi-Trefftz property \eqref{eq:TaylExp}, 
more precisely, the goal of this work is twofold:
\begin{enumerate}
\item developing algorithms for the construction of {\bf local} basis functions for quasi-Trefftz function spaces for the partial differential operator of the three-dimensional convected Helmholtz equation, $\mathbb V_h$, guaranteeing a limited computational cost for the practical construction,
see Sections \ref{sec:QTfams} and \ref{sec:norm},
\item studying the {\bf local} approximation property of the resulting spaces $\mathbb V_h$ in the following sense ; given $n\in\mathbb N$, there is a space $\mathbb V_h$  satisfying:
\begin{equation}\label{IntPb}
%\begin{array}{l}
\forall u%\in \mathcal C^\infty 
\text{ satisfying the governing PDE, }%\mathcal Lu=0,
\exists u_a\in\mathbb V_h \text{ s. t. }
\forall \varx\in\mathbb R^3,
|u(\varx)-u_a(\varx)|\leq C | \varx-\pG|^{n+1},
%\end{array}
\end{equation}
see Section \ref{sec:Approx}.
\end{enumerate}

Inspired by classical PWs $\exp \Lambda \cdot\Big(\varx-\pG \Big)$, for some $\Lambda\in\mathbb C^3$, we will focus on three different families of quasi-Trefftz functions:
\begin{itemize}
\item phase-based GPWs, following the original ansatz proposed in \cite{LMD} via the introduction of higher order terms in the phase of a PW, see \eqref{eq:c/GPW},
\item amplitude-based GPWs, following the ansatz proposed in \cite{AbGPW} via the introduction of higher order terms in the amplitude of a PW,
\item purely polynomial quasi-Trefftz functions, which so far we have only used for time-dependent wave propagation in \cite{IGMS}.
\end{itemize}
In each of the wave-based cases the ansatz is an extension of cases studied previously in two-dimensions, whereas the situation is different for the polynomial case. This is the first time that polynomial quasi-Trefftz functions are proposed for time-harmonic problems.

We will pursue the announced goals for these three families of quasi-Trefftz functions, highlighting the similarities and differences between the three cases. The fundamental contribution of this work is to show that these three families of quasi-Trefftz functions achieve the approximation properties with exactly the same number of degrees of freedom as their Trefftz function (wave-based) counterpart do for the constant-coefficient cases studied in the literature.
{The numerical results presented here show that it is possible to find quasi-Trefftz bases that, unlike GPW bases, {\bf do not suffer from the ill-conditioning problem inherent to wave-related bases}. As a consequence, it is a promising path for future development of Trefftz-like methods in the field of frequency-domain wave simulation.}
%While the corresponding quasi-Trefftz approximations of any given exact solution to the governing PDE will consequently behave similarly in the asymptotic regime, we expect them to behave very differently in the pre-asymptotic regime due to the very nature of polynomials

It is important to note that general time-harmonic wave-propagation equations have no exact polynomial solution, in other words there exist no polynomial Trefftz function in this case. However, there are more quasi-Trefftz functions than Trefftz functions, since the former are defined by a less restrictive constraint,  and as we will see it is possible to construct polynomial quasi-Trefftz functions.

\subsection{
%Hypothesis and notation/
Preliminaries}
\label{ssec:prelim}
Throughout this article, we will use the following notation.
The set of positive integers is denoted $\mathbb N$
and  
the set of non-negative integers is denoted $\mathbb N_0 :=\mathbb N\cup \{0\}$,
while the zero multi-index is denoted $\bfz = (0,0,0)$.
The canonical basis of $\mathbb R^3$ or $(\mathbb N_0)^3$ is denoted $\{e_k,k\in\{1,2, 3\}\}$, and $|\cdot |$ denotes the euclidean norm on $\mathbb R^3$.
Multi-indices in $(\mathbb N_0)^3$ are denoted $\alpha,\beta,\gamma,i,j$, the sum of multi-indices is defined as $i+j = (i_1+j_1,i_2+j_2,i_3+j_3)$ for all $i$ and $j $ in $(\mathbb N_0)^3$ 
while $|\cdot |$ also denotes the length for a multi-index,  that is $|i|=i_1+i_2+i_3$ for all $i\in(\mathbb N_0)^3$,
the factorial is denoted $i! = i_1!i_2!i_3!$,
and 
for the sake of compactness
$i\leq j$ means that $i_k\leq j_k$ for $k\in\{1,2,3\}$,
$i< j$ means that $i\neq j$ and $i_k\leq j_k$ for $k\in\{1,2,3\}$,
and  the linear order $\prec$ on $\mathbb N_0^3$ is defined by	
$$
\forall (\mu,\nu)\in\left(\mathbb N_0^3\right)^2,
 \mu\prec\nu 
 \Leftrightarrow
\left\{
 \begin{array}{l}
 |\mu|<|\nu|\ , \text{ or }\\
 |\mu|=|\nu| \text{ and } \mu_1<\nu_1\ , \text{ or }\\
 |\mu|=|\nu|,\  \mu_1=\nu_1\text{ and } \mu_2<\nu_2.
 \end{array}
 \right.
$$
A set of all indices with a common length $\ell$ will be referred to as the layer $\ell$; the layer $\ell=0$ contains one index, namely $\bfz$, and each layer $\ell>0$ contains exactly $(\ell+1)(\ell+2)/2$ indices since
$$
\sum_{|i|=\ell}1
 = \sum_{i_1=0}^\ell \sum_{i_2=0}^{\ell-i_1}1
 = \sum_{i_1=0}^\ell (\ell-i_1+1) 
% = \sum_{i_1=1}^{\ell+1} i_1
 = \frac{(\ell+1)(\ell+2)}2.
$$
We will make use of a particular type of index numbering, denoted $\mathcal N$, in order to evidence the structure of different objects of interest. Given $n\in\mathbb B$, we want this numbering to satisfy $|i|<|j|$ implies that $\mathcal N(i)<\mathcal N(j)$, hence we write it as:
$$
\forall i\in\left(\mathbb N_0\right)^3\text{ with } |i|\leq n,
\mathcal N(i) = \frac{|i|(|i|+1)(|i|+2)}{6}+\mathcal N_{|i|}(i),
$$
for some $\mathcal N_m$ providing a numbering of indices of length m. 
For instance we can choose $\mathcal N_m$ to count indices according to the linear order $\prec$ within the layer $m$ and in this case the numbering corresponds to $\mathcal N(i) = \sum_{j\prec i} 1$, 
or choose $\mathcal N_m(i)=(i_2+i_3)(i_2+i_3+1)/2+i_3+1$ for all $i$ such that $|i|=m$. 
The generic point of interest in the domain of the equation is denoted $\pG$.
The coefficient of a Taylor expansion in the neighborhood of $\pG$, for any $n\in\mathbb N$, are denoted $\TC f[\varb] = \frac1{\varb!}\partial_{\varx}^{\varb} f(\pG)$ for all function $f\in\mathcal C^n$ at $\pG$ with $\varb\in\mathbb N_0^3$, $|\varb|\leq n$. We chose to avoid an unnecessary explicit mention of $\pG$ in the $T$ notation because all Taylor expansion will be performed at $\pG$.
\begin{rmk}\label{rmk:SimpRules}
 Given $n\in\mathbb N_0$ and $\pG\in\mathbb R^3$, a few simple derivative rules can then be expressed in a compact way as follows:
 $$
 \left\{\begin{array}{l}
 \displaystyle
\TC{fg}[\varb] = \sum_{\gamma\in\mathbb N_0^3;\gamma\leq \varb} \TC f[\varb-\gamma] \TC g[\gamma]
\text{ if }f\text{ and }g\text{ are }\mathcal C^n\text{ at }\pG\text{ with }|\varb|\leq n,
\\
 \displaystyle
\TC{fgh}[\varb] 
= \sum_{\gamma\in\mathbb N_0^3;\gamma\leq \varb} \sum_{\eta\in\mathbb N_0^3;\eta\leq \gamma} \TC f[\varb-\gamma] \TC g[\eta] \TC h[\gamma-\eta],
\\
 \displaystyle
\TC{\partial_\varx^\alpha f}[\beta]=\frac{(\alpha+\varb)!}{{\beta}!} \TC{f}[\alpha+\varb] 
\text{ if  }f\text{ is }\mathcal C^n\text{ at }\pG\text{ with }|\alpha|+|\varb|\leq n,
\\ 
\TC{(X-\pG)^i}[\varb] = \delta(\varb-i).
\end{array}\right.
$$
%$\displaystyle
%\TC{fg}[\varb] = \sum_{\gamma\leq \varb} \TC f[\varb-\gamma] \TC g[\gamma]
%$ if $f$ and $g$ are $\mathcal C^n$ at $\pG$ with $|\varb|\leq n$,
%%\textcolor{blue}{
%$\displaystyle
%\TC{fgh}[\varb] 
%%= \sum_{\gamma\leq \varb} \TC f[\gamma] \TC{gh}[\varb-\gamma]
%= \sum_{\gamma\leq \varb} \sum_{\eta\leq \gamma} \TC f[\varb-\gamma] \TC g[\eta] \TC h[\gamma-\eta]
%$,
%%}
%$\TC{\partial_\varx^\alpha f}[\beta]=\frac{(\alpha+\varb)!}{{\beta}!} \TC{f}[\alpha+\varb] $ if  $f$ is $\mathcal C^n$ at $\pG$ with $|\alpha|+|\varb|\leq n$, 
%$\TC{(X-\pG)^i}[\varb] = \delta(\varb-i)$.
%%while $\coefP(\varb) = \lambda_\varb$, $\coefQ(\varb) = \mu_\varb$, as long as $|\beta|\leq \deg \polP$, $|\beta|\leq \deg \polQ$,
\end{rmk}

\begin{dfn}
A linear partial differential operator of order $2$, in three dimensions, with a given set of complex-valued functions  $c = \{\funcci ; i\in\mathbb N_0^3, |i|\leq 2\}$ will be denoted hereafter as
$$ \Lc := \sum_{i\in\mathbb N_0^3; |i|\leq 2} \funcci\left( \varx  \right) \partial_{\varx}^i,$$ 
where $\varx=(x_1,x_2,x_3)\in\mathbb R^3$ and $\partial_{\varx}^i = \partial_{x_1}^{i_1} \partial_{x_2}^{i_2} \partial_{x_3}^{i_3}$.
\end{dfn}
%
%The set $\{i\in\mathbb N_0^3, |i|\leq 2\}$ can be decomposed into subsets as follow:
%\begin{itemize}
%\item $
%\{|i|=2 \text{ with } i_k=2 \}$
%\item $
%\{|i|=2 \text{ with } i_k=1,i_l=1,k\neq l \}$
%\item $
%\{|i|=1 \text{ with } i_k = 1\}$
%\item $
%\{ \bfz \}$
%\end{itemize}
%
We will make use of the fact that the set $\{i\in\mathbb N_0^3, |i|\leq 2\}$ can be split as:
$$
\{i\in\mathbb N_0^3, |i|\leq 2\}
=
\{2e_k, 1\leq k\leq 3\}
\cup
\{e_k+e_{k'}, 1\leq k<k'\leq 3\}
\cup
\{e_k, 1\leq k\leq 3\}
\cup
\{ \bfz \}
$$ 
For instance, in the case of the convected Helmholtz equation, the variable coefficients of the partial differential operator can then be defined as follows:
%\begin{itemize}
%%\item $\displaystyle c_{e_k+e_l}=\rho(\bfm_k\bfm_l -\delta_{kl})$
%\item $\displaystyle c_{2e_k}=\rho\left((\bfm_k)^2-1\right)$ for $1\leq k\leq 3$,
%\item $\displaystyle c_{e_k+e_l}=\rho\bfm_k\bfm_l $ for $1\leq k<k'\leq 3$,
%\item $\displaystyle c_{e_k} =  \sum_l  \bfm_l\dxl\bfm_k +\nabla\cdot(\rho\bfm)\bfm_k-\dxk\rho-2\mi\kappa\rho\bfm_k$ for $1\leq k\leq 3$,
%\item  $\displaystyle c_ {\bfz} = -\nabla\cdot(\mi \kappa \rho\bfm)-\rho\kappa^2$.
%\end{itemize}
$$
\left\{\begin{array}{ll}
\displaystyle c_{2e_k}=\rho\left((\bfm_k)^2-1\right) &\text{for }1\leq k\leq 3,\\
\displaystyle c_{e_k+e_l}=\rho\bfm_k\bfm_l &\text{for }1\leq k<k'\leq 3,\\
\displaystyle c_{e_k} = \rho \sum_l  \bfm_l\dxl\bfm_k +\nabla\cdot(\rho\bfm)\bfm_k-\dxk\rho-2\mi\kappa\rho\bfm_k &\text{for }1\leq k\leq 3,\\
\displaystyle c_ {\bfz} = -\mi \kappa \nabla\cdot(\rho\bfm)-\rho\kappa^2.
\end{array}\right.
$$
We will also make use of the layer structure in multi-index space following the multi-index length, in particular we will leverage the fact that the set $\{i\in\mathbb N_0^3, |i|\leq d\}=\cup_{\ell=0}^d \{ i\in\mathbb N_0^3, |i|=\ell\}$ can be split as:
$$
\{i\in\mathbb N_0^3, |i|\leq d\}
=
{\color{magenta}\{\mathbf 0,e_k, 1\leq k\leq 3\}}
\cup
\cup_{\ell=2}^d
\left\{
{\color{red}\{ i\in\mathbb N_0^3, |i|=\ell, i_1\in\{0,1\}\}}
\cup
{\color{blue}\{ i\in\mathbb N_0^3, |i|=\ell, i_1\geq 2\}}
\right\}.
$$ 
For instance, Figure \ref{fig:indices} illustrates with the same color code some indices $i\in\mathbb N_0^3$ with $ |i|\leq 7$.
\begin{figure}[htb]\centering
\resizebox{.49\linewidth}{!}{
\begin{tikzpicture}
[rotate around x=-90,rotate around y=0,rotate around z=-30,grid/.style={very thin,gray}]
        \foreach \x in {0,1,...,7}
        \foreach \y in {0,1,...,7}
        \foreach \z in {0,1,...,7}
        {
            \draw[grid,lightgray] (\x,0,0) -- (\x,7,0);
            \draw[grid,lightgray] (0,\y,0) -- (7,\y,0);
            \draw[grid,lightgray] (0,\y,0) -- (0,\y,7);
            \draw[grid,lightgray] (0,0,\z) -- (0,7,\z);
        }
        \draw [->] (0,0) -- (8,0,0) node [right] {$\indn_{1}$};
        \draw [->] (0,0) -- (0,8,0) node [above] {$\indn_{2}$};
        \draw [->] (0,0) -- (0,0,8) node [below left] {$\indn_3$};
\draw(0,0)node[left]{$(0,0,0)$};
%% le rouge
         \draw (0,0,0) node[circle,fill=magenta,inner sep=2pt] {};
         \draw (1,0,0) node[circle,fill=magenta,inner sep=2pt] {};
         \draw (0,1,0) node[circle,fill=magenta,inner sep=2pt] {};
         \draw (0,0,1) node[circle,fill=magenta,inner sep=2pt] {};
        \draw[fill=magenta,opacity=.5] (0,1,0)--(0,0,1)--(1,0,0);
%le bleu
        \draw[fill=gray,opacity=.3] (0,7,0)--(0,0,7)--(7,0,0);
	\draw (0,0,7) node[left] {$7$};
	\draw (0,7,0) node [above] {$7$};
	\draw (7,0,0) node[below left] {$7$};
        \foreach \x in {0,1,...,7}
        {
        \pgfmathsetmacro{\int}{7-\x}
        \foreach \y in {0,...,\int}
        {
            \draw (\x,\y,7-\x-\y) node[circle,fill=blue,inner sep=2pt] {};
        }
        }
%un peu plus de rouge
        \foreach \x in {0,1}
        {
        \pgfmathsetmacro{\int}{7-\x}
        \foreach \y in {0,...,\int}
        {
            \draw (\x,\y,7-\x-\y) node[circle,fill=red,inner sep=2pt] {};
        }
        }
\end{tikzpicture}
}
\caption{Illustration of a multi-index layer in  $(\mathbb N_0)^3$. Within the layer $\{\indn\in\mathbb N_0^3, |\indn|=\ell\}$ for $\ell =7$, represented in gray, indices with $i_1\in\{0,1\}$ are represented in red while those with $i_1\geq 2$ are represented in blue. Indices with $|i|\leq 1$ are also represented in magenta.}
\label{fig:indices}
\end{figure}

%Announce important properties satisfied by the convected Helmholtz operator that will be key to our study: 
Beyond the convected Helmholtz equation, the work proposed in this article relies on a set of minimal hypotheses for the partial differential operator $\mathcal L_c$.
The first aspect will lead to the well-posedness of a subproblem in the construction of quasi-Trefftz functions, whereas the second will lead to the construction of a set of linearly independent quasi-Trefftz functions.
We gather the two in the following Hypothesis.
\begin{hyp}
\label{hyp:PDop}
Given a point $\pG\in\mathbb R^3$ and a set of complex-valued functions  $c = \{\funcci,i\in\mathbb N_0^3, |i|\leq 2\}$,
%regularity
the functions are 
assumed to be $\mathcal C^\infty$ at the point $\pG$, with
$c_{2e_1}(\pG)\neq 0$ and the matrix defined by:
%$$
%\mathcal C:=\begin{bmatrix}
%c_{2e_1}(\pG)& \frac12 c_{e_1+e_2}(\pG)  & \frac12 c_{e_1+e_3}(\pG) \\
% \frac12 c_{e_1+e_2}(\pG) & c_{2e_2}(\pG) & \frac12 c_{e_2+e_3}(\pG) \\
%  \frac12 c_{e_1+e_3}(\pG) &\frac12 c_{e_2+e_3 } (\pG)& c_{2e_3}(\pG)
%\end{bmatrix}
%or\begin{bmatrix}
%c_{2e_1}& c_{e_1+e_2} & \frac12 c_{e_1+e_3} \\
% \frac12 c_{e_1+e_2} & c_{2e_2} & \frac12 c_{e_2+e_3} \\
%  \frac12 c_{e_1+e_3} &\frac12 c_{e_2+e_3 } & c_{2e_3}
%\end{bmatrix}(\pG)
%$$
$$
\Cmat:=\begin{bmatrix}
 \coefc{2e_1}[\bfz] & \frac12 \coefc{e_1+e_2}[\bfz] & \frac12 \coefc{e_1+e_3}[\bfz] \\
 \frac12 \coefc{e_1+e_2}[\bfz] & \coefc{2e_2}[\bfz] & \frac12 \coefc{e_2+e_3}[\bfz] \\
  \frac12 \coefc{e_1+e_3}[\bfz] &\frac12 \coefc{e_2+e_3}[\bfz]  & \coefc{2e_3}[\bfz]
\end{bmatrix}
$$
is non-singular.
As a consequence, 
%\begin{center}
there are two real matrices, an orthogonal matrix $\Pmat$ and a non-singular diagonal matrix $\Dmat$, depending only on the set of coefficients $c$ evaluated at $\pG$, such that $\Cmat = \Pmat\Dmat\Pmat^T$. 
%\end{center}
\end{hyp} 
% Follow the key properties of the convected Helmholtz differential operator for any $\pG\in\mathbb R^3$ .
\noindent
The convected Helmholtz operator satisfies \Cref{hyp:PDop} according to the following comments.
 \begin{itemize}
 \item The density $\rho(\pG)$ is positive, and there is at least one index $k\in\{1,2,3\}$ such that $\bfm_k(\pG)\neq 1$ since the Mach number is assumed to be no more than $1$, $|\bfm(\pG)|<1$. 
 %Without loss of generality we will then assume that $c_{2e_1}(\pG)\neq 0$, or equivalently $ \coefc{2e_1}[\bfz] \neq 0$ .
 Hence $\rho(\pG)\left((\bfm_k)^2(\pG)-1\right)\neq 0$, in other words in particular we indeed have $c_{2e_1}\neq 0$.
 \item The matrix $\Cmat$ defined by:
$$
\begin{array}{l}
%\mathcal C:=
\begin{bmatrix}
 \coefc{2e_1}[\bfz] & \frac12 \coefc{e_1+e_2}[\bfz] & \frac12 \coefc{e_1+e_3}[\bfz] \\
 \frac12 \coefc{e_1+e_2}[\bfz] & \coefc{2e_2}[\bfz] & \frac12 \coefc{e_2+e_3}[\bfz] \\
  \frac12 \coefc{e_1+e_3}[\bfz] &\frac12 \coefc{e_2+e_3 }[\bfz] & \coefc{2e_3}[\bfz]
\end{bmatrix}
\\
=
\rho(\pG)\begin{bmatrix}
\left( \mathbf M_1(\pG)\right)^2 -1& \frac12 \mathbf M_1(\pG)\mathbf M_2(\pG)& \frac12\mathbf M_1(\pG)\mathbf M_3(\pG)\\
 \frac12\mathbf M_1(\pG)\mathbf M_2(\pG)&\left( \mathbf M_2(\pG)\right)^2-1& \frac12\mathbf M_2(\pG)\mathbf M_3(\pG)\\
  \frac12\mathbf M_1(\pG)\mathbf M_3(\pG) &\frac12 \mathbf M_2(\pG)\mathbf M_3(\pG)&\left( \mathbf M_3(\pG)\right)^2-1
\end{bmatrix}
,
\end{array}
$$
is related to the second order terms in the partial differential operator,
and  the density $\rho(\pG)$ is positive.
Under the assumption that $|\mathbf M(\pG)|<1$ it can be shown that  $1/\rho(\pG)\Cmat$ is non-singular, see appendix \ref{app:eigval}.
Hence since the Mach number is assumed to be no greater than $1$, the matrix $\Cmat$ is indeed non-singular.
 \end{itemize}

Finally, $\left\{\lambda_i,\mu_i,\nu_i, i\in(\mathbb N_0)^3\right\}$ denote complex polynomial coefficients, and, for a given integer $d$, we denote the corresponding polynomials:
$$
P:=\sum_{i\in\mathbb N_0^3, |i|\leq d} \lambda_i \mathbf X^i,\
 Q:=\sum_{i\in\mathbb N_0^3, |i|\leq d} \mu_i \mathbf X^i \text{ and }
R:=\sum_{i\in\mathbb N_0^3, |i|\leq d} \nu_i \mathbf X^i, \text{  where }\mathbf X^i = X_1^{i_1} X_2^{i_2} X_3^{i_3}.$$
As a convention, when referring to {\it a polynomial of degree at most equal to $d$} we include the zero polynomial.

%%%%%%%%%%%%%%%%%%%%%%%%%%%%%%%%%%%%%%%%%%%%%%%%%%%
\section{Three types of quasi-Trefftz functions}
\label{sec:QTfams}
%start with local and $\pG$
%RMD local process, all in the neighborhood of a single point (maybe here insist param $q$ fixed here for now)
The quasi-Trefftz property for a function is a property of the image of this function under the action of the partial differential operator. Two fundamental aspects of this property are related to its statement in terms of a Taylor expansion approximation: (i) the fact that it is a local property, and (ii) the fact that it allows for a choice in the desired order of approximation.
Hence, as it relies on enforcing the quasi-Trefftz property, the construction of quasi-Trefftz functions is performed at a given point $\pG$ and constructed functions satisfy the property at a given order of approximation $q$.

The first question is obviously that of the existence of quasi-Trefftz functions.
A second question is nevertheless equally important for the efficient implementation of quasi-Trefftz methods, that of the computational cost of the practical construction of quasi-Trefftz bases.
Indeed, the construct the quasi-Trefftz functions is only  a pre-computation to the discretization of a quasi-Trefftz weak formulation, therefore its computational cost must be acceptable compared to that of the assembly of the discrete matrix and the resolution of the linear system. 
The former question will be addressed by the derivation of an algorithm for the construction of quasi-Trefftz functions.
The latter question will be settled by the precise steps of the construction algorithm, as the algorithm only applies explicit closed formulas while it requires no numerical resolution of any system.
% and propose algorithms for their construction. WITH LIMITED COMPUTATIONAL COST since EXPLICIT FORMULA RATHER THAN SOLVING ANY SYSTEM
%{\color{red} I'd like to emphasize that the challenge here is partly to make sure that we know how to construct quasi-Trefftz functions, and partly to make sure that their computation is not a deal breaker as it will be a precomp step within the quasi-Trefftz method}

In this section, we present three types of quasi-Trefftz functions.
\begin{itemize}
\item The original GPW ansatz, namely $\exp \left( \Lambda \cdot\Big(\varx-\pG \Big) + \text{HOT}\right)$, introduced higher order terms in the phase of a PW. The general form of this ansatz can be described as $\exp P(\varx-\pG)$, for some polynomial $P$.
\item In comparison, a new ansatz was proposed \cite{AbGPW} via the introduction of higher order terms in the amplitude of a PW as $(1+\text{HOT}) \exp  \Lambda \cdot\Big(\varx-\pG \Big)$. Therefore such an ansatz has the general form $Q(\varx-\pG) \exp  \Lambda \cdot\Big(\varx-\pG \Big)$, for some polynomial $Q$ and some $\Lambda$.
\item Moreover, we propose here to consider purely polynomial quasi-Trefftz functions, described as $R(\varx-\pG)$ for some polynomial $R$. 
\end{itemize}
To guarantee the existence of such quasi-Trefftz functions, we will study the existence of polynomials ($P$, $Q$ or $R$) such that the associated ansatz satisfies the desired quasi-Trefftz property \eqref{eq:TaylExp}.
To do so we will conveniently reformulate the problem to  evidence properties of the resulting system,  underlining the shared common structure of these systems. These properties will appear to be central to the construction of quasi-Trefftz functions.

%%%%%%%%%%%%%%%%%%%%%%%%%%
\subsection{Forming a non-linear system for amplitude-based GPWs}
\label{ssec:Ab}
Initially, the abstract problem of construction of an amplitude-based  GPW can be written, for a given $q\in\mathbb N$, as:
\begin{equation}
\label{eq:AbProb}%Amplitude Based Problem
\left\{
\begin{array}{l}
\text{Find } (\polQ, \Lambda)\in\mathbb C[X_1,X_2,X_3] \times\mathbb C^3 \text{ such that }\\
J(\varx) := \polQ(\varx-\pG)\exp \Lambda \cdot\Big(\varx-\pG \Big) \text{ satisfies }  \\
\Lc J(\varx) = O(|\varx-\pG|^q)
\end{array}
\right.
\end{equation}
The unknowns here are of two types, either polynomial, for $\polQ$, or scalar, for the three components of $\Lambda$, and the specific role of $\Lambda$ will be highlighted throughout the discussion.
To formulate a more concrete problem, we focus on the action of the differential operator $\Lc$ on the ansatz. If $J(\varx) := \polQ(\varx-\pG)\exp \Lambda \cdot\Big(\varx-\pG \Big)$ with $\polQ\in\mathbb C[X_1,X_2,X_3]$ and $ \Lambda\in\mathbb C^3$, then
$$
\begin{array}{rl}
\Lc J(\varx)
=&\displaystyle
\left(
\sum_{k=1}^3
\funcc{2e_k}(\varx)\left(\partial_\varx^{2e_k}\polQ(\varx-\pG) + 2 \Lambda_k \partial_{\varx}^{e_k} \polQ(\varx-\pG) + \Lambda_k^2 \polQ(\varx-\pG)\right) \right.
\\ &\displaystyle
+
\sum_{1\leq k<k'\leq 3}
\funcc{e_k+e_k'}(\varx)\Big(\partial_{\varx}^{e_k+e_{k'}} \polQ(\varx-\pG) +\Lambda_k\partial_{\varx}^{e_{k'}} \polQ(\varx-\pG)  +\Lambda_{k'}\partial_{\varx}^{e_{k}}  \polQ(\varx-\pG) + \Lambda_k\Lambda_{k'} \polQ(\varx-\pG)\Big)
\\ &\displaystyle\left.
+
\sum_{k=1}^3
\funcc{e_k} (\varx)\Big(\partial_{\varx}^{e_k} \polQ(\varx-\pG)  +\Lambda_k \polQ(\varx-\pG) \Big)
+
%\sum_{|\alpha|=0}
\funcczr(\varx) \polQ(\varx-\pG)
\right) \exp  \Lambda \cdot\Big(\varx-\pG \Big) .
\end{array}
$$
We can then define the partial differential operator
$$
\begin{array}{rl}
\LcAm \polQ (\varx):= &
\displaystyle
\sum_{k=1}^3
\funcc{2e_k}(\varx)\left(\partial_\varx^{2e_k}\polQ(\varx-\pG) + 2 \Lambda_k \partial_{\varx}^{e_k} \polQ(\varx-\pG) + \Lambda_k^2 \polQ(\varx-\pG)\right)
\\&\displaystyle
+
\sum_{1\leq k<k'\leq 3}
\funcc{e_k+e_k'}(\varx)\Big(\partial_{\varx}^{e_k+e_{k'}} \polQ(\varx-\pG) +\Lambda_k\partial_{\varx}^{e_{k'}} \polQ(\varx-\pG)  +\Lambda_{k'}\partial_{\varx}^{e_{k}}  \polQ(\varx-\pG) + \Lambda_k\Lambda_{k'} \polQ(\varx-\pG)\Big)
\\& \displaystyle
+
\sum_{k=1}^3
\funcc{e_k} (\varx)\Big(\partial_{\varx}^{e_k} \polQ(\varx-\pG)  +\Lambda_k \polQ(\varx-\pG) \Big)
+
%\sum_{|\alpha|=0}
\funcczr(\varx) \polQ(\varx-\pG)
\end{array}
$$
to emphasize that
$$
\Lc J(\varx)
 = \Big(\LcAm \polQ (\varx)\Big) \exp  \Lambda \cdot\Big(\varx-\pG \Big),
$$
where for $ \Lambda\in\mathbb C^3$ the exponential term is locally bounded.
As a result, for $J$ to satisfy the desired property $\Lc J(\varx) = O(|\varx-\pG|^q)$, it is then sufficient for $\Lambda$ and $\polQ$ to satisfy $\LcAm \polQ(\varx)= O(|\varx-\pG|^q)$.
Therefore we will formulate a new problem for the construction of GPWs
\begin{equation}
\label{prob:Ab}
\left\{
\begin{array}{l}
\text{Find } (\polQ, \Lambda)\in\mathbb C[X_1,X_2,X_3] \times\mathbb C^3 \text{ such that }\\
%J(\varx) := \polQ(\varx)\exp \Lambda \cdot\Big(\varx-\pG \Big) \text{ satisfies }  \\
\LcAm \polQ(\varx)
= O(|\varx-\pG|^q)
 \\\text{then }J(\varx) := \polQ(\varx-\pG)\exp \Lambda \cdot\Big(\varx-\pG \Big),
\end{array}
\right.
\end{equation}
and any solution $J$ to \eqref{prob:Ab} will be solution to the initial problem \eqref{eq:AbProb}.

We can now express a concrete problem in terms of scalar equations and scalar unknowns thanks to \Cref{rmk:SimpRules},
the equations corresponding to cancelling the Taylor expansion coefficients of  $\LcAm \polQ$ for orders from 0 to $q-1$ and the unknowns corresponding to all the free parameters defining the GPW namely the $\mu$s and $\Lambda$.
Indeed, \eqref{prob:Ab} can equivalently be stated as follows as long as $d\geq q+1$: 
\begin{equation}
\label{prob:Abref}%reformulated
\left\{
\begin{array}{l}
\text{Find } \Lambda\in\mathbb C^3 \text{ and } \{\mu_{i}\in \mathbb C,i\in\mathbb N_0^3, |i|\leq d\} \text{ satisfying } \forall \varb\in\mathbb N_0^3 \text{ such that } |\varb|<q \\
\displaystyle
\sum_{\gamma\in\mathbb N_0^3;\gamma\leq \varb}\left(
\sum_{k=1}^3
\coefc{2e_k}[\varb-\gamma]\left((\gamma_k+2)(\gamma_k+1)\mu_{\gamma+2e_k} + 2 \Lambda_k (\gamma_k+1)\mu_{\gamma+e_k} + \Lambda_k^2\mu_{\gamma}\right) \right.
\\\displaystyle
+
\sum_{1\leq k<k'\leq 3}
\coefc{e_k+e_{k'}}[\varb-\gamma]\Big((\gamma_k+1)(\gamma_{k'}+1)\mu_{\gamma+e_k+e_{k'}} +\Lambda_k(\gamma_{k'}+1)\mu_{\gamma+e_{k'}}  +\Lambda_{k'}(\gamma_{k}+1)\mu_{\gamma+e_k} + \Lambda_k\Lambda_{k'}\mu_\gamma\Big)
\\ \displaystyle\left.
+
\sum_{k=1}^3
\coefc{e_k}[\varb-\gamma]\Big((\gamma_{k}+1)\mu_{\gamma+e_k}   +\Lambda_k \mu_\gamma \Big)
+
%\sum_{|\alpha|=0}
\coefczr[\varb-\gamma] \mu_\gamma\right) = 0,
\end{array}
\right.
\end{equation}
if $\{\mu_i\}$ is the set of polynomial coefficients of $\polQ$.
The choice $d\geq q+1$ simply ensures that all equations have the same structure. Indeed, for instance, there would be no $\mu_{\varb+2e_k}$ term in $\TC{\LcAm \polQ}[\beta]$ for $|\beta|=q$ if we chose $d<q+1$.
% {\color{red} box for choice of upper bound for d - it facilitates the systematic study of the system!}
 Hence we will always consider the case:
 \begin{center}
   \fbox{ $d\geq q+1$ }
\end{center}

%%%%%%%%%%%%%%%%%%%%%%%%%%%
\subsection{Forming a non-linear system for phase-based GPWs}
\label{ssec:Pb}
Similarly, the abstract problem of construction of a phase-based  GPW can initially be written, for a given $q\in\mathbb N$, as:
\begin{equation}
\label{eq:PbProb}%Phase Based Problem
\left\{
\begin{array}{l}
\text{Find a polynomial } \polP\in\mathbb C[X_1,X_2,X_3] \text{ such that }\\
G(\varx) := \exp \polP(\varx-\pG) \text{ satisfies }  \\
\Lc G(\varx) = O(|\varx-\pG|^q).
\end{array}
\right.
\end{equation}
Thanks to the definition of the partial differential operator
$$
\begin{array}{rl}
\LcPh \polP (\varx):= &
\displaystyle
\sum_{k=1 }^3
\funcc{2e_k}(\varx)\left(\partial_{\varx}^{2e_k}\polP(\varx-\pG) + \left(\partial_{\varx}^{e_k} \polP(\varx-\pG)\right)^2\right) 
 \\&\displaystyle
+
\sum_{1\leq k<k'\leq 3}
\funcc{e_k+e_k'}(\varx)\Big(\partial_{\varx}^{e_k+e_{k'}}\polP(\varx-\pG) + \partial_{\varx}^{e_k}\polP(\varx-\pG)\partial_{\varx}^{e_{k'}} \polP(\varx-\pG)\Big)
 \\&\displaystyle
+
\sum_{k=1}^3
\funcc{e_k}(\varx) \partial_{\varx}^{e_k} \polP(\varx-\pG) 
+
%\sum_{|\alpha|=0}
\funcczr(\varx),
\end{array}
$$
we can easily verify that
$$
\Lc G(\varx)
 = \Big(\LcPh \polP (\varx)\Big) \exp  \polP(\varx-\pG).
$$
Hence any solution $G$ to the problem:
\begin{equation}
\label{prob:Pb}
\left\{
\begin{array}{l}
\text{Find a polynomial } \polP\in\mathbb C[X_1,X_2,X_3] \text{ such that }\\
%G(\varx) := \exp \polP(\varx) \text{ satisfies }  \\
\LcPh\polP(\varx)
 = O(|\varx-\pG|^q)
 \\\text{then }G(\varx) := \exp \polP(\varx-\pG) 
\end{array}
\right.
\end{equation}
will also be a solution to the initial problem  \eqref{eq:PbProb}.
In terms of scalar unknowns and equations, 
%{\color{red} peut etre faire ici aussi un commentaire comme pour la sous section precedente} 
as long as $d\geq q+1$, this is then equivalent to: 
\begin{equation}
\label{prob:Pbref}%reformulated
\left\{
\begin{array}{l}
\text{Find } \{\lambda_{i}\in \mathbb C,i\in\mathbb N_0^3, |i|\leq d\} \text{ satisfying } \forall \varb \in\mathbb N_0^3\text{ such that } |\varb|<q \\
\displaystyle
\sum_{\gamma\in\mathbb N_0^3;\gamma\leq \varb}\left(
\sum_{k=1}^3
\coefc{2e_k}[\varb-\gamma]\left((\gamma_k+2)(\gamma_k+1)\lambda_{\gamma+2e_k} +  \sum_{\eta\leq \gamma}(\gamma_{k}-\eta_{k} +1)\lambda_{\gamma-\eta+e_{k}}(\eta_k +1)\lambda_{\eta+e_k}\right) \right.
\\\displaystyle
+
\sum_{1\leq k<k'\leq 3}
\coefc{e_k+e_{k'}}[\varb-\gamma]\left((\gamma_k+1)(\gamma_{k'}+1)\lambda_{\gamma+e_k+e_{k'}} + \sum_{\eta\leq \gamma}(\gamma_{k'}-\eta_{k'} +1)\lambda_{\gamma-\eta+e_{k'}}(\eta_k +1)\lambda_{\eta+e_k}\right)
\\ \displaystyle\left.
+
\sum_{k=1}^3
\coefc{e_k}[\varb-\gamma](\gamma_{k}+1)\lambda_{\gamma+e_k}    \right) 
+
%\sum_{|\alpha|=0}
\coefczr[\varb]= 0,
\end{array}
\right.
\end{equation}
if $\{\lambda_i\}$ is the set of polynomial coefficients of $\polP$.
Similarly here the choice $d\geq q+1$ simply ensures that all equations have the same structure. 

%%%%%%%%%%%%%%%%%%%%%%%%%%%
\subsection{Forming a system for polynomial functions}
\label{ssec:Pol}
The abstract problem of construction of a purely polynomial quasi-Trefftz function can simply be written, for a given $q\in\mathbb N$, as:
\begin{equation}
\label{eq:PolProb}%Polynomial Based Problem
\left\{
\begin{array}{l}
\text{Find a polynomial } \polR\in\mathbb C[X_1,X_2,X_3] \text{ such that }  \\
H(\varx) := \polR(\varx-\pG) \text{ satisfies }  \\
\Lc H(\varx) = O(|\varx-\pG|^q).
\end{array}
\right.
\end{equation}
In terms of scalar unknowns and equations, as long as $d\geq q+1$ to ensure again that all equations have the same structure, this is equivalent to: 
\begin{equation}
\label{prob:Polref}%reformulated
\left\{
\begin{array}{l}
\text{Find } \{\nu_{i}\in \mathbb C,i\in\mathbb N_0^3, |i|\leq d\} \text{ satisfying } \forall \varb\in\mathbb N_0^3 \text{ such that } |\varb|<q \\
\displaystyle
\sum_{\gamma\in\mathbb N_0^3;\gamma\leq \varb}\left(
\sum_{k=1}^3
\coefc{2e_k}[\varb-\gamma](\gamma_k+2)(\gamma_k+1)\nu_{\gamma+2e_k} \right.
\\\displaystyle
+
\sum_{1\leq k<k'\leq 3}
\coefc{e_k+e_{k'}}[\varb-\gamma](\gamma_k+1)(\gamma_{k'}+1)\nu_{\gamma+e_k+e_{k'}}
\\ \displaystyle\left.
+
\sum_{k=1}^3
\coefc{e_k}[\varb-\gamma](\gamma_{k}+1)\nu_{\gamma+e_k}   
+
%\sum_{|\alpha|=0}
\coefczr[\varb-\gamma]\nu_\gamma \right) = 0,
\end{array}
\right.
\end{equation}
if $\{\nu_i\}$ is the set of polynomial coefficients of $\polR$.

%%%%%%%%%%%%%%%%%%%%%%%%%%%
\subsection{Common structure: a hierarchy of linear subsystems}
These systems share common aspects but also exhibit differences. 
We will leverage the former to derive very similar construction algorithms for the three types of quasi-Trefftz functions. The key-point here will be to reformulate each of these systems of equations (for $|\beta|<q$) as a hierarchy of linear subsystems $\ell$ (for $|\beta|=\ell$); this hierarchical structure for increasing $\ell$ is precisely shared by the three types of systems derived respectively in Subsections \ref{ssec:Ab}, \ref{ssec:Pb} and \ref{ssec:Pol}. Hence the full systems will later be solved by solving the linear subsystems as a hierarchy for increasing values of $\ell$.

The unknowns in both systems \eqref{prob:Abref}, \eqref{prob:Pbref} and \eqref{prob:Polref} include the $(d+1)(d+2)(d+3)/6$ polynomial coefficients, respectively  $\{\mu_{i}\in \mathbb C,i\in\mathbb N_0^3, |i|\leq d\}$,$\{\lambda_{i}\in \mathbb C,i\in\mathbb N_0^3, |i|\leq d\}$, and $\{\nu_{i}\in \mathbb C,i\in\mathbb N_0^3, |i|\leq d\}$, while only in the amplitude-based case there are three additional scalar unknowns, $\Lambda\in\mathbb C^3$. In the polynomial case the system is linear, whereas in both GPW cases the systems are non-linear. However in the amplitude-based case the non-linear terms are limited to products of one $\mu_i$ and powers of $\Lambda_k$s.
%In order to emphasize similarities in the structure of the systems, we will now impose the value of $\Lambda$ to be given, so that the unknowns in both systems are now only the polynomial coefficients, $\{\lambda_{i}\in \mathbb C,i\in\mathbb N_0^3, |i|\leq d\}$ and $\{\mu_{i}\in \mathbb C,i\in\mathbb N_0^3, |i|\leq d\}$, there are  $\frac{(d+1)(d+2)(d+3)}{6}$ of them in each case. System \eqref{prob:Abref} is then linear, as $\Lambda$ isn't an unknown anymore, while System \eqref{prob:Pbref} is non-linear. 

Besides, each system has $q(q+1)(q+2)/6$ equations, and we will now describe their common layer structure. A close inspection of the equations reveals an underlying structure linked to the unknowns' and equations' multi-indices. Indeed, for any equation $\beta$, unknowns $\mu_i$ or $\lambda_i$ with $|i|\leq|\beta|+1$ may appear in non-linear terms, whereas unknowns $\mu_i$, $\lambda_i$ or $\nu_i$ with $|i|=|\beta|+2$ can only appear in linear terms. This is summarized in the following two tables.
\begin{center}
{\renewcommand{\arraystretch}{1.5}
\begin{tabular}{|c|c|c|c|}
\hline
Amplitude-based &  Phase-based& Indices&Comments\\\hline\hline
$\Lambda_{k}\mu_{\gamma+e_k}$& $\lambda_{\gamma-\eta+e_{k}}\lambda_{\eta+e_k}$& 
$\eta\leq \gamma\leq \varb$ and $1\leq k\leq 3$&$\begin{array}{c}|\gamma+e_{k }|\leq |\varb|+1\\|\eta+e_{k }|\leq|\varb|+1\\|\gamma-\eta+e_{k }|\leq|\varb|+1\end{array}$
 \\\hline
 $ \Lambda_k^2\mu_{\gamma}$ && 
$\gamma\leq \varb$ and $1\leq k\leq 3$&$|\gamma|\leq |\varb|$
 \\\hline
 $\begin{array}{c}\Lambda_{k}\mu_{\gamma+e_{k'}}\\  \Lambda_{k'}\mu_{\gamma+e_k} \end{array}$&$\lambda_{\gamma-\eta+e_{k'}}\lambda_{\eta+e_k}$&
$\eta\leq \gamma\leq \varb$ and $1\leq k<k'\leq 3$&$\begin{array}{c}|\gamma+e_{k }|\leq|\varb|+1\\|\gamma+e_{k^{'} }|\leq|\varb|+1\\|\gamma-\eta+e_{k^{'} }|\leq|\varb|+1\\|\eta+e_{k^{'} }|\leq|\varb|+1\end{array}$
 \\\hline
$\Lambda_k\Lambda_{k'}\mu_\gamma$& &
$\gamma\leq \varb$ and $1\leq k<k'\leq 3$&$|\gamma|\leq|\varb|$
 \\\hline
$\Lambda_k \mu_\gamma$& & $ \gamma\leq \varb$ and $1\leq k\leq 3$&$|\gamma|\leq|\varb|$
 \\\hline
 \end{tabular}}\\
{Non-linear terms in amplitude-base and phase-base \Cref{prob:Abref,prob:Pbref} for a given $\varb$. 
}
\end{center}
\begin{center}
{\renewcommand{\arraystretch}{1.5}
\begin{tabular}{|c|c|c|c|c|}
\hline
Ampl.-based &  Phase-based&Polynomial& Indices&Comments\\\hline\hline
 $\mu_{\gamma+2e_k} $&
  $\lambda_{\gamma+2e_k} $&
   $\nu_{\gamma+2e_k} $&
 $ \gamma\leq \varb$ and $1\leq k\leq 3$&$|\gamma+2e_k|=|\varb|+2$ if $\gamma=\varb$\\ 
 &&&& otherwise $|\gamma+2e_k|\leq|\varb|+1$
 \\\hline
 $\mu_{\gamma+e_k+e_k'}$&$\lambda_{\gamma+e_k+e_k'}$& $\nu_{\gamma+e_k+e_k'}$&
 $\gamma\leq \varb$ and $1\leq k<k'\leq 3$&$|\gamma+e_k+e_k'|=|\varb|+2$ if $\gamma=\varb$\\ 
 &&&& otherwise $|\gamma+2e_k|\leq|\varb|+1$
\\\hline
$\mu_{\gamma+e_k}$&$\lambda_{\gamma+e_k}$ & $\nu_{\gamma+e_k}$&
 $ \gamma\leq \varb$ and $1\leq k\leq 3$&$|\gamma+e_k|\leq|\varb|+1$
\\\hline
$\mu_{\gamma}$&& $\nu_{\gamma}$&
 $ \gamma\leq \varb$&$|\gamma|\leq|\varb|+1$
\\\hline
 \end{tabular}}\\
{Linear terms in amplitude-base, phase-base and polynomial \Cref{prob:Abref,prob:Pbref,prob:Polref} for a given $\varb$.}
\end{center}
%
%\textcolor{blue}{
%\begin{center}
%\begin{tabular}{|c|c|c|c|}
%\hline
% &Phase-based& Amplitude based & Indices\\\hline\hline
%non-linear terms&
% $\lambda_{\gamma-\eta+e_{k'}}\lambda_{\eta+e_k}$&
% $\Lambda_{k'}\mu_{\gamma+e_k}$&
%$\eta\leq \gamma\leq \varb$ and $1\leq k<k'\leq 3$
% \\\hline
% linear terms &
% $\lambda_{\gamma+2e_k} $&
% $\mu_{\gamma+2e_k} $&
% $ \gamma\leq \varb$ and $1\leq k\leq 3$
% \\\hline
% \end{tabular}\\
% pas mieux je pense
%\end{center}}
We have seen that in each %both the amplitude-based and the phase-based 
 case  choosing $d\geq q+1$ ensures that all equations of the system share a common structure. It is straightforward to see from these tables that none of the $\mu$, $\lambda$ or $\nu$ unknowns with indices $i\in\mathbb N_0^3$ such that $|i|>q+1$ appear in the system since $|\varb|\leq q-1$. Hence these unknowns are not constrained by the system: their values do not affect the system - and therefore neither do they affect the quasi-Trefftz property - even though they would of course affect the definition of the corresponding quasi-Trefftz function.
As a result, it is sufficient to seek a polynomial $P$, $Q$ and $R$ of degree $d$ satisfying:
 \begin{center}
   \fbox{ $d=q+1$ }
\end{center}

As we can see from the {\it Comments} columns of the previous tables, all non-linear terms involve unknowns with multi-indices of length at most equal to $|\varb|+1$, while the only unknowns with multi-indices of length $|\varb|+2$ are  $\mu_{\varb+e_k+e_k'}$, $\lambda_{\varb+e_k+e_k'}$ or $\nu_{\varb+e_k+e_k'}$ with $1\leq k\leq k'\leq 3$.
The length of multi-indices then plays an important role in the structure of the system. %: it allows for gathering all non-linear terms within some layers, and focus on the remaining (linear) terms in another layer. 
Hence, in the index space $(\mathbb N_0)^3$, we describe the set of multi-indices of a given length $\ell$ as a layer, as illustrated in \Cref{fig:IndexLayer}. %explain layer in index space
In order to take advantage of the systems layer structure, we will now split the sets of equations and unknowns according to their multi-index lengths.
\newcommand{\DDT}{\draw[dashed,thick]}
\newcommand{\NodeThreeD}[3]{\draw (#1,#2,#3) node[circle,fill,inner sep=2pt] {};
\DDT(#1,#2,0)--(#1,#2,#3);\DDT(0,#2,0)--(#1,#2,0);\DDT(#1,0,0)--(#1,#2,0);
\DDT(#1,0,#3)--(#1,#2,#3);\DDT(0,0,#3)--(#1,0,#3);\DDT(#1,0,0)--(#1,0,#3);
\DDT(0,#2,#3)--(#1,#2,#3);\DDT(0,#2,0)--(0,#2,#3);\DDT(0,0,#3)--(0,#2,#3);}
\begin{figure}[htb]\centering
\resizebox{.49\linewidth}{!}{
\begin{tikzpicture}
[rotate around x=-90,rotate around y=0,rotate around z=-30,grid/.style={very thin,gray}]
        \foreach \x in {0,1,...,7}
        \foreach \y in {0,1,...,7}
        \foreach \z in {0,1,...,7}
        {
            \draw[grid,lightgray] (\x,0,0) -- (\x,7,0);
            \draw[grid,lightgray] (0,\y,0) -- (7,\y,0);
            \draw[grid,lightgray] (0,\y,0) -- (0,\y,7);
            \draw[grid,lightgray] (0,0,\z) -- (0,7,\z);
        }
        \draw [->] (0,0) -- (8,0,0) node [right] {$\indn_{1}$};
        \draw [->] (0,0) -- (0,8,0) node [above] {$\indn_{2}$};
        \draw [->] (0,0) -- (0,0,8) node [below left] {$\indn_3$};
\draw(0,0)node[left]{$(0,0,0)$};
        \draw[fill=blue,opacity=.5] (0,7,0)--(0,0,7)--(7,0,0);
\NodeThreeD313;
        \draw[] (3,1,3) node[color=white,circle,fill] {};
\draw (3,0,0) node[below left] {$3$};
\draw (0,1,0) node[above left] {$1$};
\draw (0,0,3) node[left] {$3$};
\draw (0,0,7) node[left] {$7$};
\draw (0,7,0) node [above] {$7$};
\draw (7,0,0) node[below left] {$7$};
        \foreach \x in {0,1,...,7}
        {
        \pgfmathsetmacro{\int}{7-\x}
        \foreach \y in {0,...,\int}
        {
            \draw (\x,\y,7-\x-\y) node[circle,fill=blue,inner sep=2pt] {};
        }
        }
%        \draw[fill=gray,opacity=.5] (0,3,0)--(0,0,3)--(3,0,0);
%        \foreach \x in {0,1,...,3}
%        {
%        \pgfmathsetmacro{\int}{3-\x}
%        \foreach \y in {0,...,\int}
%        {
%            \draw (\x,\y,3-\x-\y) node[circle,fill=gray,inner sep=2pt] {};
%        }
%        }
\end{tikzpicture}
}
\caption{Illustration of a multi-index layer in  $(\mathbb N_0)^3$. The layer $\{\indn\in\mathbb N_0^3, |\indn|=\ell\}$ for $\ell =7$ is represented in blue. All elements in the layer are represented as blue dots, the element $\indn=(3,1,3)$ is highlighted in white.}
\label{fig:IndexLayer}
\end{figure}

Let's consider,  for $\ell\in\mathbb N_0$ with $\ell<q$, the subset of equations corresponding to $|\varb|=\ell$. From our previous observations we know that all the terms involving unknowns with multi-index length equal to $\ell+2$, namely $\mu_{\varb+e_k+e_k'}$, $\lambda_{\varb+e_k+e_k'}$ or $\nu_{\varb+e_k+e_k'}$ with $1\leq k\leq k'\leq 3$, are linear. 
Hence if unknowns with a shorter multi-index -- and the $\Lambda$ unknowns in the amplitude-based case -- were already known, it would suggest, for $\ell\in\mathbb N_0$ with $\ell<q$, to define a linear underdetermined subsystem with:
\begin{itemize}
\item $\frac{(\ell+1)(\ell+2)}{2}$ equations, namely the equations corresponding to $\varb$ with $|\varb|=\ell$,
\item $\frac{(\ell+3)(\ell+4)}{2}$ unknowns, namely the unknowns $\{ \mu_i,|i|=\ell+2 \}$,  $\{  \lambda_i,|i|=\ell+2 \}$ or $\{ \nu_i,|i|=\ell+2 \}$,
\item a right hand side depending on $\{ \mu_i,|i|<\ell+2 \}\cup \Lambda$, $\{  \lambda_i,|i|<\ell+2 \}$, or $\{ \nu_i,|i|<\ell+2 \}$.
\end{itemize}
To ensure that the right hand side is known at each layer $\ell$, it is then only natural to proceed layer by layer for increasing values of $\ell$ from $0$ to $q-1$.  
%KINDA ALREADY SAID EARLEIR In the index spaces, either $i\in\mathbb N_0^3$ or $\varb\in\mathbb N_0^3$, the set of indices of a given length can be seen as layer, therefore we refer to the system's structure as a layer structure.

The construction of a solution to the initial system, \eqref{prob:Abref}, \eqref{prob:Pbref} or \eqref{prob:Polref}, then boils down to the successive construction of a solution to each subsystem.
In each case, a set of subsystems gathers $\displaystyle 1/2 \sum_{\ell=0}^{q-1} (\ell+1)(\ell+2) = q(q+1)(q+2)/6$ equations, so that is exactly the full set of equations of the initial system. 
From the point of view of unknowns the situation is different. Aside from the $\displaystyle 1/2\sum_{\ell=0}^{q-1} (\ell+3)(\ell+4)= q(q^2+9q+26)/6$ unknowns appearing in the combined subsystems,
we immediately notice that  the unknowns $\{ \mu_i,|i|\leq 1 \}$, $\Lambda$, $\{  \lambda_i,|i|\leq 1 \}$, and $\{ \nu_i,|i|\leq 1 \}$ do not belong to any set of subsystem unknowns, but only appear in right hand sides of the subsystems. 
So the subsystems' solvability won't  be affected by these terms, yet their values need to be fixed in order for the subsystem's right hand sides to be known.
The construction of a solution to the initial system will hence start from setting the values of  $\{ \mu_i,|i|\leq 1 \}$ and $\Lambda$, the values of $\{  \lambda_i,|i|\leq 1 \}$ or the values of $\{ \nu_i,|i|\leq 1 \}$ before turning to the hierarchy of subsystems for increasing values of $\ell$ from $0$ to $q-1$.

Not only do the amplitude and phase based cases share the same layer structure, but their subsystems also share the same structure. Indeed, independently of the case, for a given layer $\ell$, the subsystem reads as
\begin{equation}
\label{gensubsys}%general subsystem
\left\{
\begin{array}{l}
\text{Find } \{\xi_{i}\in \mathbb C,i\in\mathbb N_0^3, |i|=\ell+2\} \text{ satisfying } \forall \varb \text{ such that } |\varb|=\ell \\
\displaystyle
\sum_{k=1}^3
(\varb_k+2)(\varb_k+1)\coefc{2e_k}[\bfz]\xi_{\varb+2e_k} 
\\\displaystyle
+
\sum_{1\leq k<k'\leq 3}
(\varb_k+1)(\varb_{k'}+1)\coefc{e_k+e_{k'}}[\bfz]\xi_{\varb+e_k+e_{k'}}  = \mathsf B_\varb,
\end{array}
\right.
\end{equation}
where the right hand side $\mathsf B$ depends not only on the case but also on variable coefficients of the PDE. 
This is because the three quantities $\LcAm \polQ(\varx)$, $\LcPh \polP (\varx)$ and $\Lc R(\varx)$ all  include terms coming from the operators
$$
\displaystyle
\sum_{k=1 }^3
\funcc{2e_k}(\pG)\partial_{\varx}^{2e_k}
+
\sum_{1\leq k<k'\leq 3}
\funcc{e_k+e_k'}(\pG)\partial_{\varx}^{e_k+e_{k'}}
=
\sum_{1\leq k\leq k'\leq 3}
\funcc{e_k+e_k'}(\pG)\partial_{\varx}^{e_k+e_{k'}}.
$$
As a consequence, the study of existence of solutions to these subsystems is independent of the case.
These subsystems will be the backbone of the construction algorithm for both families of GPWs.

\begin{rmk}
\label{rem:pol}
Gathering unknowns according to the length of their index, $|i|$, 
is related to splitting the unknowns from the polynomial $P$, $Q$ or $R$ according  to the total degree of each monomial:
\begin{equation*}
P = \sum_{\ell' = 0}^{q+1} 
\left(\sum_{ i\in\mathbb N_0^3,
 |i| = \ell'} \lambda_i \mathbf X^i \right),\
 Q = \sum_{\ell' = 0}^{q+1} 
\left(\sum_{ i\in\mathbb N_0^3,
 |i| = \ell'} \mu_i \mathbf X^i \right),
 \text{ or }
 R = \sum_{\ell' = 0}^{q+1} 
\left(\sum_{ i\in\mathbb N_0^3,
 |i| = \ell'} \nu_i \mathbf X^i \right),
\end{equation*}
and the $\ell$th subsystem is related to certain derivatives of homogeneous polynomials of degree $\ell'=\ell+2$.
\end{rmk}

The subsystems are linear and underdetermined. Their right hand sides depend on the Taylor expansion coefficients of the  set of complex-valued PDE coefficients  $c = \{\funcci,i\in\mathbb N_0^3, |i|\leq 2\}$, as well as other unknowns hopefully previously computed thanks to a recursion on the layer $\ell$.
Let's now turn to the question of existence of solutions to each subsystem.

\subsection{Echelon form of the subsystems}

To proceed and prove the existence of solutions to each subsystem, following \Cref{rem:pol}, we will denote by $\CoPo= \mathbb C[\mathbf X]$ the space of complex polynomials in three variables, $\mathbf X = (X,Y,Z)$, and by $\HoPo{d}\subset \CoPo$ the space of  homogeneous polynomials of degree $d$.
In order to prove the existence of a solution to each linear subsystem, we therefore introduce the partial differential operator
\begin{equation*}
%\label{eq:defPDop}
\begin{array}{rccc}
\Delta_{c,\ell}: &\HoPo{\ell+2}&\rightarrow &\HoPo{\ell}\\
&f &\mapsto& \Delta_c f
\end{array}
\end{equation*}
where, given  the  set of complex-valued PDE coefficients  $c = \{\funcci,i\in\mathbb N_0^3, |i|\leq 2\}$, the linear operator $\Delta_c$ is defined on $\mathbb A$ by
%\begin{equation}
%\label{eq:defPDop}
%\begin{array}{rl}
%\displaystyle
%\Delta_c  \left[\sum_{|i|=\ell+2} p_i \mathbf X^i\right] &\displaystyle
%= 
%\sum_{|\varb|=\ell}\left( \sum_{k=1}^3 (\varb_k+2)(\varb_k+1) \coefc{2e_k}[\bfz] p_{\varb+2e_k}  \right.
%\\&\displaystyle \left.
%+
% \sum_{1\leq k< k'\leq 3} (\varb_k+1)(\varb_{k'}+1) \coefc{e_k+e_{k'}}[\bfz] p_{\varb+e_k+e_{k'}} 
%  \right)\mathbf X^\varb
%\end{array}
%\end{equation}
\begin{equation*}
%\label{eq:defPDop}
\left\{
\begin{array}{l}
\begin{array}{rl}
\displaystyle
\Delta_c  \left[\sum_{i\in(\mathbb N_0)^3,|i|<2} \xi_i \mathbf X^i\right] &\displaystyle
= 
0,% \text{ and if } d\geq 2:
\end{array}
\\
\begin{array}{rl}
\displaystyle
 \text{if } d\geq 2,
\Delta_c  \left[\sum_{i\in(\mathbb N_0)^3, |i|\leq d} \xi_i \mathbf X^i\right] &\displaystyle
= 
\sum_{\varb\in\mathbb N_0^3;|\varb|\leq d-2}\left( \sum_{k=1}^3 (\varb_k+2)(\varb_k+1) \coefc{2e_k}[\bfz] \xi_{\varb+2e_k}  \right.
\\&\displaystyle \left.
+
 \sum_{1\leq k< k'\leq 3} (\varb_k+1)(\varb_{k'}+1) \coefc{e_k+e_{k'}}[\bfz] \xi_{\varb+e_k+e_{k'}} 
  \right)\mathbf X^\varb
\end{array}
\end{array}
\right.
\end{equation*}
For a given $\ell$, the existence of solutions to \eqref{gensubsys} is equivalent to the surjectivity of the operator $\Delta_{c,\ell}$. While $\dim \HoPo{\ell+2} = (\ell+3)(\ell+4)/2 $ and $\dim \HoPo{\ell}= (\ell+1)(\ell+2)/2$, unfortunately, it is not simple here to express explicitly the kernel of $\Delta_{c,\ell}$ to find its dimension, unlike for the 2D Laplacian operator. However, we can evidence the operator's full-rank by a careful choice of appropriate bases so that the corresponding matrix representation of $\Delta_{c,\ell}$ is in echelon form\footnote{In general the echelon form could be evidenced starting from any matrix representation of the operator, choosing any bases of $\HoPo{\ell+2} $ and $\HoPo{\ell}$,  and performing Gaussian elimination. Instead, here, the bases are carefully chosen so that the corresponding matrix is already in echelon form and no Gaussian elimination is needed. }.%by exhibiting the echelon structure of a particular matrix representation of $\Delta_{c,\ell}$ in appropriate bases.
To  describe such a matrix, we will use the canonical bases of $\HoPo{L}$ for $L\in\{\ell,\ell+2\}$, namely $\{ \mathbf X^i, i\in\mathbb N_0^3 \text{ and }|i|=L \}$, and we will number the columns for $L=\ell+2$ and the rows for $L=\ell$ according to the linear order $\prec$ defined for multi-indices $i\in\mathbb N_0^3$ with a given length $L$ by
 $$\forall (i,j)\in\left(\mathbb N_0^3\right)^2,i\prec j\text{ if } 
\left\{\begin{array}{l}
i_1<j_1 ; \text{ or}\\
i_1=j_1 \text{ and }i_2<j_2,
\end{array}\right. 
$$
while $|i|=L$ and $|j|=L$, so $i_3 = L-i_1-i_2$ and similarly $j_3=L-j_1-j_2$. On each row $\varb$ of the matrix, the only non-zero terms are
\begin{equation*}
{\renewcommand{\arraystretch}{1.7}
\begin{array}{|c|c|}
\hline
\text{Column index} & \text{Matrix entry}
\\\hline
\beta+2e_k \text{ with } 1\leq k \leq 3 & (\varb_k+2)(\varb_k+1) \coefc{2e_k}[\bfz] 
\\\hline
\varb+e_k+e_{k'} \text{ with } 1\leq k< k'\leq 3 &  (\varb_k+1)(\varb_{k'}+1) \coefc{e_k+e_{k'}}[\bfz] 
\\\hline
\end{array}
}\
% \textcolor{blue}{choose one or the other of these two tables}
% {\renewcommand{\arraystretch}{1.7}
%\begin{array}{|c|c|}
%\hline
%\text{Column index} & \text{Matrix entry}
%\\\hline
%\beta+2e_1 & (\varb_1+2)(\varb_1+1) \coefc{2e_1}[\bfz] 
%\\\hline
%\beta+2e_2 & (\varb_2+2)(\varb_2+1) \coefc{2e_2}[\bfz] 
%\\\hline
%\beta+2e_3 & (\varb_3+2)(\varb_3+1) \coefc{2e_3}[\bfz] 
%\\\hline
%\varb+e_1+e_{2} &  (\varb_1+1)(\varb_{2}+1) \coefc{e_1+e_{2}}[\bfz] 
%\\\hline
%\varb+e_1+e_{3} &  (\varb_1+1)(\varb_{3}+1) \coefc{e_1+e_{3}}[\bfz] 
%\\\hline
%\varb+e_2+e_{3} &  (\varb_2+1)(\varb_{3}+1) \coefc{e_2+e_{3}}[\bfz] 
%\\\hline
%\end{array}
%}
\end{equation*}
One can easily verify that
\begin{equation}
\label{eq:orderIndices}
\varb + 2e_3
\prec \varb + e_2+e_3
\prec \varb + 2e_2
\prec \varb + e_1+e_3
\prec \varb + e_1+e_2
\prec \varb + 2e_1,
\end{equation}
therefore the last nonzero entry on each row $\varb$ stands in column $\varb+2e_1$, and the entry is $ (\varb_1+2)(\varb_1+1) \coefc{2e_1}[\bfz] $. So, under the simple assumption that $ \coefc{2e_1}[\bfz] \neq 0$ (as included in \Cref{hyp:PDop}), the last nonzero entries for any two rows $\varb$ and $\widetilde \varb$ with $\varb \prec \widetilde \varb$ are respectively in columns $\varb+2e_1$ and $\widetilde \varb+2e_1$. Since $\varb+2e_1 \prec \widetilde \varb+2e_1$ it shows the echelon structure of the matrix, implying that the matrix has full rank.

\begin{rmk}
\label{rmk:nonSym}
By deciding for a numbering scheme  we broke the symmetry between the three component indices. Of course, under the assumption that $ \coefc{2e_2}[\bfz] \neq 0$, or $ \coefc{2e_3}[\bfz] \neq 0$, we could choose an index numbering scheme for which the corresponding matrix would have a similar echelon structure.
\end{rmk}

Back to Subsystem \eqref{gensubsys}, this discussion can be summarized as follows.
\begin{lmm}
Given $\ell\in\mathbb N_0$ a point $\pG\in\mathbb R^3$ and a set of complex-valued functions  $c = \{\funcci,i\in\mathbb N_0^3, |i|\leq 2\}$, with
$c_{2e_1}(\pG)\neq 0$, a matrix representation of the partial differential operator $\Delta_{c,\ell}$ is in echelon form, hence the operator is surjective.
 
Moreover, for Subsystem \eqref{gensubsys}, if both the indices $\beta$ of equations and the indices $i$ of unknowns are numbered according to the linear order $\prec$,
 this is the matrix of the Subsystem, hence it guarantees the existence of solutions for any right hand side $B\in \mathbb C^{(\ell+1)(\ell+2)/2}$.
\end{lmm}

In practice,
% unique solution - rather in order to construct a solution
the space spanned by the $\frac{(\ell+3)(\ell+4)}{2}$ columns of the matrix is of dimension $\frac{(\ell+1)(\ell+2)}{2}$. The columns displaying a $ (\varb_1+2)(\varb_1+1) \coefc{2e_1}[\bfz] $ entry form a linearly independent set, since they are in echelon form. Their indices are $\{\varb+2e_1,\varb\in(\mathbb N_0)^3,|\varb|=\ell\}$, which is equivalent to $\{i\in(\mathbb N_0)^3,|i|=\ell+2, i_1>1 \}$ as represented in \Cref{fig:layer}.
Hence the indices of the remaining columns, i.e. the columns which do not display a $ (\varb_1+2)(\varb_1+1) \coefc{2e_1}[\bfz] $ entry, are simply $\{i\in(\mathbb N_0)^3,|i|=\ell+2, i_1\leq1 \}$.

Accordingly, in ordre to compute a solution to the subsystem \eqref{gensubsys} for a given right hand side, we can take advantage of the echelon structure, by simply fixing first the values of the $\{\xi_{i}, i\in(\mathbb N_0)^3,|i|=\ell+2, i_1\leq1\}$ unknowns, and then solving by substitution the resulting square triangular system for the remaining unknowns  $\{\xi_{i}, i\in(\mathbb N_0)^3,|i|=\ell+2, i_1>1\}$.  See \Cref{Algo:SubS}. Hence, due to the particular echelon form of the matrix, Gaussian elimination is not necessary, and \Cref{Algo:SubS} is simply the back-substitution that follows it. Note that from the matrix point of view this is equivalent to turning the matrix into a square triangular matrix by adding rows of the identity.
\begin{figure}
\centering
\resizebox{.8\linewidth}{!}{
\begin{tikzpicture}
% if I need to push the three triangles away from each other
\def\delta{1}

\def\costhirty{0.866025}
\def\sinthirty{0.5}
\def\cosonetwenty{-0.5}
\def\sinonetwenty{0.866025}
\def\lev{5}%level l considered
\pgfmathsetmacro{\levt}{\lev+2}

% initial point for the North portion of the graph (North third)
\def\inix{0}
\pgfmathsetmacro{\iniy}{\delta}

% North-third axes, starting from the North-directed vector [0;\lev+2]
% "j like" axis for this portion of the graph (North third)
\pgfmathsetmacro{\iendb}{\inix+\sinthirty*\levt}
\pgfmathsetmacro{\jendb}{\iniy+\costhirty*\levt}
%\draw[thick,->] (\inix,\iniy) -- (\iendb,\jendb);
% "i like" axis for this portion of the graph (North third)
\pgfmathsetmacro{\ienda}{\inix-\sinthirty*\levt}
\pgfmathsetmacro{\jenda}{\iniy+ \costhirty*\levt}
%\draw[thick,->] (\inix,\iniy) -- (\ienda,\jenda);

%% axes labels
%\draw (\iendb,\jendb) node[anchor=south,font=\huge] {$\tilde j$};
%\draw (\ienda,\jenda) node[anchor=south,font=\huge] {$\tilde i$};

%%plot title
%\pgfmathsetmacro{\ttx}{\inix}
%\pgfmathsetmacro{\tty}{\iniy+\levt}
%\draw (\ttx,\tty) node[anchor=south,font=\huge] {$\lambda^x$-indices};

% alpha component in the vertical direction - North third
%defining the x-axis direction for the alpha component starting from direction [0;1]
\pgfmathsetmacro{\axdirx}{\sinthirty}
\pgfmathsetmacro{\axdiry}{\costhirty}
%defining the y-axis direction for the alpha component starting from direction [0;1]
\pgfmathsetmacro{\aydirx}{-\sinthirty}
\pgfmathsetmacro{\aydiry}{\costhirty}

\pgfmathsetmacro{\fina}{\lev+1}
\foreach \x in {0,...,\fina}
{ % North grid
  \pgfmathsetmacro{\xs}{\inix+\x*\axdirx}%x component of starting point
  \pgfmathsetmacro{\ys}{\iniy+\x*\axdiry}
  \pgfmathsetmacro{\xe}{\xs+(\levt-\x)*\aydirx}
  \pgfmathsetmacro{\ye}{\ys+(\levt-\x)*\aydiry}
  \draw[thick,{color=orange}] (\xs,\ys) -- (\xe,\ye);
}
\foreach \y in {0,...,\fina}
{ % North grid
  \pgfmathsetmacro{\xs}{\inix+\y*\aydirx}%x component of starting point
  \pgfmathsetmacro{\ys}{\iniy+\y*\aydiry}
  \pgfmathsetmacro{\xe}{\xs+(\levt-\y)*\axdirx}
  \pgfmathsetmacro{\ye}{\ys+(\levt-\y)*\axdiry}
  \draw[thick,dotted,{color=red}] (\xs,\ys) -- (\xe,\ye);
}
\foreach \z in {0,...,\fina}
{ % North grid
  \pgfmathsetmacro{\xs}{\inix+(\levt-\z)*\axdirx}%x component of starting point
  \pgfmathsetmacro{\ys}{\iniy+(\levt-\z)*\axdiry}
  \pgfmathsetmacro{\xe}{\inix-(\levt-\z)*\axdirx}
  \pgfmathsetmacro{\ye}{\ys}
  \draw[thick,dashed,{color=brown}] (\xs,\ys) -- (\xe,\ye);
}
  
  %highlighting the stencil
\pgfmathsetmacro{\betax}{3}
\pgfmathsetmacro{\betay}{1}
\pgfmathsetmacro{\betaz}{1}
% first the three "outer" points
\pgfmathsetmacro{\xs}{\inix +(\betax+2)*\axdirx+(\betay+0)*\aydirx}
\pgfmathsetmacro{\ys}{\iniy+ (\betax+2)*\axdiry+(\betay+0)*\aydiry}
        \draw[] (\xs,\ys)  node[color=blue,diamond,fill,inner sep=4pt] {};
        \draw[] (\xs,\ys)  node[color=black,circle,fill,inner sep=2pt] {};
\pgfmathsetmacro{\xs}{\inix +(\betax+0)*\axdirx+(\betay+2)*\aydirx}
\pgfmathsetmacro{\ys}{\iniy+ (\betax+0)*\axdiry+(\betay+2)*\aydiry}
        \draw[] (\xs,\ys)  node[color=black,circle,fill,inner sep=2pt] {};
\pgfmathsetmacro{\xs}{\inix +(\betax+0)*\axdirx+(\betay+0)*\aydirx}
\pgfmathsetmacro{\ys}{\iniy+ (\betax+0)*\axdiry+(\betay+0)*\aydiry}
        \draw[] (\xs,\ys)  node[color=black,circle,fill,inner sep=2pt] {};
% then the three "inner" points
\pgfmathsetmacro{\xs}{\inix +(\betax+1)*\axdirx+(\betay+0)*\aydirx}
\pgfmathsetmacro{\ys}{\iniy+ (\betax+1)*\axdiry+(\betay+0)*\aydiry}
        \draw[] (\xs,\ys)  node[color=black,circle,fill,inner sep=2pt] {};
\pgfmathsetmacro{\xs}{\inix +(\betax+0)*\axdirx+(\betay+1)*\aydirx}
\pgfmathsetmacro{\ys}{\iniy+ (\betax+0)*\axdiry+(\betay+1)*\aydiry}
        \draw[] (\xs,\ys)  node[color=black,circle,fill,inner sep=2pt] {};
\pgfmathsetmacro{\xs}{\inix +(\betax+1)*\axdirx+(\betay+1)*\aydirx}
\pgfmathsetmacro{\ys}{\iniy+ (\betax+1)*\axdiry+(\betay+1)*\aydiry}
        \draw[] (\xs,\ys)  node[color=black,circle,fill,inner sep=2pt] {};

%% markers
%\pgfmathsetmacro{\finaa}{\lev+2}
%\foreach \xx in {0,...,\finaa}
%{
%\pgfmathsetmacro{\finb}{\finaa-\xx}
% \foreach \yy in {0,...,\finb}
% {
%\pgfmathsetmacro{\xs}{\inix +\xx*\axdirx+\yy*\aydirx}
%\pgfmathsetmacro{\ys}{\iniy+ \xx*\axdiry+\yy*\aydiry}
%  \draw (\xs,\ys) node[circle,fill=black,inner sep=2pt] {};%\draw[mark=o,fill=black,mark size=3pt,mark options={color=black}] plot coordinates {(\xs,\ys)}; 
%  }}

%legend
  \pgfmathsetmacro{\x}{\inix+6*\axdirx}%x position of the legend  
  \pgfmathsetmacro{\xm}{\x-0.025cm}%x position of the legend  
  \pgfmathsetmacro{\xp}{\x-0.005cm}%x position of the legend  
  \pgfmathsetmacro{\ya}{\iniy+3*\axdiry}%y position for component a 
  \pgfmathsetmacro{\yb}{\iniy+2*\axdiry}%y position for component b
  \pgfmathsetmacro{\yc}{\iniy+1*\axdiry}%y position for component c
  
\draw[thick,{color=orange}] (\xm,\ya)  -- (\xp,\ya);
\draw (\x,\ya) node[anchor=west,font=\huge,{color=orange}] {$i_1$ constant};
\draw[thick,{color=red},dotted] (\xm,\yb)  -- (\xp,\yb);
\draw (\x,\yb) node[anchor=west,font=\huge,{color=red}] {$i_2$ constant};
\draw[thick,{color=brown},dashed] (\xm,\yc)  -- (\xp,\yc);
\draw (\x,\yc) node[anchor=west,font=\huge,{color=brown}] {$i_3$ constant};

% initial point for the North portion of the graph (North third)
\def\inix{9}

% North-third axes, starting from the North-directed vector [0;\lev+2]
% "j like" axis for this portion of the graph (North third)
\pgfmathsetmacro{\iendb}{\inix+\sinthirty*\levt}
\pgfmathsetmacro{\jendb}{\iniy+\costhirty*\levt}
%\draw[thick,->] (\inix,\iniy) -- (\iendb,\jendb);
% "i like" axis for this portion of the graph (North third)
\pgfmathsetmacro{\ienda}{\inix-\sinthirty*\levt}
\pgfmathsetmacro{\jenda}{\iniy+ \costhirty*\levt}
%\draw[thick,->] (\inix,\iniy) -- (\ienda,\jenda);

%% axes labels
%\draw (\iendb,\jendb) node[anchor=south,font=\huge] {$\tilde j$};
%\draw (\ienda,\jenda) node[anchor=south,font=\huge] {$\tilde i$};

%%plot title
%\pgfmathsetmacro{\ttx}{\inix}
%\pgfmathsetmacro{\tty}{\iniy+\levt}
%\draw (\ttx,\tty) node[anchor=south,font=\huge] {$\lambda^x$-indices};

% alpha component in the vertical direction - North third
%defining the x-axis direction for the alpha component starting from direction [0;1]
\pgfmathsetmacro{\axdirx}{\sinthirty}
\pgfmathsetmacro{\axdiry}{\costhirty}
%defining the y-axis direction for the alpha component starting from direction [0;1]
\pgfmathsetmacro{\aydirx}{-\sinthirty}
\pgfmathsetmacro{\aydiry}{\costhirty}

\pgfmathsetmacro{\fina}{\lev+1}
\foreach \x in {0,...,\fina}
{ % North grid
  \pgfmathsetmacro{\xs}{\inix+\x*\axdirx}%x component of starting point
  \pgfmathsetmacro{\ys}{\iniy+\x*\axdiry}
  \pgfmathsetmacro{\xe}{\xs+(\levt-\x)*\aydirx}
  \pgfmathsetmacro{\ye}{\ys+(\levt-\x)*\aydiry}
  \draw[thick,{color=orange}] (\xs,\ys) -- (\xe,\ye);
}
\foreach \y in {0,...,\fina}
{ % North grid
  \pgfmathsetmacro{\xs}{\inix+\y*\aydirx}%x component of starting point
  \pgfmathsetmacro{\ys}{\iniy+\y*\aydiry}
  \pgfmathsetmacro{\xe}{\xs+(\levt-\y)*\axdirx}
  \pgfmathsetmacro{\ye}{\ys+(\levt-\y)*\axdiry}
  \draw[thick,dotted,{color=red}] (\xs,\ys) -- (\xe,\ye);
}
\foreach \z in {0,...,\fina}
{ % North grid
  \pgfmathsetmacro{\xs}{\inix+(\levt-\z)*\axdirx}%x component of starting point
  \pgfmathsetmacro{\ys}{\iniy+(\levt-\z)*\axdiry}
  \pgfmathsetmacro{\xe}{\inix-(\levt-\z)*\axdirx}
  \pgfmathsetmacro{\ye}{\ys}
  \draw[thick,dashed,{color=brown}] (\xs,\ys) -- (\xe,\ye);
}

% markers
\pgfmathsetmacro{\finaa}{\lev+2}
\foreach \xx in {2,...,\finaa}
{
\pgfmathsetmacro{\finb}{\lev+2-\xx}
 \foreach \yy in {0,...,\finb}
 {
\pgfmathsetmacro{\xs}{\inix +\xx*\axdirx+\yy*\aydirx}
\pgfmathsetmacro{\ys}{\iniy+ \xx*\axdiry+\yy*\aydiry}
  \draw (\xs,\ys) node[circle,fill=blue,diamond,inner sep=2pt] {};%\draw[mark=o,fill=black,mark size=3pt,mark options={color=black}] plot coordinates {(\xs,\ys)}; 
  }}

%legend
  \pgfmathsetmacro{\x}{\inix+6*\axdirx}%x position of the legend  
  \pgfmathsetmacro{\xm}{\x-0.025cm}%x position of the legend  
  \pgfmathsetmacro{\xp}{\x-0.005cm}%x position of the legend  
  \pgfmathsetmacro{\ya}{\iniy+3*\axdiry}%y position for component a 
  
%\draw (\x,\ya) node[anchor=west,font=\huge,{color=white}] {$i_1$ constant};

\end{tikzpicture}
}
\caption{Two representations of a layer of indices $i$ of the unknowns of the $\ell$th subsystem for $\ell = 5$, corresponding to the blue layer in \Cref{fig:IndexLayer}. Each grid point corresponds to one index  $i\in(\mathbb N_0)^3$ with  $|i|=\ell+2$.
%Solid, dotted and dashed lines respectively join grid points with constant value of $i_1$, $i_2$ and $i_3$.
Left: Indices of the unknowns involved in equation $\varb = (3,1,1)$ are highlighted with black circles.
The index $\varb+2e_1 = (5,1,1)$ is highlighted as a blue diamond.
Right: Indices corresponding to $\varb+2e_1 $ for all $\varb\in(\mathbb N_0)^3$ such that $|\varb|=\ell$. %linearly independent columns of the subsystem's matrix are highlighted with blue diamonds under the assumption that $ \coefc{2e_1}[\bfz] \neq 0$.
%{\color{red} \`a toi de voir si \c{c}a te convient}
}
\label{fig:layer}
\end{figure}

%%%%%%%%%%%%%%%%%%%%%%%%%%%%%%%%%%%%
\subsection{Construction of 
quasi-Trefftz functions%Amplitude and Phase based GPWs
}
%Given a point $\pG\in\mathbb R^3$, a set of complex-valued functions  $c = \{\funcci,i\in\mathbb N_0^3, |i|\leq 2\}$ is assumed to satisfy \Cref{hyp:PDop}. 
%
Given a point $\pG\in\mathbb R^3$ and a set of complex-valued functions  $c = \{\funcci,i\in\mathbb N_0^3, |i|\leq 2\}$ satisfying  \Cref{hyp:PDop}, we can now turn back to the construction of solutions to Systems \eqref{prob:Abref}, \eqref{prob:Pbref} and \eqref{prob:Polref}, and hence the construction of quasi-Trefftz functions.
%Proceeding for each layer $\ell$ from $0$ to $q-1$, the right hand side of each subsystem is known from previous layers.
 \Cref{Algo:SubS} summarizes one way to compute $\xi= \{\xi_{i}\in \mathbb C;i\in(\mathbb N_0)^3, |i|=\ell+2\}$ solution of a subsytem \eqref{gensubsys} for a given right hand side $\mathsf B^\ell= \{\mathsf B^\ell_{\varb}\in \mathbb C;\varb\in(\mathbb N_0)^3, |\varb|=\ell\}$.

Remarkably, \Cref{Algo:PbS,Algo:AbS,Algo:PolS} build solutions to the non-linear problems \eqref{prob:Abref}, \eqref{prob:Pbref} and \eqref{prob:Polref} while relying exclusively on explicit closed formulas.
\begin{algorithm}[H] 
\caption{$\xi=$ solve\_subsystem$\left(\ell,\mathsf B^\ell,\left\{\coefc{e_k+e_{k'}}[\bfz], 1\leq k\leq k'\leq 3\right\}\right)$}
\label{Algo:SubS}
\begin{algorithmic}[1]
%\For{$i_1 \gets 0$ to $1$}
%\For{$i_2 \gets 0$ to $\ell+2-i_1$}
%\State $i = (i_1,i_2,\ell+2-i_1-i_2)$
%\State Fix $\xi_i$
%\EndFor
%\EndFor
\State Fix $\{\xi_{i}\in \mathbb C;i\in(\mathbb N_0)^3, |i|=\ell+2, i_1\in\{0,1\}\}$
\For{$\varb_1 \gets 0$ to $ \ell$}
\For{$\varb_2 \gets 0$ to $\ell-\varb_1$}
\State $\varb := (\varb_1,\varb_2,\ell-\varb_1-\varb_2)$
\State
$\displaystyle
\begin{array}{rr}
 \xi_{\varb+2e_1} := &\displaystyle
 \frac{1}{(\varb_1+2)(\varb_1+1)\coefc{2e_1}[\bfz]}\Bigg(\mathsf {B^\ell_{\varb}}-\sum_{k=2}^3
(\varb_k+2)(\varb_k+1)\coefc{2e_k}[\bfz]\xi_{\varb+2e_k}
 \\&\displaystyle
-\sum_{1\leq k<k'\leq 3}
(\varb_k+1)(\varb_{k'}+1)\coefc{e_k+e_{k'}}[\bfz]\xi_{\varb+e_k+e_{k'}}
\Bigg)
\end{array}
$ 
\EndFor
\EndFor
\end{algorithmic}
\end{algorithm}

Thanks to \Cref{Algo:SubS}, \Cref{Algo:AbS,Algo:PbS} will compute a solution to systems \eqref{prob:Abref} and \eqref{prob:Pbref} and then construct the associated GPWs $J$ and $G$, solutions to the initial problems \eqref{eq:AbProb} and \eqref{eq:PbProb}, while \Cref{Algo:PolS} will compute a solution to system \eqref{prob:Polref} and construct the associated polynomial quasi-Trefftz function $H$. 
\begin{algorithm}[H] 
\caption{Amplitude based}
\label{Algo:AbS}
\begin{algorithmic}[1]
\State Given $\pG\in\mathbb R^3$, $q\in\mathbb N$ and $c = \{\funcci,i\in\mathbb N_0^3, |i|\leq 2\}$ satisfying  \Cref{hyp:PDop} 
\State Fix $\{\mu_{i}\in \mathbb C;i\in(\mathbb N_0)^3, |i|\in\{0,1\}\}$ and $\Lambda\in\mathbb C^3$
\For{$\ell \gets 0$ to $q-1$}
\State Compute $\mathsf B^\ell= \{\mathsf B^\ell_{\varb}\in \mathbb C;\varb\in(\mathbb N_0)^3, |\varb|=\ell\}$ according to
$$
\begin{array}{rl}
\mathsf B^\ell_{\varb}=&
\displaystyle
-
%{\color{red} \mathbf 1_{ \beta\neq \mathbf 0}}
\sum_{\gamma\in\mathbb N_0^3;\gamma< \varb}
\sum_{k=1}^3
\coefc{2e_k}[\varb-\gamma](\gamma_k+2)(\gamma_k+1)\mu_{\gamma+2e_k}\\
&\displaystyle
-\sum_{\gamma\in\mathbb N_0^3;\gamma\leq \varb}\left(
\sum_{k=1}^3
\coefc{2e_k}[\varb-\gamma]\left(%(\gamma_k+2)(\gamma_k+1)\mu_{\gamma+2e_k} + 
2 \Lambda_k (\gamma_k+1)\mu_{\gamma+e_k} + \Lambda_k^2\mu_{\gamma}\right) \right.
\\&\displaystyle \phantom{+\sum_{\gamma\leq \varb}\Bigg(}
+
\sum_{1\leq k<k'\leq 3}
\coefc{e_k+e_{k'}}[\varb-\gamma]\Big((\gamma_k+1)(\gamma_{k'}+1)\mu_{\gamma+e_k+e_{k'}} +\Lambda_k(\gamma_{k'}+1)\mu_{\gamma+e_{k'}} 
\\&\displaystyle \phantom{+\sum_{\gamma\leq \varb}\Bigg(
+
\sum_{1\leq k<k'\leq 3}
\coefc{e_k+e_{k'}}[\varb-\gamma]\Big(}
 +\Lambda_{k'}(\gamma_{k}+1)\mu_{\gamma+e_k} + \Lambda_k\Lambda_{k'}\mu_\gamma\Big)
\\ &\displaystyle \phantom{+\sum_{\gamma\leq \varb}\Bigg(}\left.
+
\sum_{k=1}^3
\coefc{e_k}[\varb-\gamma]\Big((\gamma_{k}+1)\mu_{\gamma+e_k}   +\Lambda_k \mu_\gamma \Big)
+
%\sum_{|\alpha|=0}
\coefczr[\varb-\gamma] \mu_\gamma\right) 
\end{array}
$$
\State Compute $\mu^\ell= \{\mu_{i}\in \mathbb C;i\in(\mathbb N_0)^3, |i|=\ell+2\}$ via \Cref{Algo:SubS}: $$\mu^\ell = \text{ solve\_subsystem}\left(\ell,\mathsf B^\ell,\left\{\coefc{e_k+e_{k'}}[\bfz], 1\leq k\leq k'\leq 3\right\}\right)$$
\EndFor
\State $\displaystyle Q:=\sum_{i\in\mathbb N_0^3, |i|\leq q+1} \mu_i \mathbf X^i $ with $\mathbf X^i = X_1^{i_1} X_2^{i_2} X_3^{i_3}$
\State $J(\varx) := \polQ(\varx-\pG)\exp \Lambda \cdot\Big(\varx-\pG \Big)$
\end{algorithmic}
\end{algorithm}
%%%%%%%%
\begin{algorithm}[H] 
\caption{Phase based}
\label{Algo:PbS}
\begin{algorithmic}[1]
\State Given $\pG\in\mathbb R^3$, $q\in\mathbb N$ and $c = \{\funcci,i\in\mathbb N_0^3, |i|\leq 2\}$ satisfying  \Cref{hyp:PDop} 
\State Fix $\{\lambda_{i}\in \mathbb C;i\in(\mathbb N_0)^3, |i|\in\{0,1\}\}$
\For{$\ell \gets 0$ to $q-1$}
\State Compute $\mathsf B^\ell= \{\mathsf B^\ell_{\varb}\in \mathbb C;\varb\in(\mathbb N_0)^3, |\varb|=\ell\}$ according to
$$
\begin{array}{rl}
\mathsf B^\ell_{\varb}=&
\displaystyle
-
%{\color{red} \mathbf 1_{ \beta\neq \mathbf 0}}
\sum_{\gamma\in\mathbb N_0^3;\gamma< \varb}
\sum_{k=1}^3
\coefc{2e_k}[\varb-\gamma](\gamma_k+2)(\gamma_k+1)\lambda_{\gamma+2e_k} \\&
\displaystyle
-\sum_{\gamma\in\mathbb N_0^3;\gamma\leq \varb}\left(
\sum_{k=1}^3
\coefc{2e_k}[\varb-\gamma]\left( \sum_{\eta\leq \gamma}(\gamma_{k}-\eta_{k} +1)\lambda_{\gamma-\eta+e_{k}}(\eta_k +1)\lambda_{\eta+e_k}\right) \right.
\\&\displaystyle \phantom{+\sum_{\gamma\leq \varb}\Bigg(}
+
\sum_{1\leq k<k'\leq 3}
\coefc{e_k+e_{k'}}[\varb-\gamma]\Bigg((\gamma_k+1)(\gamma_{k'}+1)\lambda_{\gamma+e_k+e_{k'}} 
\\&\displaystyle \phantom{+\sum_{\gamma\leq \varb}\Bigg(+\sum_{1\leq k<k'\leq 3}\coefc{e_k+e_{k'}}[\varb-\gamma]\Bigg(}
+ \sum_{\eta\in\mathbb N_0^3;\eta\leq \gamma}(\gamma_{k'}-\eta_{k'} +1)\lambda_{\gamma-\eta+e_{k'}}(\eta_k +1)\lambda_{\eta+e_k}\Bigg)
\\& \displaystyle \phantom{+\sum_{\gamma\leq \varb}\Bigg(}\left.
+
\sum_{k=1}^3
\coefc{e_k}[\varb-\gamma](\gamma_{k}+1)\lambda_{\gamma+e_k}   
 \right) 
 -
%\sum_{|\alpha|=0}
\coefczr[\varb]
\end{array}
$$
\State Compute $\lambda^\ell= \{\lambda_{i}\in \mathbb C;i\in(\mathbb N_0)^3, |i|=\ell+2\}$ via \Cref{Algo:SubS}: $$\lambda^\ell = \text{ solve\_subsystem}\left(\ell,\mathsf B^\ell,\left\{\coefc{e_k+e_{k'}}[\bfz], 1\leq k\leq k'\leq 3\right\}\right)$$
\EndFor
\State $\displaystyle P:=\sum_{i\in\mathbb N_0^3, |i|\leq q+1} \lambda_i \mathbf X^i $ with $\mathbf X^i = X_1^{i_1} X_2^{i_2} X_3^{i_3}$
\State $G(\varx) :=\exp  \polP(\varx-\pG)$
\end{algorithmic}
\end{algorithm}
%%%%%%
\begin{algorithm}[H] 
\caption{Polynomial}
\label{Algo:PolS}
\begin{algorithmic}[1]
\State Given $\pG\in\mathbb R^3$, $q\in\mathbb N$ and $c = \{\funcci,i\in\mathbb N_0^3, |i|\leq 2\}$ satisfying  \Cref{hyp:PDop} 
\State Fix $\{\nu_{i}\in \mathbb C;i\in(\mathbb N_0)^3, |i|\in\{0,1\}\}$
\For{$\ell \gets 0$ to $q-1$}
\State Compute $\mathsf B^\ell= \{\mathsf B^\ell_{\varb}\in \mathbb C;\varb\in(\mathbb N_0)^3, |\varb|=\ell\}$ according to
$$
\begin{array}{rl}
\mathsf B^\ell_{\varb}=&
\displaystyle
-
%{\color{red} \mathbf 1_{ \beta\neq \mathbf 0}}
\sum_{\gamma\in\mathbb N_0^3;\gamma< \varb}
\sum_{k=1}^3
\coefc{2e_k}[\varb-\gamma](\gamma_k+2)(\gamma_k+1)\nu_{\gamma+2e_k}\\
&\displaystyle
-\sum_{\gamma\in\mathbb N_0^3;\gamma\leq \varb}\left(
% \phantom{+\sum_{\gamma\leq \varb}\Bigg(}+
\sum_{1\leq k<k'\leq 3}
\coefc{e_k+e_{k'}}[\varb-\gamma](\gamma_k+1)(\gamma_{k'}+1)\nu_{\gamma+e_k+e_{k'}} \right.
\\ &\displaystyle \phantom{+\sum_{\gamma\leq \varb}\Bigg(}\left.
+
\sum_{k=1}^3
\coefc{e_k}[\varb-\gamma](\gamma_{k}+1)\nu_{\gamma+e_k} 
+
%\sum_{|\alpha|=0}
\coefczr[\varb-\gamma] \nu_\gamma\right) 
\end{array}
$$
\State Compute $\nu^\ell= \{\nu_{i}\in \mathbb C;i\in(\mathbb N_0)^3, |i|=\ell+2\}$ via \Cref{Algo:SubS}: $$\nu^\ell = \text{ solve\_subsystem}\left(\ell,\mathsf B^\ell,\left\{\coefc{e_k+e_{k'}}[\bfz], 1\leq k\leq k'\leq 3\right\}\right)$$
\EndFor
\State $\displaystyle R:=\sum_{i\in\mathbb N_0^3, |i|\leq q+1} \nu_i \mathbf X^i $ with $\mathbf X^i = X_1^{i_1} X_2^{i_2} X_3^{i_3}$
\State $H(\varx) :=\polR(\varx-\pG)$
\end{algorithmic}
\end{algorithm}
%%%%%%%%

Interestingly, the quasi-Trefftz property of the quasi-Trefftz functions, $J$, $G$ and $H$, built from \Cref{Algo:PbS,Algo:AbS,Algo:PolS} are satisfied independently of  the fixed values throughout these algorithms: $$\Lc J(\varx) = O(|\varx-\pG|^q),\ \Lc G(\varx) = O(|\varx-\pG|^q) \text{ and }\Lc H(\varx) = O(|\varx-\pG|^q).$$
We will refer to the choice of these values as the initialization process. 
However an appropriate choice of initialization will be crucial to prove approximation properties of the resulting sets of quasi-Trefftz functions.
\begin{rmk}
While the governing PDE considered in the previous discussion had a zero right-hand side, the common structure of the systems formed for the three types of quasi-Trefftz functions 
For a PDE with a non-zero smooth right-hand-side,
both in the linear polynomial case and the non-linear GPW case, the quasi-Trefftz property 
$$
\Lc \varphi(\varx) -f= O(|\varx-\pG|^q)
\text{ instead of }
\Lc \varphi(\varx) = O(|\varx-\pG|^q)$$
\end{rmk}

%%%%%%%%%%%%%%%%%%%%%%%%%%%%%%%%%%%
\section{
{ Initialization process and quasi-Trefftz spaces}
}
\label{sec:norm}
Beyond the construction of individual quasi-Trefftz functions, we now turn to the construction of quasi-Trefftz spaces. Given the construction Algorithms from the developed section, it is natural to leverage the initialization process to do so. The values to be chosen for each kind of quasi-Trefftz function can be listed as follows.
\begin{center}
{\renewcommand{\arraystretch}{1.5}
\begin{tabular}{|c|c|c|}
\hline
Amplitude-based &  Phase-based& Polynomial\\\hline\hline
$\Lambda\in\mathbb C^3$ & &
\\\hline
$\mu_{\mathbf 0}$&$\lambda_{\mathbf 0}$ & $\nu_{\mathbf 0}$
\\\hline
$\left\{\mu_{i}  \in\mathbb C;i\in(\mathbb N_0)^3, |i|=1\right\}$ & $\left\{\lambda_{i} \in\mathbb C;i\in(\mathbb N_0)^3, |i|=1\right\}$& $\left\{\nu_{i} \in\mathbb C;i\in(\mathbb N_0)^3, |i|=1\right\}$ \\\hline
   for $\ell$ from $0$ to $q-1$
   &
   for $\ell$ from $0$ to $q-1$
   &
   for $\ell$ from $0$ to $q-1$
   \\
   $\{\mu_{i}\in \mathbb C;i\in(\mathbb N_0)^3,\qquad$
  &
  $\{\lambda_{i}\in \mathbb C;i\in(\mathbb N_0)^3,\qquad$% |i|=\ell+2, i_1\in\{0,1\}\}$
  &
  $\{\nu_{i}\in \mathbb C;i\in(\mathbb N_0)^3,\qquad$
   \\
   $ |i|=\ell+2, i_1\in\{0,1\}\}$&$ |i|=\ell+2, i_1\in\{0,1\}\}$& $ |i|=\ell+2, i_1\in\{0,1\}\}$\\
   \hline
 \end{tabular}}\\
 {
% Initialization summary %of choice of free parameters/initialization 
% for the three types of quasi-Trefftz functions.
 %
 Top rows: see Step 1 in \Cref{Algo:PbS,Algo:AbS,Algo:PolS}.
 
 Bottom row: see Step 1 in \Cref{Algo:SubS}. }
%{Lists of free parameter for the construction of GPWs.}
\end{center}

Once  the initialization procedure chosen, the corresponding set of quasi-Trefftz functions will span spaces, and it is the approximation property of  these spaces that will be studied later.
One significant difference between the different types of ansatz (with the same fixed maximal degree for the polynomial term) lies in a fundamental property of the corresponding space of all functions satisfying the quasi-Trefftz property: while in the GPW cases this space is infinite dimensional, in the polynomial case it has a finite dimension.

\subsection{GPW spaces}
Keeping in mind the motivation for the design of GPWs, that is adding higher order terms either in the phase or the amplitude of a PW: 
%{\color{red} maybe not clear enough, maybe add notation for the polynomial $P$ and $Q$ as a reminder, maybe also write $i\kappa\mathbf k\cdot \mathbf x$ ou quelque chose comme ca}
\begin{equation}
\label{eq:GPWint}
(1+\text{HOT})\exp \Lambda \cdot\Big(\varx-\pG \Big)
\text{ or }
\exp \Big[\Lambda \cdot\Big(\varx-\pG \Big)+\text{HOT}\Big],
\end{equation}
we now turn to the GPW initialization process.
From \Cref{Algo:SubS,Algo:PbS,Algo:AbS}, the free parameters in the construction of a GPW are $\left\{\lambda_{i} \in\mathbb C \text{ for }i\in(\mathbb N_0)^3, |i|\leq q+1,i_1\in\{0,1\}\right\}$ for a Phase-based GPW, and for an Amplitude-based GPW $\left\{\mu_{i}  \in\mathbb C\text{ for }i\in(\mathbb N_0)^3, |i|\leq q+1,i_1\in\{0,1\}\right\}$ plus $\Lambda\in\mathbb C^3$. In both cases, we follow the intuition that lead to the choice of ansatz \eqref{eq:GPWint} as a generalization of PW functions to build a family of GPWs. In order to do so, only a few free parameters are sufficient, corresponding to the linear terms in the phase, and except for the constant coefficient of the amplitude for an Amplitude-based GPW, we will naturally set the remaining parameters to zero to reduce the amount of computation associated with the construction of each GPW. 
The next table summarizes the situation. 
\begin{center}
{\renewcommand{\arraystretch}{1.5}
\begin{tabular}{|c|c|c|}
\hline
Amplitude-based &  Phase-based& Comment\\\hline\hline
$\Lambda\in\mathbb C^3$ & $\left\{\lambda_{i} \in\mathbb C;i\in(\mathbb N_0)^3, |i|=1\right\}$& $[\![1]\!]$ \\\hline
$\mu_{\mathbf 0}$& &  Set to $1$
\\\hline
 &$\lambda_{\mathbf 0}$&Set to $0$
\\\hline
  %$\left\{\lambda_{i} \in\mathbb C \text{ for }i\in(\mathbb N_0)^3, |i|\in\{0,1\}\right\}$
 $\left\{\mu_{i}  \in\mathbb C;i\in(\mathbb N_0)^3, |i|=1\right\}$& & Set to $0$
 \\\hline\hline
   for $\ell$ from $0$ to $q-1$
   &
   for $\ell$ from $0$ to $q-1$
   &
   \\
   $\{\mu_{i}\in \mathbb C;i\in(\mathbb N_0)^3,\qquad$
  &
  $\{\lambda_{i}\in \mathbb C;i\in(\mathbb N_0)^3,\qquad$% |i|=\ell+2, i_1\in\{0,1\}\}$
   & Set to $0$\\
   $ |i|=\ell+2, i_1\in\{0,1\}\}$&$ |i|=\ell+2, i_1\in\{0,1\}\}$& \\
   \hline
 \end{tabular}}\\
 {Initialization summary %of choice of free parameters/initialization 
 for the two types of GPW quasi-Trefftz functions.
 
 Top rows: corresponding to Step 1 in \Cref{Algo:PbS,Algo:AbS}.
 
 Bottom row: corresponding to Step 1 in \Cref{Algo:SubS}. }
%{Lists of free parameter for the construction of GPWs.}
\end{center}
In order to completely define our choices of quasi-Trefftz functions, it is then sufficient to describe how are chosen the three parameters corresponding to linear terms of the phase $[\![1]\!]$.

$[\![1]\!]$ In order to build not a single but rather a set of quasi-Trefftz functions, we now have three non-zero free parameters in each case, namely: %indicated in red in the first line of the previous table
\begin{equation}
\label{eq:3params}
\text{either }\Lambda=
\begin{bmatrix}\Lambda_1\\\Lambda_2\\\Lambda_3\end{bmatrix}\text{ in }\mathbb C^3
\text{ or } 
\begin{bmatrix}\lambda_{e_1}\\\lambda_{e_2}\\\lambda_{e_3}\end{bmatrix} %,
%\begin{bmatrix}\nu_{e_1}\\\nu_{e_2}\\\nu_{e_3}\end{bmatrix} 
\text{ in }\mathbb C^3.
\end{equation} 
In the constant-coefficient Helmholtz case, the matrix $\Cmat$ introduced in \Cref{hyp:PDop} is the identity  $I_3$ and it is then natural to fix these as $\cst\dir$, with $\cst = \mathrm i\kappa$ to obtain a PW exact solution since $(\cst\dir)^T (I_3\cst\dir)=-\kappa^2$ is independent of $\dir$. Yet in the general case, this matrix $\Cmat$ is associated with anisotropy in the second order terms of the partial differential operator, and it is then natural to introduce (i) the orthonormal basis of eigenvectors of $\Cmat$ via $\Pmat$ and (ii) the anisotropic scaling by the eigenvalues of $\Cmat$ via $\Dmat$.
Hence for each quasi-Trefftz function, under \Cref{hyp:PDop}, we will fix these as $\cst \Pmat\Dmat^{-1/2}\dir$, where $\cst\in\mathbb C^*$ and $\dir\in\mathbb R^3$ with $|\dir|=1$.
 To define a set of $p$ distinct -- and linearly independent under appropriate assumptions as we will see later -- quasi-Trefftz functions, we will choose distinct directions $\{ \dir_l\in\mathbb R^3 \text{ for } l \text{ from } 1 \text{ to } p; |\dir_l|=1 \text{ with } \dir_L\neq\dir_k \text{ if }k\neq l \}$ while we will choose a common value for $\cst$ for each of the $p$ functions in the set. %each function will be defined by a unique direction $\dir$. 
 Each direction $ \mathbf d_l$ will be defined by two angles $(\theta_l,\varphi_l)$ as follows:
 $$
 \mathbf d_l  = 
 \begin{bmatrix}
 \sin\theta_l\sin\varphi_l\\\sin\theta_l\cos\varphi_l\\\cos\theta_l
 \end{bmatrix}
 $$
with $\theta_l \in [0, \pi], \varphi_l\in[0,2\pi($. Hence the direction of propagation of a GPW is parametrized by the two (spherical) angles $(\theta_l,\varphi_l)$.
Under \Cref{hyp:PDop}, for any $X=\cst \Pmat\Dmat^{-1/2}\dir$, we can easily verify that:
 %could mention that it guarantees that  for any $X=[X_1,X_2,X_3]^T$ in \eqref{eq:3params} then
 $$
\begin{array}{rl} X^T( \Cmat X)& = \cst^2 (\Pmat\Dmat^{-1/2}\dir)^T(\Pmat\Dmat\Pmat^T)(\Pmat\Dmat^{-1/2}\dir)\\
& = \cst^2\dir^T\dir \\
&=\cst^2,
\end{array}
%( \text{ except for }\cst\text{ everything is real in here})
 $$
 or equivalently:
 $$
 \sum_{k=1}^3
\coefc{2e_k}[\mathbf 0]\left( X_{{k}}\right)^2 
+
\sum_{1\leq k<k'\leq 3}
\coefc{e_k+e_{k'}}[\mathbf 0] X_{{k'}}X_{k}
=
\cst^2.
$$
 As a result, this quantity does not depend on the direction $\dir$, but rather has the same value for the whole set of functions, as we discussed in the Helmholtz case. This crucial fact will be key to prove approximation properties of the quasi-Trefftz functions.
 
\begin{rmk}
As a result of this choice, the GPW functions boil  down to classical PWs in the case of a constant coefficient Helmholtz equation.
\end{rmk}

Following this remark, it is interesting to notice that the space of all either amplitude-based or phase-based GPW functions satisfying the quasi-Trefftz property at a given order is infinite dimensional, just like the space of PWs with a given wavenumber and any direction of propagation is infinite dimensional.
Moreover the discrete space defined here for a fixed value of $p$ depends of the choice of angles $\theta_l$ and $\varphi_l$, . 
This is indeed similar to the case of PW functions.

\subsection{Polynomial space}
%The case of polynomial quasi-Trefftz spaces differs from the GPW case in that the polynomial space is uniquely defined, independently of the initialization. 
On the contrary, as we will now see, the dimension of the polynomial quasi-Trefftz space is finite and equal to the number of values fixed in the initialization process.

\begin{lmm}\label{lmm:PolSpace}
Given $\pG\in\mathbb R^3$, $q\in\mathbb N$ and any differential operator $\Lc$ defined by its coefficients $c = \{\funcci,i\in(\mathbb N_0)^3, |i|\leq 2\}$ satisfying  \Cref{hyp:PDop}, the corresponding polynomial quasi-Trefftz space is the space of polynomials $H$ of degree at most equal to $q+1$ satisfying the quasi-Trefftz property $\Lc H(\varx) = O(|\varx-\pG|^q)$. This is a space of dimension $(q+2)^2$.
\end{lmm}
\begin{proof}
Given $\pG\in\mathbb R^3$, $q\in\mathbb N$ and any differential operator $\Lc$ defined by its coefficients $c = \{\funcci,i\in(\mathbb N_0)^3, |i|\leq 2\}$, the corresponding polynomial quasi-Trefftz space is the kernel of the following linear operator:
$$
\begin{array}{rrcl}
\mathrm L_q:&\mathbb P_{q+1}&\to&\mathbb C^{q(q+1)(q+2)/6}\\
&H&\mapsto&\{T_{\Lc H} [\beta], \beta\in(\mathbb N_0)^3, |\beta|\leq q-1 \},
\end{array}
$$
where $\mathbb P_{q+1}$ denotes the space of polynomials in three variables of degree at most equal to $q+1$.
In order to find the dimension of the kernel, we will consider a convenient $q(q+1)(q+2)/6\times (q+2)(q+3)(q+4)/6$ matrix of the operator $\mathrm L_q$. We choose the canonical basis of $\mathbb C^{q(q+1)(q+2)/6}$, the basis $\{  p_i,  i\in(\mathbb N_0)^3, |i|\leq q+1\}$ of $\mathbb P_{q+1}$ defined by
$$
\forall i\in(\mathbb N_0)^3, |i|\leq q+1, p_i: = (\mathbf X-\pG)^i ,
$$
and the numbering  $\mathcal N$ of multi-indices introduced in \Cref{ssec:prelim},
so the corresponding matrix $\mathsf R$ of  $\mathrm L_q$ has the following properties:
$$
\forall i\in\left(\mathbb N_0\right)^3, |i|\leq q-1, 
\left\{
\begin{array}{l}
\mathsf R_{\mathcal N(i)\mathcal N\left(i+2e_1\right)} = c_{2e_1}(\pG)% \text{ so } \mathsf R_{\mathcal N(i)\mathcal N\left(i+2e_1\right)} \neq 0 \text{ by \Cref{hyp:PDop}}
,\\
\mathsf R_{\mathcal N(i)\mathcal N(j)}=0 \text{ if } j>i+2e_1,
\end{array}
\right.
$$
If the differential operator $\Lc$ satisfies  \Cref{hyp:PDop},
 then the rows of this matrix are clearly linearly independent, since $j>i+2e_1$ implies $\mathcal N(j)>\mathcal N(i+2e_1)$.
 Hence the matrix $\mathsf R$ is full-rank, and by the rank theorem its kernel has dimension equal to:
 $$
 \frac{(q+2)(q+3)(q+4)}{6} -\frac{q(q+1)(q+2)}{6} = (q+2)^2.
 $$
\end{proof}

Since there are exactly $(q+2)^2$ free parameters to choose in the initialization process to construct a quasi-Trefftz function, a natural basis of the polynomial quasi-Trefftz space could be constructed by choosing, for each function, one and only one value to be one while all the others are set to zero. To summarize, each function in this basis, indexed by any $ j\in\left(\mathbb N_0\right)^3$, with $|j|\leq q+1$ and $j_1\in\{0,1\}$, is defined thanks to the following choice of initialization:
$$
\forall i\in\left(\mathbb N_0\right)^3, |i|\leq q+1, i_1\in\{0,1\}, \nu_i = \delta (i-j).
$$

\section{Approximation properties}
\label{sec:Approx}
%%%%%%%%%%%%%%%%%%%%%%%%%%%%%%%%%%%
The construction of quasi-Trefftz functions is based on Taylor expansions, it is therefore natural to use similar tools to study their approximation properties. The central idea here is precisely to approximate a given exact PDE solution $u$ by a linear combination $u_a$ of quasi-Trefftz functions by matching their respective Taylor expansions at $\pG$. Indeed, for any order of approximation $n$, we have:
\begin{equation}
\label{eq:TEapp}
\left\{\begin{array}{l}
\forall i\in\left(\mathbb N_0\right)^3, |i|\leq n,\\
\partial_i u (\pG) = \partial_i u_a(\pG),
\end{array}\right.
\Rightarrow u_a(x)-u(x) = O\left(\left|\varx-\pG\right|^{n+1}\right).
%{\color{red} verif:n or n+1?}
\end{equation}
This in turn leads to the convergence of $u-u_a$ in various norms of interest in the regime $\left|\varx-\pG\right|\rightarrow 0$, moreover higher order convergence follows from increasing the value of the order of approximation $n$ in the Taylor expansion.

%%%%%%%%%%%%%%%%%%%%
%describe the resulting linear system
Matching the Taylor expansion of a linear combination $u_a$ of quasi-Trefftz functions to that of a given function $u$ leads to a linear system.
\begin{itemize}
\item Each unknown is a weight of the desired linear combination, and is indexed by $l$;
 there are as many unknowns as there are functions in the quasi-Trefftz set.
\item Each equation corresponds to one Taylor expansion coefficient, and is indexed by $i\in\left(\mathbb N_0\right)^3$ with $ |i|\leq n$;
there are $(n+1)(n+2)(n+3)/6$ equations.
\end{itemize}
%%%%%%%%%%%%%%%%%%%%
The system's matrix can then be defined for a given list of quasi-trefftz function thanks to a numbering of the equations. 
The entries of the system's matrix are the  partial derivatives of quasi-Trefftz functions evaluated at $\pG$.
Hence, given any numbering $\mathcal N$ of multi-indices, for any family of $p$ functions $\{b_l,\text{ for } l\in\mathbb N, l\leq p\}$, the $(\mathcal N(i),l)$ entry of the corresponding $\frac{(n+1)(n+2)(n+3)}{6} \times p$ matrix $M^{n,p}$ is:
\begin{equation}
\label{eq:genmat}
M^{n,p}_{\mathcal N(i),l} = 
\frac{\partial_x^i b_l(\pG)}{i!} = T_{b_l}[i]
\end{equation}
Moreover, the value $p=(n+1)^2$ will be of particular interest in what follows. So to simplify the notation, the matrix corresponding to $p=(n+1)^2$ will be denoted with the superscript $[n]$, for instance $M^{[n]}$ instead of $M^{n,(n+1)^2}$.
%$M^{n,p}\in\mathbb C^{\frac{(n+1)(n+2)(n+3)}{6} \times p}$
 % {\color{red} faudrait surtout une notation pour cette matrice qui inclut $p$ et $n$  pour que tout ca reste clair}
 For three families of quasi-Trefftz functions introduced in the previous section, we will use the following notation for the corresponding matrices.
\begin{center}
{\renewcommand{\arraystretch}{1.5}
\begin{tabular}{|c|c|c|}
\hline
Amplitude-based functions &  Phase-based functions & Polynomial functions\\\hline
$\mata^{n,p}$ and $\mata^{[n]}$&
$\matp^{n,p}$ and $\matp^{[n]}$&
$\matq^{n,p}$ and $\matq^{[n]}$\\
   \hline
 \end{tabular}}
\end{center}
In order to evidence the structure of the linear system, the equations will be numbered as follows. 
%$\mathcal N(i)=$ 
To leverage the loop structure of \Cref{Algo:PbS,Algo:AbS,Algo:PolS}, we will use the numbering $\mathcal N$ introduced in \Cref{ssec:prelim}.
%we want $\mathcal N$  to satisfy $|i|<|j|$ implies that $\mathcal N(i)<\mathcal N(j)$, hence we write it as:
%$$
%\forall i\in\left(\mathbb N_0\right)^3\text{ with } |i|\leq n,
%\mathcal N(i) = \frac{|i|(|i|+1)(|i|+2)}{6}+\mathcal N_{|i|}(i),
%$$
%for some $\mathcal N_m$ providing a numbering of indices of length m. 
%For instance we can choose $\mathcal N_m$ to count indices according to the linear order $\prec$ within the layer $m$ and in this case the numbering corresponds to $\mathcal N(i) = \sum_{j\prec i} 1$, 
%or choose $\mathcal N_m(i)=(i2+i3)(i2+i3+1)/2+i3+1$ for all $i$ such that $|i|=m$. 

A road map was proposed in \cite{IGS} to prove approximation properties of  %GPWs, and apply in general for 
GPW functions. It can be summarized as follows: % The following outlines the different steps 
\begin{enumerate}
\item\label{step:prelim} for each quasi-Trefftz function, express all the basic parameters in \Cref{Algo:PbS,Algo:AbS,Algo:PolS} in terms of the free parameters that are not set to $0$;
\item\label{step:idref} identify a reference case, here a classical PW case;
\item\label{step:propref} study useful properties of the reference matrix; %basis {\color{red} or maybe the $Z_l^m$ functions (actually more of the corresponding matri $M$ - that's all we need)};
\item\label{step:link} establish a link between each of the quasi-Trefftz cases and the reference case;
\item\label{step:proof}  prove the approximation properties of quasi-Trefftz bases.
\end{enumerate}
While \Cref{step:idref,step:propref,step:proof}  are case-independent, \Cref{step:prelim,step:link} will be treated separately for each family of GPW functions.
These two key points rely on understanding how the entries of the linear system matrices depend on the initialization of our GPW functions, emphasizing their properties shared by corresponding entries on a given row as well as their differences.
Two important questions about these matrices concern their rank. (i) How large of a rank can they have? (ii) What particular choice of angles in the initialization parameters can guarantee the maximal rank? These will lead the choice (i) of how many different quasi-Trefftz functions to define, and (ii) of how to choose the initialization angles.
%{\color{red} there are two questions about the matrices: 1) how large of a rank can they have, 2) what particular choice of angles can reach the maximal rank}
%
%{\color{red} maybe say why we care about the rank: maximal rank (leading to maximal rank of the QT matrices) means we can find a linear combination of our basis functions to approximate any sol of the PDE OR MAYBE NOT}
As a by-product, the resulting families of GPW functions will be proved to be linearly independent.
It seems important to underline the fundamental part that the interplay of \Cref{hyp:PDop} and the choice of initialization will play in the rest of this section.

By contrast, thanks to the choice of initialization for the polynomial quasi-Trefftz functions, the polynomial quasi-Trefftz can be studied directly.

Given an order $n$ for the approximation property \eqref{eq:TEapp}, the order $q$ of the quasi-Trefftz property will be chosen to guarantee a similar construction for all the polynomial coefficients of quasi-Trefftz basis functions that will appear in the Taylor expansion  \eqref{eq:TEapp}. This will require to set $q\geq n-1$, it is then sufficient to construct the quasi-Trefftz basis functions with the parameter $q$ satisfying:
 \begin{center}
   \fbox{ $ q=\max(n-1,1)$ }
\end{center}
This will be particularly helpful to describe all polynomial coefficients of the quasi-Trefftz basis functions in terms of the initialization parameters,  see Section \ref{ssec:Relate}.

%%%%%%%%%%%%%%%%%%%%%%%%%%%%%%%%%%%
\subsection{Preliminary results}
The goal is to investigate how the terms computed in \Cref{Algo:AbS,Algo:SubS,Algo:PbS}, %Algo:PolS}, 
namely $\mu_{\beta+2e_1}$ and $\lambda_{\beta+2e_1}$, % and $\nu_{\beta+2e_1}$, 
depend on the three free parameters  from the initialization process \eqref{eq:3params}.
 We will proceed by induction with respect to the layer $\ell$. In each case the result will rely on a careful inspection of the right-hand side $B^\ell$ of the subsystems.

\subsubsection{For amplitude-based GPWs}
For an amplitude-based GPW, we focus on investigating properties of $\{\mu_{i}\in \mathbb C,i\in\mathbb N_0^3, |i|\leq q+1\}$.
Here, the three non-zero free parameters in the initialization procedure are $\Lambda_1,\Lambda_2,\Lambda_3$.
All $\mu$s computed from \Cref{Algo:AbS,Algo:SubS} clearly appear to be polynomials with respect to these three free parameters, that is they are elements of $\mathbb C[\Lambda_1, \Lambda_2, \Lambda_3]$.
Moreover, as first noted in \cite{LMinterp}, the initialization ensures that
$$
Q_N = 0, \text{ where } 
Q_N(\Lambda_1,\Lambda_2,\Lambda_3):= \sum_{k=1}^3
\coefc{2e_k}[\mathbf 0]\left(\Lambda_k\right)^2 
+
\sum_{1\leq k<k'\leq 3}
\coefc{e_k+e_{k'}}[\mathbf 0] \Lambda_{k'}\Lambda_k-\cst^2
$$
which turns our attention to elements of $\mathbb C[\Lambda_1,\Lambda_2,\Lambda_3]/ (Q_N)$ instead.
We will therefore investigate how the other $\mu$s can be expressed in terms of the three free parameters,  $\Lambda_1,\Lambda_2,\Lambda_3$.

\begin{lmm}\label{lmm:l2e1Ab}
Given $q\in\mathbb N$, a point $\pG\in\mathbb R^3$, a set of complex-valued functions  $c = \{\funcci,i\in\mathbb N_0^3, |i|\leq 2\}$ is assumed to satisfy \Cref{hyp:PDop}.

Consider any amplitude-based GPW associated to differential operator $\Lc$, $J(\varx) := \polQ(\varx-\pG)\exp \Lambda \cdot\Big(\varx-\pG \Big)$ with $\displaystyle Q:=\sum_{i\in\mathbb N_0^3, |i|\leq q+1} \mu_i \mathbf X^i $, constructed via \Cref{Algo:AbS,Algo:SubS}, with the initialization introduced in \Cref{sec:norm} for $\cst\in\mathbb C^*$ a unit vector $\dir\in\mathbb S^2$.
Then $\mu_{2e_1}$ can be expressed as a polynomial of degree at most equal to $1$ in $\mathbb C[\Lambda_1,\Lambda_2,\Lambda_3]$, with coefficients  depending on $\cst$ yet independent of $\dir$.
\end{lmm}
\begin{proof}
From the formulas in \Cref{Algo:SubS,Algo:AbS} for $\ell=0$ and $\varb=\bfz$ we get:
$$\left\{
\begin{array}{l}
 \mu_{2e_1} = \displaystyle
 \frac{1}{2\coefc{2e_1}[\bfz]}\Bigg(\mathsf {B^0_{\mathbf 0}}-\sum_{k=2}^3
2\coefc{2e_k}[\bfz]\mu_{2e_k}-\sum_{1\leq k<k'\leq 3}
\coefc{e_k+e_{k'}}[\bfz]\mu_{e_k+e_{k'}}
\Bigg),\\
\mathsf B^0_{\bfz} =  \displaystyle-\sum_{k=1}^3 \coefc{2e_k}[\bfz](2\Lambda_k\mu_{e_k}+\Lambda_k^2 \mu_\bfz)
-\sum_{1\leq k < k'\leq 3} \coefc{e_k+e_{k'}}[\bfz]\Big(\mu_{e_k+e_{k'}}+\Lambda_k \mu_{e_{k'}}+\Lambda_{k'} \mu_{e_{k}} + \Lambda_k\Lambda_{k'}\mu_\bfz\Big)\\
\phantom{\mathsf B^0_{\bfz} = }
-\displaystyle \sum_{k=1}^3\coefc{e_k}[\bfz](\mu_{e_k}+\Lambda_k \mu_\bfz)
- \coefc{\bfz}[\bfz]\mu_\bfz,
\end{array}\right.
$$
so the initialization implies:
$$\left\{
\begin{array}{l}
 \mu_{2e_1} = \displaystyle
 \frac{1}{2\coefc{2e_1}[\bfz]}\mathsf {B^0_{\bfz}},\\
\mathsf B^0_{\bfz} =  \displaystyle-\sum_{k=1}^3 \coefc{2e_k}[\bfz]\Lambda_k^2
-\sum_{1\leq k < k'\leq 3} \coefc{e_k+e_{k'}}[\bfz] \Lambda_k\Lambda_{k'}
-\displaystyle \sum_{k=1}^3\coefc{e_k}[\bfz]\Lambda_k 
- \coefc{\bfz}[\bfz],\\
\phantom{\mathsf B^0_{\bfz} }=  -\displaystyle Q_N(\Lambda_1,\Lambda_2,\Lambda_3)-\cst^2
-\displaystyle \sum_{k=1}^3\coefc{e_k}[\bfz]\Lambda_k 
- \coefc{\bfz}[\bfz],
\end{array}\right.
$$
Therefore, since $ Q_N(\Lambda_1,\Lambda_2,\Lambda_3)=0$, we obtain
$$
 \mu_{2e_1} = 
 \frac{1}{2\coefc{2e_1}[\bfz]}
 \left(-\cst^2
-\displaystyle \sum_{k=1}^3\coefc{e_k}[\bfz]\Lambda_k 
- \coefc{\bfz}[\bfz]\right),
$$
which proves the claim since $\cst$ is a fixed constant according to the initialization.
\end{proof}
\begin{prop}
\label{prop:Abmus}
Given $q\in\mathbb N$ and a point $\pG\in\mathbb R^3$, a set of complex-valued functions  $c = \{\funcci,i\in\mathbb N_0^3, |i|\leq 2\}$ is assumed to satisfy \Cref{hyp:PDop}.

Consider any amplitude-based GPW, $J(\varx) := \polQ(\varx-\pG)\exp \Lambda \cdot\Big(\varx-\pG \Big)$ with $\displaystyle Q:=\sum_{i\in\mathbb N_0^3, |i|\leq q+1} \mu_i \mathbf X^i $, constructed via \Cref{Algo:AbS,Algo:SubS}, with the initialization introduced in \Cref{sec:norm} for $\cst\in\mathbb C^*$ a unit vector $\dir\in\mathbb S^2$.
Then,
 for all $\ell$ from 0 to $q-1$ and all $\beta\in(\mathbb N_0)^3$ such that $|\beta|=\ell$, $\mu_{\varb+2e_1}$ 
can be expressed as a polynomial of degree at most equal to 
$|\varb|+1$
 in $\mathbb C[\Lambda_1,\Lambda_2,\Lambda_3]$, with coefficients depending on $\cst$ yet independent of $\dir$.
\end{prop}
\begin{proof}
The echelon form of the system is crucial here.
We will proceed by induction with respect to $\varb$ according to the linear order $\prec$, which corresponds to the order in which the $\mu_{\varb+2e_k} $ terms are computed in the algorithms backward substitution. 

The case $\varb=\mathbf 0$ is precisely \Cref{lmm:l2e1Ab}.

Assume $\varb\in(\mathbb N_0)^3$ with $0\prec\varb$ is such that the result holds for all $\varb'\prec \varb$: $\mu_{\varb'+2e_1}$ can be expressed as a polynomial of degree at most equal to $|\varb'|+1$ in $\mathbb C[\Lambda_1,\Lambda_2,\Lambda_3]$. Then $ \mu_{\varb+2e_1}$ is computed according to \Cref{Algo:SubS,Algo:AbS}. Hence, since 
$$
2 \Lambda_k (\gamma_k+1)\mu_{\gamma+e_k} + \Lambda_k^2\mu_{\gamma} = 
\Lambda_k(\gamma_{k'}+1)\mu_{\gamma+e_{k'}}
+\Lambda_{k'}(\gamma_{k}+1)\mu_{\gamma+e_k} + \Lambda_k\Lambda_{k'}\mu_\gamma
\text{ for } k'=k,
$$
we can gather these terms in a sum over $1\leq k\leq k'\leq 3$ and
 $ \mu_{\varb+2e_1}$ can be written as:
\begin{equation}\label{eq:mb+2e1}
\begin{array}{ll}
 \mu_{\varb+2e_1} = &\displaystyle
 \frac{1}{(\varb_1+2)(\varb_1+1)\coefc{2e_1}[\bfz]}\Bigg(
-\sum_{\gamma\in\mathbb N_0^3;\gamma\in\mathbb N_0^3;\gamma< \varb}
\sum_{k=1}^3
\coefc{2e_k}[\varb-\gamma](\gamma_k+2)(\gamma_k+1)\mu_{\gamma+2e_k} \\&\displaystyle
-\sum_{\gamma\in\mathbb N_0^3;\gamma\leq \varb}
\sum_{1\leq k<k'\leq 3}
\coefc{e_k+e_{k'}}[\varb-\gamma](\gamma_k+1)(\gamma_{k'}+1)\mu_{\gamma+e_k+e_{k'}}  \\&\displaystyle
-\sum_{\gamma\in\mathbb N_0^3;\gamma\in\mathbb N_0^3;\gamma\leq \varb}
\sum_{1\leq k\leq k'\leq 3}
\coefc{e_k+e_{k'}}[\varb-\gamma]
\Big(\Lambda_k(\gamma_{k'}+1)\mu_{\gamma+e_{k'}}  +\Lambda_{k'}(\gamma_{k}+1)\mu_{\gamma+e_k} + \Lambda_k\Lambda_{k'}\mu_\gamma\Big) \\&\displaystyle
-\sum_{\gamma\leq \varb}
\left(
\sum_{k=1}^3
\coefc{e_k}[\varb-\gamma]\Big((\gamma_{k}+1)\mu_{\gamma+e_k}   +\Lambda_k \mu_\gamma \Big)
+
\coefczr[\varb-\gamma] \mu_\gamma
\right)
 \\&\displaystyle-\sum_{k=2}^3
(\varb_k+2)(\varb_k+1)\coefc{2e_k}[\bfz]\mu_{\varb+2e_k}
 \\&\displaystyle
-\sum_{1\leq k<k'\leq 3}
(\varb_k+1)(\varb_{k'}+1)\coefc{e_k+e_{k'}}[\bfz]\mu_{\varb+e_k+e_{k'}}
\Bigg).
\end{array}
\end{equation}
On the right hand side we observe that the $\mu$ terms fall in one of two categories as elements of $\mathbb C[\Lambda_1,\Lambda_2,\Lambda_3]$:
\begin{enumerate}
\item $\mu_i$ with $i_1\in\{0,1\}$, chosen in the initialization process, either $\mu_\bfz=1$ or otherwise the terms set to zero, 
\item $\mu_i$ with $i_1\geq 2$, computed at a previous iteration for $\varb'=i-2e_1\prec\varb$.
\end{enumerate}
The linear terms can be listed as follows.
\begin{center}
{\renewcommand{\arraystretch}{1.1}
\begin{tabular}{|c||c|c|c|c|c|c|}
\hline
Terms &  Indices & $0$ & $1$ & $\mu_{\varb'+2e_1}$ \\\hline\hline
$\mu_{\gamma+2e_k}$& 
$\begin{array}{c}
\gamma<\varb(\neq\mathbf 0),k\in\{1,2,3\}
\\ (\gamma+2e_k)_1\in\{0,1\}
\end{array}$
& 
\ding{51}
 & & \\\hline
$\mu_{\gamma+2e_k}$& 
$\begin{array}{c}
\gamma<\varb(\neq\mathbf 0),k\in\{1,2,3\}
\\ (\gamma+2e_k)_1>1
\end{array}$
& & & 
$ \varb' =\gamma+2(e_k-e_1) $ \\\hline
$\mu_{\gamma+e_k+e_{k'}}$ & 
$\begin{array}{c}
\gamma\leq\varb,1\leq k<k'\leq 3
\\ (\gamma+e_k+e_{k'})_1\in\{0,1\}
\end{array}$
 &
 \ding{51}&&\\\hline
$\mu_{\gamma+e_k+e_{k'}}$ & 
$\begin{array}{c}
\gamma\leq\varb,1\leq k<k'\leq 3
\\ (\gamma+e_k+e_{k'})_1>1
\end{array}$
 & & & 
$\varb'=\gamma+e_k+e_{k'}-2e_1$ 
\\\hline
$\mu_{\gamma+e_k}$ & 
$\begin{array}{c}
\gamma\leq\varb,k\in\{1,2,3\}
\\ ({\gamma+e_k})_1\in\{0,1\}
\end{array}$
 & \ding{51} & & 
\\\hline
$\mu_{\gamma+e_k}$ & 
$\begin{array}{c}
\gamma\leq\varb,k\in\{1,2,3\}
\\ ({\gamma+e_k})_1>1
\end{array}$
 & & & $\varb'=\gamma+e_k-2e_1$
\\\hline
$\mu_{\gamma}$ & 
$\gamma=\mathbf 0$
 & & \ding{51} & 
\\\hline
$\mu_{\gamma}$ & 
$\begin{array}{c}
\gamma\leq\varb
\\ ({\gamma})_1\in\{0,1\},\gamma\neq\mathbf 0
\end{array}$
 & \ding{51} & & 
\\\hline
$\mu_{\gamma}$ & 
$\begin{array}{c}
\gamma\leq\varb,
\\ ({\gamma})_1>1
\end{array}$
 & & & $\varb'=\gamma-2e_1$
\\\hline
$\mu_{\varb+2e_k}$& 
$\begin{array}{c}
k\in\{2,3\}
\\ (\varb+2e_k)_1\in\{0,1\}
\end{array}$
& 
\ding{51}
 & & \\\hline
$\mu_{\varb+2e_k}$& 
$\begin{array}{c}
k\in\{2,3\}
\\ (\varb+2e_k)_1>1
\end{array}$
& & & 
$ \varb' =\varb+2(e_k-e_1) $ \\\hline
$\mu_{\varb+e_k+e_{k'}}$ & 
$\begin{array}{c}
1\leq k<k'\leq 3
\\ (\varb+e_k+e_{k'})_1\in\{0,1\}
\end{array}$
 &
 \ding{51}&&\\\hline
$\mu_{\varb+e_k+e_{k'}}$ & 
$\begin{array}{c}
1\leq k<k'\leq 3
\\ (\varb+e_k+e_{k'})_1>1
\end{array}$
 & & & 
$\varb'=\varb+e_k+e_{k'}-2e_1$ 
\\\hline
 \end{tabular}}
\end{center}
Hence the linear terms can be expressed as elements of $\mathbb C[\Lambda_1,\Lambda_2,\Lambda_3]$ as either $0$, $1$, or by induction hypothesis  as a polynomial of degree at most equal to $|\varb'|+1$. Moreover, in this last case, the values of $\varb'$ identified in the previous table are such that:
$$
\left\{\begin{array}{l}
\gamma<\varb\Rightarrow |\gamma+2(e_k-e_1)|+1<|\varb|+1,\\
\gamma\leq\varb\Rightarrow |\gamma+e_k+e_{k'}-2e_1|+1\leq |\varb|+1,\\
\gamma\leq\varb\Rightarrow |\gamma+e_k-2e_1|+1\leq |\varb|,
\end{array}\right.
\text{ and }
\left\{\begin{array}{l}
 |\varb+2(e_k-e_1)|+1=|\varb|+1,\\
 |\varb+e_k+e_{k'}-2e_1|+1= |\varb|+1.
\end{array}\right.
$$
In summary all the linear terms in the right hand side of \eqref{eq:mb+2e1} can be expressed as polynomials of degree at most equal to $|\varb|+1$ in $\mathbb C[\Lambda_1,\Lambda_2,\Lambda_3]$.

As for the non-linear terms, they all appear for indices $\gamma\leq\varb$ and $1\leq k \leq k'\leq 3$, and they can be identified as follows.
\begin{center}
{\renewcommand{\arraystretch}{1.1}
\begin{tabular}{|c||c|c|c|c|c|c|c|c|}
\hline 
Terms &  Indices & $0$ & $1$& $\Lambda_k$ or $\Lambda_{k'}$ & $\mu_{\varb'+2e_1}$&Case \\\hline\hline
$\Lambda_k\mu_{\gamma+e_{k'}}$& 
$(\gamma+e_{k'})_1\in\{0,1\}$& \ding{51} & &\ding{51}& & \\\hline
$\Lambda_k\mu_{\gamma+e_{k'}}$& 
$ (\gamma+e_{k'})_1>1$& & &\ding{51}  &$\varb'=\gamma+e_{k'}-2e_1$&1 \\\hline
$\Lambda_{k'}\mu_{\gamma+e_k}$ & 
$(\gamma+e_k)_1\in\{0,1\}$ & \ding{51}&&\ding{51}&&\\\hline
$\Lambda_{k'}\mu_{\gamma+e_k}$ & 
$(\gamma+e_k)_1>1$ & & &\ding{51}& $\varb'=\gamma+e_k-2e_1$ &2
\\\hline
$\Lambda_k\Lambda_{k'}\mu_\gamma$ & 
$\gamma=\mathbf 0$ & & \ding{51} &  \ding{51} \ding{51}&&3
\\\hline
$\Lambda_k\Lambda_{k'}\mu_\gamma$ & 
$ ({\gamma})_1\in\{0,1\}$, $\gamma\neq\mathbf 0$& \ding{51} & &   \ding{51} \ding{51}&&
\\\hline
$\Lambda_k\Lambda_{k'}\mu_\gamma$ & 
$ ({\gamma})_1>1$& & &   \ding{51} \ding{51}& $\varb'=\gamma-2e_1$ &4
\\\hline
$\Lambda_k\mu_\gamma$ & 
$\gamma=\mathbf 0$ & & \ding{51} &  \ding{51} &&5
\\\hline
$\Lambda_k\mu_\gamma$ & 
$ ({\gamma})_1\in\{0,1\}$, $\gamma\neq\mathbf 0$& \ding{51} & &   \ding{51} &&
\\\hline
$\Lambda_k\mu_\gamma$ & 
$ ({\gamma})_1>1$& & &   \ding{51} & $\varb'=\gamma-2e_1$ &6
\\\hline
 \end{tabular}}
\end{center}
Hence, as elements of  $\mathbb C[\Lambda_1,\Lambda_2,\Lambda_3]$, these  non-linear terms can be expressed either as zero or by induction hypothesis as a polynomial of degree at most equal to:
\begin{itemize}
\item $|\beta'|+2 \leq |\varb|+1$ in cases 1 and 2,
\item $2$ in case 3,
\item $|\varb'|+3\leq |\varb|+1$ in case 4,
\item $1$ in case 5,
\item $|\varb'|+2\leq|\varb|$ in case 6.
\end{itemize}
In summary all the non-linear terms in the right hand side of \eqref{eq:mb+2e1} can be expressed as polynomials of degree at most equal to $|\varb|+1$ in $\mathbb C[\Lambda_1,\Lambda_2,\Lambda_3]$.

Therefore $ \mu_{\varb+2e_1}$ in \eqref{eq:mb+2e1} can be expressed as a polynomial of degree at most equal to $|\varb|+1$ in $\mathbb C[\Lambda_1,\Lambda_2,\Lambda_3]$. This concludes the proof.
\end{proof}

\subsubsection{For phase-based GPWs}
For a phase-based GPW, we focus on investigating properties of $\{\lambda_{i}\in \mathbb C,i\in\mathbb N_0^3, |i|\leq q+1\}$.
Here, the three non-zero free parameters in the initialization procedure are $\lambda_{e_1},\lambda_{e_2},\lambda_{e_3}$.
All $\lambda$s computed from \Cref{Algo:PbS,Algo:SubS} clearly appear to be polynomials with respect to these three free parameters, that is they are elements of $\mathbb C[\lambda_{e_1},\lambda_{e_2},\lambda_{e_3}]$.
Moreover, as first noted in \cite{LMinterp}, the initialization ensures that
$$
P_N = 0, \text{ where } 
P_N(\lambda_{e_1},\lambda_{e_2},\lambda_{e_3}):= \sum_{k=1}^3
\coefc{2e_k}[\mathbf 0]\left(\lambda_{e_k}\right)^2 
+
\sum_{1\leq k<k'\leq 3}
\coefc{e_k+e_{k'}}[\mathbf 0] \lambda_{e_{k'}}\lambda_{e_k}-\cst^2
$$
which turns our attention to elements of $\mathbb C[\lambda_{e_1},\lambda_{e_2},\lambda_{e_3}]/ (P_N)$ instead.
We will therefore investigate how the other $\lambda$s can be expressed in terms of the three free parameters,  $\lambda_{e_1},\lambda_{e_2},\lambda_{e_3}$.
\begin{lmm}\label{lmm:l2e1Pb}
Given $q\in\mathbb N$ and a point $\pG\in\mathbb R^3$, a set of complex-valued functions  $c = \{\funcci,i\in\mathbb N_0^3, |i|\leq 2\}$ is assumed to satisfy \Cref{hyp:PDop}.

Consider any phase-based GPW associated to the partial differential operator $\Lc$, $G(\varx) :=\exp  \polP(\varx-\pG)$ with $\displaystyle P:=\sum_{i\in\mathbb N_0^3, |i|\leq q+1} \lambda_i \mathbf X^i $, constructed via \Cref{Algo:PbS,Algo:SubS}, with the initialization introduced in \Cref{sec:norm} for $\cst\in\mathbb C^*$ and a unit vector $\dir$.
Then $\lambda_{2e_1}$ can be expressed as a polynomial of degree at most equal to 1 in $\mathbb C[\lambda_{e_1},\lambda_{e_2},\lambda_{e_3}]$, with coefficients depending on $\cst$ yet independent of $\dir$.
\end{lmm}
\begin{proof}
Since  $\lambda_{2e_1}$ is computed at iteration $\ell=0$ in \Cref{Algo:PbS}, according to \Cref{Algo:SubS}, we have:
$$
%\begin{array}{rr}
 \lambda_{2e_1} = %&\displaystyle
 \frac{1}{2\coefc{2e_1}[\bfz]}\Bigg(\mathsf {B^0_{\mathbf 0}}-\sum_{k=2}^3
2\coefc{2e_k}[\bfz]\lambda_{2e_k}
% \\&\displaystyle
-\sum_{1\leq k<k'\leq 3}
\coefc{e_k+e_{k'}}[\bfz]\lambda_{e_k+e_{k'}}
\Bigg),
%\end{array}
$$
so the initialization then implies:
\begin{equation}
\label{eq:lambda2e1Pb}
%\begin{array}{rr}
 \lambda_{2e_1} = %&\displaystyle
 \frac{1}{2\coefc{2e_1}[\bfz]}\mathsf {B^0_{\mathbf 0}},
%\end{array}
\end{equation}
By definition, the right hand side $B^0=\left[B^0_{\mathbf 0}\right]\in\mathbb C$ is:
$$
\begin{array}{rl}
B^0_{\mathbf 0}
&\displaystyle=
-\sum_{k=1}^3
\coefc{2e_k}[\mathbf 0]\left( \lambda_{e_{k}}\right)^2 
-
\sum_{1\leq k<k'\leq 3}
\coefc{e_k+e_{k'}}[\mathbf 0]\Bigg(\lambda_{e_k+e_{k'}} 
+ \lambda_{e_{k'}}\lambda_{e_k}\Bigg)
-
\sum_{k=1}^3
\coefc{e_k}[\mathbf 0]\lambda_{e_k}   
 -
\coefczr[\mathbf 0],\\
&\displaystyle=
-P_N(\lambda_{e_1},\lambda_{e_2},\lambda_{e_3})-\cst^2
-
\sum_{1\leq k<k'\leq 3}
\coefc{e_k+e_{k'}}[\mathbf 0]\lambda_{e_k+e_{k'}} 
-
\sum_{k=1}^3
\coefc{e_k}[\mathbf 0]\lambda_{e_k}   
 -
\coefczr[\mathbf 0],
\end{array}
$$
so from the initialization, imposing all the $\lambda_{e_k+e_{k'}} $ to be zero as well as $P_N=0$, we can express the right hand side as:
\begin{equation}
\label{eq:B0Pb}
B^0_{\mathbf 0}
=
-\cst^2
-
\sum_{k=1}^3
\coefc{e_k}[\mathbf 0]\lambda_{e_k}   
 -
%\sum_{|\alpha|=0}
\coefczr[\mathbf 0].
\end{equation}
As a result, combining \eqref{eq:lambda2e1Pb} and  \eqref{eq:B0Pb} yields:
$$
 \lambda_{2e_1} =
 \frac{1}{2\coefc{2e_1}[\bfz]}\left( - \cst^2
-
\sum_{k=1}^3
\coefc{e_k}[\mathbf 0]\lambda_{e_k}   
 -
%\sum_{|\alpha|=0}
\coefczr[\mathbf 0]\right),
$$
which proves the claim since $\cst$ is a fixed constant according to the initialization.
\end{proof}
\begin{prop}
\label{prop:Pblambdas}
Given $q\in\mathbb N$ and a point $\pG\in\mathbb R^3$, a set of complex-valued functions  $c = \{\funcci,i\in\mathbb N_0^3, |i|\leq 2\}$ is assumed to satisfy \Cref{hyp:PDop}.

Consider any phase-based GPW associated to the partial differential operator $\Lc$, $G(\varx) :=\exp  \polP(\varx-\pG)$ with $\displaystyle P:=\sum_{i\in\mathbb N_0^3, |i|\leq q+1} \lambda_i \mathbf X^i $, constructed via \Cref{Algo:PbS,Algo:SubS}, with the initialization introduced in \Cref{sec:norm} for $\cst\in\mathbb C^*$ and a unit vector $\dir$.
Then, for all $\ell$ from 0 to $q-1$ and all $\beta\in(\mathbb N_0)^3$ such that $|\beta|=\ell$, $\lambda_{\varb+2e_1}$ can be expressed as a polynomial of degree at most equal to $|\varb|+1$ in $\mathbb C[\lambda_{e_1},\lambda_{e_2},\lambda_{e_3}]$, with coefficients depending on $\cst$ yet independent of $\dir$.
\end{prop}

\begin{proof}
Here again the echelon form of the system is crucial here.
We will proceed by induction with respect to $\varb$ according to the linear order $\prec$, which corresponds to the order in which the $\lambda_{\varb+2e_k} $ terms are computed by backward substitution in the algorithms. 

The case $\varb=\mathbf 0$ is precisely \Cref{lmm:l2e1Pb}.

Assume $\varb\in(\mathbb N_0)^3$ with $0\prec\varb$ is such that the result holds for all $\varb'\prec \varb$: $\lambda_{\varb'+2e_1}$ can be expressed as a polynomial of degree at most equal to $|\varb'|+1$ in $\mathbb C[\lambda_{e_1},\lambda_{e_2},\lambda_{e_3}]$. Then $ \lambda_{\varb+2e_1}$ is computed according to \Cref{Algo:SubS,Algo:PbS}. Hence, since 
$$
(\gamma_{k}-\eta_{k} +1)\lambda_{\gamma-\eta+e_{k}}(\eta_k +1)\lambda_{\eta+e_k} = 
(\gamma_{k'}-\eta_{k'} +1)\lambda_{\gamma-\eta+e_{k'}}(\eta_k +1)\lambda_{\eta+e_k}
\text{ for } k'=k,
$$
we can gather these terms in a sum over $1\leq k\leq k'\leq 3$ and
 $ \lambda_{\varb+2e_1}$ can be written as:
\begin{equation}\label{eq:lb+2e1}
\begin{array}{rl}
 \lambda_{\varb+2e_1} = &\displaystyle
 \frac{1}{(\varb_1+2)(\varb_1+1)\coefc{2e_1}[\bfz]}\Bigg(
-\sum_{\gamma\in\mathbb N_0^3;\gamma< \varb}
\sum_{k=1}^3
\coefc{2e_k}[\varb-\gamma](\gamma_k+2)(\gamma_k+1)\lambda_{\gamma+2e_k} 
 \\&
\displaystyle
-\sum_{\gamma\in\mathbb N_0^3;\gamma\leq \varb} 
\sum_{1\leq k<k'\leq 3}
\coefc{e_k+e_{k'}}[\varb-\gamma](\gamma_k+1)(\gamma_{k'}+1)\lambda_{\gamma+e_k+e_{k'}} \\&
\displaystyle
-\sum_{\gamma\in\mathbb N_0^3;\gamma\leq \varb}
\sum_{1\leq k\leq k'\leq 3}
\coefc{e_k+e_{k'}}[\varb-\gamma]
 \sum_{\eta\in\mathbb N_0^3;\eta\leq \gamma}(\gamma_{k'}-\eta_{k'} +1)\lambda_{\gamma-\eta+e_{k'}}(\eta_k +1)\lambda_{\eta+e_k}
 \\&\displaystyle
 -
\sum_{\gamma\in\mathbb N_0^3;\gamma\leq \varb} \sum_{k=1}^3
\coefc{e_k}[\varb-\gamma](\gamma_{k}+1)\lambda_{\gamma+e_k}   
-
 \coefczr[\varb]
 \\&\displaystyle
 -\sum_{k=2}^3
(\varb_k+2)(\varb_k+1)\coefc{2e_k}[\bfz]\lambda_{\varb+2e_k}
 \\&\displaystyle
-\sum_{1\leq k<k'\leq 3}
(\varb_k+1)(\varb_{k'}+1)\coefc{e_k+e_{k'}}[\bfz]\lambda_{\varb+e_k+e_{k'}}
\Bigg),
\end{array}
\end{equation}
On the right hand side we observe that the $\lambda$ terms fall in one of two categories as elements of $\mathbb C[\lambda_{e_1},\lambda_{e_2},\lambda_{e_3}]$:
\begin{enumerate}
\item $\lambda_i$ with $i_1\in\{0,1\}$, chosen in the initialization process, either a $\lambda_{e_k}$ or otherwise the terms set to zero, %so in particular it can each be expressed as a polynomial of degree at most equal to $|\varb|+1$ in $\mathbb C[\lambda_{e_1},\lambda_{e_2},\lambda_{e_3}]$,
\item $\lambda_i$ with $i_1\geq 2$, computed at a previous iteration for $\varb'=i-2e_1\prec\varb$.
\end{enumerate}
The linear terms can be listed as follows.
\begin{center}
{\renewcommand{\arraystretch}{1.1}
\begin{tabular}{|c||c|c|c|c|c|}
\hline
Terms &  Indices & $0$ & $\lambda_{e_k}$& $\lambda_{\varb'+2e_1}$ \\\hline\hline
$\lambda_{\gamma+2e_k}$& 
$\begin{array}{c}
\gamma<\varb(\neq\mathbf 0),k\in\{1,2,3\}
\\ (\gamma+2e_k)_1\in\{0,1\}
\end{array}$
& 
\ding{51}
 & & \\\hline
$\lambda_{\gamma+2e_k}$& 
$\begin{array}{c}
\gamma<\varb(\neq\mathbf 0),k\in\{1,2,3\}
\\ (\gamma+2e_k)_1>1
\end{array}$
& & &
$ \varb' =\gamma+2(e_k-e_1) $ \\\hline
$\lambda_{\gamma+e_k+e_{k'}}$ & 
$\begin{array}{c}
\gamma\leq\varb,1\leq k<k'\leq 3
\\ (\gamma+e_k+e_{k'})_1\in\{0,1\}
\end{array}$
 &
 \ding{51}&&\\\hline
$\lambda_{\gamma+e_k+e_{k'}}$ & 
$\begin{array}{c}
\gamma\leq\varb,1\leq k<k'\leq 3
\\ (\gamma+e_k+e_{k'})_1>1
\end{array}$
 & & &
$\varb'=\gamma+e_k+e_{k'}-2e_1$ 
\\\hline
$\lambda_{\gamma+e_k}$ & 
$\gamma=\mathbf 0,k\in\{1,2,3\}$
 & & \ding{51} &
\\\hline
$\lambda_{\gamma+e_k}$ & 
$\begin{array}{c}
\gamma\leq\varb,k\in\{1,2,3\}
\\ ({\gamma+e_k})_1\in\{0,1\},\gamma\neq\mathbf 0
\end{array}$
 & \ding{51} & &
\\\hline
$\lambda_{\gamma+e_k}$ & 
$\begin{array}{c}
\gamma\leq\varb,k\in\{1,2,3\}
\\ ({\gamma+e_k})_1>1
\end{array}$
 & & & $\varb'=\gamma+e_k-2e_1$
\\\hline
$\lambda_{\varb+2e_k}$& 
$\begin{array}{c}
k\in\{2,3\}
\\ (\varb+2e_k)_1\in\{0,1\}
\end{array}$
& 
\ding{51}
 & & \\\hline
$\lambda_{\varb+2e_k}$& 
$\begin{array}{c}
k\in\{2,3\}
\\ (\varb+2e_k)_1>1
\end{array}$
& & &
$ \varb' =\varb+2(e_k-e_1) $ \\\hline
$\lambda_{\varb+e_k+e_{k'}}$ & 
$\begin{array}{c}
1\leq k<k'\leq 3
\\ (\varb+e_k+e_{k'})_1\in\{0,1\}
\end{array}$
 &
 \ding{51}&&\\\hline
$\lambda_{\varb+e_k+e_{k'}}$ & 
$\begin{array}{c}
1\leq k<k'\leq 3
\\ (\varb+e_k+e_{k'})_1>1
\end{array}$
 & & &
$\varb'=\varb+e_k+e_{k'}-2e_1$ 
\\\hline
 \end{tabular}}
\end{center}
Hence the linear terms can be expressed as elements of $\mathbb C[\lambda_{e_1},\lambda_{e_2},\lambda_{e_3}]$ as either $0$, or a $\lambda_{e_k}$, or by induction hypothesis  as a polynomial of degree at most equal to $|\varb'|+1$. Moreover, in this last case, the values of $\varb'$ identified in the previous table are such that:
$$
\left\{\begin{array}{l}
\gamma<\varb\Rightarrow |\gamma+2(e_k-e_1)|+1<|\varb|+1,\\
\gamma\leq\varb\Rightarrow |\gamma+e_k+e_{k'}-2e_1|+1\leq |\varb|+1,\\
\gamma\leq\varb\Rightarrow |\gamma+e_k-2e_1|+1\leq |\varb|,
\end{array}\right.
\text{ and }
\left\{\begin{array}{l}
 |\varb+2(e_k-e_1)|+1=|\varb|+1,\\
 |\varb+e_k+e_{k'}-2e_1|+1= |\varb|+1.
\end{array}\right.
$$
In summary all the linear terms in the right hand side of \eqref{eq:lb+2e1} can be expressed as polynomials of degree at most equal to $|\varb|+1$ in $\mathbb C[\lambda_{e_1},\lambda_{e_2},\lambda_{e_3}]$.

As for the non-linear terms, all of the form $\lambda_{\gamma-\eta+e_{k'}}\lambda_{\eta+e_k}$ with $\eta\leq\gamma\leq\varb$, they can be described as follows.
\begin{center}
{\renewcommand{\arraystretch}{1.1}
\begin{tabular}{|c|c|c|c|c|}
\hline
 Indices & $0$ & $\lambda_{e_k}$ or $\lambda_{e_{k'}}$& $\lambda_{\varb'+2e_1}$ \\\hline\hline
$\eta=\gamma=\mathbf 0$
&  & 
\ding{51} \ding{51}&\\\hline
$\begin{array}{c}
\eta=\mathbf 0, \gamma\neq\mathbf 0
\\ ({\gamma-\eta+e_{k'}})_1\in\{0,1\}
\end{array}$ & \ding{51} & &\\\hline
$\begin{array}{c}
\eta=\mathbf 0
\\ ({\gamma-\eta+e_{k'}})_1>1
\end{array}$ & & \ding{51}&$\varb'=\gamma+e_{k'}-2e_1$ \\\hline
$\begin{array}{c}
\eta\neq\mathbf 0
\\ ({\eta+e_{k}})_1\in\{0,1\}
\end{array}$ &\ding{51}& &\\\hline
$\begin{array}{c}
 ({\eta+e_{k}})_1>1
\\ \eta=\gamma
\end{array}$& &\ding{51}&$\varb'=\eta+e_k-2e_1$ \\\hline
$\begin{array}{c}
 ({\eta+e_{k}})_1>1
\\ ({\gamma-\eta+e_{k'}})_1\in\{0,1\}
\end{array}$ & \ding{51}&& \\\hline
$\begin{array}{c}
 ({\eta+e_{k}})_1>1
\\ ({\gamma-\eta+e_{k'}})_1>1
\end{array}$ & & &$\begin{array}{c}\varb'=\eta+e_k-2e_1\\\varb''=\gamma-\eta+e_{k'}-2e_1\end{array}$ \\\hline
 \end{tabular}}
\end{center}
Hence the non-linear terms can be expressed as elements of $\mathbb C[\lambda_{e_1},\lambda_{e_2},\lambda_{e_3}]$ as either $0$, or a $\lambda_{e_k}\lambda_{e_{k'}}$, or by induction hypothesis  as a polynomial of degree at most equal to $|\varb'|+2$ or $|\varb'|+|\varb''|+2$. Moreover, in these last two cases, the values of $\varb'$ and $\varb''$ identified in the previous table are such that:
$$
\left\{\begin{array}{l}
\gamma\leq\varb\Rightarrow |\gamma+e_{k'}-2e_1|+2\leq|\varb|+1\\
\eta\leq\varb\Rightarrow |\eta+e_k-2e_1|+2\leq |\varb|+1\\
\gamma\leq\varb\Rightarrow |\eta+e_k-2e_1|+ | \gamma-\eta+e_{k'}-2e_1 | +2\leq |\varb|
\end{array}\right.
$$
In summary all the non-linear terms in the right hand side of \eqref{eq:lb+2e1} can be expressed as polynomials of degree at most equal to $|\varb|+1$ in $\mathbb C[\lambda_{e_1},\lambda_{e_2},\lambda_{e_3}]$.

Therefore $ \lambda_{\varb+2e_1}$ in \eqref{eq:lb+2e1} can be expressed as a polynomial of degree at most equal to $|\varb|+1$ in $\mathbb C[\lambda_{e_1},\lambda_{e_2},\lambda_{e_3}]$. This concludes the proof.
\end{proof}

\subsection{The reference matrix} 

Following our choice of initialization and the previous results presented in \Cref{prop:Abmus,prop:Pblambdas}, %,prop:Polnus}, 
we consider
%the natural choice of reference case is that of PWs (or exponential functions depending on propagating or evanescent medium). For 
a set of $p$ directions $\{ \mathbf d_l\}_{1\leq l\leq p}$ defined by two angles $(\theta_l,\varphi_l)$ as follows:
 $$
 \mathbf d_l  = 
 \begin{bmatrix}
 \sin\theta_l\sin\varphi_l\\\sin\theta_l\cos\varphi_l\\\cos\theta_l
 \end{bmatrix}.
 $$
To describe the natural choice of reference case, we then define - for a common value $\cst\in\mathbb C$ and the matrices $\Pmat$ and $\Dmat$ from \Cref{hyp:PDop} - the functions:
$$
\chi_l:
\varx\mapsto \exp \cst \Pmat\Dmat^{-1/2}\dir_l\cdot\Big(\varx-\pG \Big),$$
as well as the associated matrices $\mate^{n,p}$ in $\mathbb C^{{(n+1)(n+2)(n+3)}/{6} \times p}$:
$$
\mate^{n,p}_{\mathcal N(i),l}=T_{\chi_l}[i].
$$

Each entry $(\mathcal N(i),l)$ being a multiple of $(\sin\varphi_l) ^{i_1} (\cos\varphi_l)^{i_2}(\sin\theta_l)^{i_1+i_2}(\cos\theta_l)^{i_3}$, hence we define the associated matrices $\matr^{n,p}$ in $\mathbb C^{{(n+1)(n+2)(n+3)}/{6} \times p}$:
 $$
\matr^{n,p}_{\mathcal N(i),l}=\frac{(\sin\varphi_l) ^{i_1} (\cos\varphi_l)^{i_2}(\sin\theta_l)^{i_1+i_2}(\cos\theta_l)^{i_3}}{i!}.
 $$
 and show how they are related in the following result.
 \begin{lmm}Consider any set of $p$ directions $\{ \mathbf d_l\}_{1\leq l\leq p}$ as well as the matrices $\Pmat$ and $\Dmat$ from \Cref{hyp:PDop}, together with the associated ${(n+1)(n+2)(n+3)}/{6} \times p$ complex matrices $\mate^{n,p}$ and $\matr^{n,p}$.  
 There exists a block diagonal %non-singular 
 matrix $D^n\in \mathbb C^{ {(n+1)(n+2)(n+3)}/{6} \times {(n+1)(n+2)(n+3)}/{6}  }$ such that $\mate^{n,p}=D^n\matr^{n,p}$ and $\matr^{n,p}=D^n\mate^{n,p}$, independently of their number $p$ of columns.
 \end{lmm}
 \begin{proof}
The entries of the reference matrix $\matr^{n,p}$ are:
 $$
 T_{\chi_l}[i] = \frac{1}{i!}
 \left( \left(  \cst \Pmat\Dmat^{-1/2} \dir_l \right)_1\right)^{i_1}
 \left( \left(  \cst \Pmat\Dmat^{-1/2} \dir_l \right)_2\right)^{i_2}
 \left( \left(  \cst \Pmat\Dmat^{-1/2} \dir_l \right)_3\right)^{i_3},
 $$
 and we consider row blocks corresponding to increasing values of $|i|$: for $m$ from $0$ to $n$ we denote by $\mate^{n,p}_{(m)}$ and $\matr^{n,p}_{(m)}$ the  ${(m+1)(m+2)}/{2}\times p$ blocks of rows corresponding to all $i$ such that $|i|=m$, that is $\left(\mate^{n,p}_{(m)}\right)_{\mathcal N_m(i) l} = \left(\mate^{n,p}\right)_{\mathcal N(i)l}$.
 As a reminder, $ \cst\in\mathbb C$ is a constant,  while $\Pmat$ and $\Dmat$ depend only on the value of the PDE coefficients $c$ evaluated at $\pG$, hence the entries of the $3\times 3$ matrix  $ \cst \Pmat\Dmat^{-1/2}$ are independent of $l$; moreover the matrix $\cst \Pmat\Dmat^{-1/2}$ is non-singular under \Cref{hyp:PDop}.
Thus, according to \Cref{lmm:prodsvecents} below for $A:=\cst \Pmat\Dmat^{-1/2}$, for all $m$ from $0$ to $n$ there exist an %non-singular
 ${(m+1)(m+2)}/{2}\times{(m+1)(m+2)}/{2}$ matrix $\mathfrak C^{(m)}$, defined entry-wise by $\left(\mathfrak C^{(m)}\right)_{\mathcal N_m(i)\mathcal N_m(j)}=\mathfrak c_{ij}$, such that 
%for any vector defined for any $(\theta,\varphi)\in\mathbb R^2$ as:
% $$
% \mathbf d  = 
% \begin{bmatrix}
% \sin\theta\sin\varphi\\\sin\theta\cos\varphi\\\cos\theta
% \end{bmatrix},
% $$
% then 
$\mate^{n,p}_{(m)}=\mathfrak C^{(m)}\matr^{n,p}_{(m)}$
and $\matr^{n,p}_{(m)}=\mathfrak C^{(m)}\mate^{n,p}_{(m)}$.

Defining the block diagonal matrix $D^n:=\text{diag}\left(\mathfrak C^{(0)},\mathfrak C^{(1)},...,\mathfrak C^{(n)}\right)$, %which is non-singular since each of the diagonal blocks is non-singular, 
we then obtain the desired property: $\mate^{n,p}=D^n\matr^{n,p}$ and $\matr^{n,p}=D^n\mate^{n,p}$.
 \end{proof}
\begin{lmm}
\label{lmm:prodsvecents}
Let $A\in\mathbb C^{3\times 3}$ be non-singular. Then, for any vector $Y\in\mathbb C^3$,  any product of powers of the entries of $AY$ with a total power $m\in \mathbb N$ can be written as a linear combination - depending only on $A$ but independent of $Y$ - of products of powers of the entries of $Y$ each one of the products having a total power equal to $m$. Mathematically speaking:
\begin{equation}
\label{eq:AYvsY}
\forall i\in(\mathbb N_0)^3, 
\exists\ \{\mathfrak c_{ij}\in\mathbb C,j\in(\mathbb N_0)^3, |j|=|i|\} \text{ such that }
\forall Y\in\mathbb C^3,\
\prod_{k=1}^3 \big((AY)_k\big)^{i_k} 
=
\sum_{j\in\mathbb N_0^3;|j|=|i|} \mathfrak c_{ij} \prod_{l=1}^3 (Y_l)^{j_l} .
\end{equation}

Moreover: %if $A$ is non-singular,
\begin{equation}
\label{eq:XvsA-1X}
\forall X\in\mathbb C^3,\
\prod_{k=1}^3 \big(X_k\big)^{i_k} 
=
\sum_{j\in\mathbb N_0^3;|j|=|i|} \mathfrak c_{ij} \prod_{l=1}^3 \left((A^{-1}X)_l\right)^{j_l} .
\end{equation}
% then the ${(m+1)(m+2)}/{2}\times{(m+1)(m+2)}/{2}$ matrix defined entry-wise by $\left(\mathfrak C^{(m)}\right)_{\mathcal N_m(i)\mathcal N_m(j)}=\mathfrak c_{ij}$, is non-singular.
%
 %We can then define the  ${(m+1)(m+2)}/{2}\times{(m+1)(m+2)}/{2}$ matrix $\mathfrak C^{(m)}$ entry-wise by $\left(\mathfrak C^{(m)}\right)_{\mathcal N_m(i)\mathcal N_m(j)}=\mathfrak c_{ij}$.
\end{lmm}
\begin{rmk}
The result still holds for $n\in\mathbb N$ and $A\in\mathbb C^{n\times n}$, the proof is more tedious as it requires the introduction of more indices. In this article we only use the $n=3$ case hence we do not prove the more general case.
\end{rmk}
\begin{proof}
Since, for $k\in\{1,2,3\}$ and any $3\times 3$ matrix $A$, we have:
$$
\big((AY)_k\big)^{i_k}
 = \left(\sum_{l=1}^3 A_{kl}Y_l\right)^{i_k}
 = \sum_{j\in\mathbb N_0^3;|j|=i_k} \frac{(i_k)!}{j!} \prod_{l=1}^3(A_{kl}Y_l)^{j_l}
 = \sum_{j\in\mathbb N_0^3;|j|=i_k}\left( \frac{(i_k)!}{j!} \prod_{\tilde l=1}^3(A_{k\tilde l})^{j_{\tilde l}}\right) \prod_{l=1}^3(Y_l)^{j_l},
$$
then we get:
$$
\prod_{k=1}^3\big((AY)_k\big)^{i_k}  
= \sum_{j\in\mathbb N_0^3;|j|=i_1}\sum_{j'\in\mathbb N_0^3;|j'|=i_2}\sum_{j''\in\mathbb N_0^3;|j''|=i_3} \frac{i!}{j!(j')!(j'')!}
\left(\prod_{\tilde l=1}^3(A_{1\tilde l})^{j_{\tilde l}}(A_{2\tilde l})^{j'_{\tilde l}}(A_{3\tilde l})^{j''_{\tilde l}}\right) \prod_{l=1}^3(Y_l)^{j_l+j'_l+j''_l}.
$$
This concludes the proof of \eqref{eq:AYvsY}, as $\displaystyle \sum_{l=1}^3\big(j_l+j'_l+j''_l\big) = |j|+|j'|+|j''| = i_1+i_2+i_3 = |i| $.

Moreover, if $A$ is non-singular, then, for any $X\in\mathbb C^3$, \eqref{eq:AYvsY} applied to $Y=A^{-1}X$ shows \eqref{eq:XvsA-1X}.%that $\mathfrak C^{(m)}$ is indeed non-singular.
\end{proof}
As a direct consequence we have the following results.
\begin{cor}
\label{cor:rkmatexp}
Given $(n,p)\in\mathbb N^2$ and any choice of $p$ directions $\{ \mathbf d_l\in\mathbb S^2\}_{1\leq l\leq p}$,
the rank of $\mate^{n,p}$ is equal to the rank of $\matr^{n,p}$. 
\end{cor}
\begin{cor}
\label{cor:submatexp}
 Similarly, any pair of sub-matrices of $\mate^{n,p}$ and $\matr^{n,p}$ corresponding to removing the \underline{same} sets of all rows $i$ with $|i|=m$ for a given $m$ - in particular removing all rows $i$ with $|i|>1$ - have the same rank.
\end{cor}
The next step is to study the rank of these matrices, and the following step will be to relate it to the rank of the GPW matrices $\mata$ and $\matp$.

%%%%%%%%%%%%%%%%%%%%%%%%%%%%%%%%%%%
\subsection{Properties of the reference matrix}
We are interested here in the rank of the reference matrix, $\matr^{n,p}$, or equivalently the rank of the exponential matrix $\mate^{n,p}$ according to \Cref{cor:rkmatexp},  as well as the rank of some of their useful sub-matrices. In particular the rank's value depends on the number of columns $p$, and we will see that:
\begin{itemize}
\item independently of the number of columns $p$, the rank is at most $(n+1)^2$,
\item while there exist sets of $p=(n+1)^2$ directions that guarantee the rank of the corresponding matrix to be maximal, that is $(n+1)^2$.
\end{itemize}

%In order to choose $p$, w
To study this reference matrix, we first remark that its entries are a set of functions, on each row, evaluated at a set of points, on each column. The following result provides the maximum possible rank of such matrices depending on the dimension of the space generated by the function set.
\begin{lmm}
\label{lmm:matsfuncsevals}
Consider a set $\mathcal D$ of $N_f$ complex-valued functions defined on a domain $\omega$ (in any dimension), denoted $d_k$ for $k$ from $1$ to $N_f$, while  the dimension of $span \ \mathcal D$ is $N_d<N_f$. We consider any $N_f\times N_p$ matrix $\mathsf D$ defined element wise by evaluating the elements of the function space $\mathcal D$ at a set of $N_p$ points $\{\Theta_l\in\omega\}_{1\leq l\leq N_p}$, namely:
$$
\mathsf D_{kl} = d_k(\Theta_l).
$$
Then the rank of $\mathsf D$ is at most equal to $N_d$.
\footnote{
This argument was previously presented in our roadmap paper \cite{IGS} for a particular case defined by the function space $ \mathcal F_n=\{\theta\mapsto\cos^k \theta\sin^{K-k}\theta/(k!(K-k)!), 0\leq k\leq K\leq n\}$,  containing $\# \mathcal F_n=(n+1)(n+2)/2$ functions and spanning a space of dimension $\dim span \mathcal F_n=2n+1$.
}

Moreover, if ${\mathcal B}=\{y_k, 1\leq k\leq N_b\}$ is any generating set of $span\ \mathcal D$, then the corresponding $N_b\times N_p$ matrix ${\mathsf B}$, namely:
$$
{\mathsf B}_{kl} = {y}_k(\Theta_l),
$$
has the same rank as $\mathsf D$.
% if $\overline{\mathcal D}_n\subset\mathcal D_n$ is a generating set for $span\ \mathcal D_n$, the submatrix of $\mathsf D$, denoted $\overline{\mathsf D}$ and obtained by removing the rows corresponding to each $d_k\in\mathcal D_n\backslash\overline{\mathcal D}_n$, has the same rank as $\mathsf D$ itself. And similarly, if $\widetilde{\mathcal D}_n$ is any generating set for $span\ \mathcal D_n$, then the corresponding matrix $\widetilde{\mathsf D}$, namely:
%$$
%\widetilde{\mathsf D}_{kl} = \widetilde{d}_k(\Theta_l)
%$$
%also has the same rank.
\end{lmm}
\begin{proof}
Because of the number of elements in $\mathcal D$ and the dimension of its span, there exists a matrix $\mathsf C\in\mathbb C^{(N_f-N_d)\times N_f}$, of rank $N_f-N_d$ such that
$$
\forall k'\in\mathbb N, k'\leq N_f-N_d,\sum_{k=1}^{N_f}\mathsf C_{k'k}d_k = 0.
$$
In particular, independently of the number $N_p$ of columns of $\mathsf D$, this yields: 
$$
\mathsf C\mathsf D = 0_{(N_f-N_d)\times N_p}.
$$
Hence the $N_p$ columns of $\mathsf D$ belong to the kernel of $\mathsf C$, which is of dimension $N_d$ according to the rank-nullity theorem.
So indeed the rank of $\mathsf D$ is at most equal to $N_d$.

Moreover, consider a set $\widetilde{\mathcal D}\subset\mathcal D$, being a generating set for $span\ \mathcal D$, and the submatrix of $\mathsf D$, denoted $\widetilde{\mathsf D}$ and obtained by keeping only the rows corresponding to each $d_k\in\widetilde{\mathcal D}$.
%removing the rows corresponding to each $d_k\in\mathcal D_n\backslash\widetilde{\mathcal D}_n$. 
Since the space generated by the rows of $\mathsf D$ is the same as the space generated by the rows of $\widetilde{\mathsf D}$, %Since $\widetilde{\mathsf D}_n$ is a submatrix of $\mathsf D_n$, we have $rk \ \widetilde{\mathsf D}_n\leq rk\ \mathsf D_n$. Besides, since each function in $\mathcal D_n$ belongs to $span \ \widetilde{\mathcal D}_n$, then each row in $\mathsf D$ can be written as linear combination of the rows of $\widetilde{\mathsf D}$, hence $rk \ \widetilde{\mathsf D}\geq rk\ \mathsf D$. 
 $\mathsf D$ and $\widetilde{\mathsf D}$ have the same rank.
Similarly, consider ${\mathcal B}=\{y_k, 1\leq k\leq N_b\}$ any generating set of $span\ \mathcal D$, then consider the corresponding matrix ${\mathsf B}$.
Each of its rows can be written as a linear combination of the rows of $\widetilde{\mathsf D}$, since $\widetilde{\mathcal D}$ is a generating set of $span\ \mathcal D$, while each row of $\widetilde{\mathsf D}$ can be written as a linear combination of the rows of ${\mathsf B}$, since ${\mathcal B}$ is also a generating set of $span\ \mathcal D$.
This proves the second claim.
\end{proof}
 
To address the particular case of the matrix $\matr$, we then define the following functions and function space:
$$
\forall i\in(\mathbb N_0)^3, \ f_i(\theta,\varphi) = (\sin\varphi) ^{i_1} (\cos\varphi)^{i_2}(\sin\theta)^{i_1+i_2}(\cos\theta)^{i_3}/i!
\text{ and }
\mathcal F_n :=\{f_i, i\in(\mathbb N_0)^3, |i|\leq n\},
$$
$$
\widetilde{\mathcal F}_n :=\{f_i, i\in(\mathbb N_0)^3, |i|\leq n,i_1\in\{0,1\}\}.
$$
The set $\mathcal F_n$ contains ${(n+1)(n+2)(n+3)}/{6} $ elements, let's now identify the dimension of $Span\mathcal F_n$.
We will make use of the following functions and function spaces:
$$
Z_l^m(\theta,\varphi) = e^{\mi m\varphi}(\sin\theta)^{|m|}(\cos\theta)^{l-|m|}
\text{ and }
\mathcal G_n :=\{Z_l^m,0\leq l\leq n,-l\leq m\leq l\},
$$
and for the spherical harmonics with
$
C_l^m = \sqrt{\frac{2l+1}{4\pi}\frac{(l-m)!}{(l+m)!}}
$ and the
Legendre polynomials $P_l^m$ (see appendix \ref{sec:SHrem}):
$$
Y_l^m(\theta,\varphi):=C_l^m P_l^m(\cos\theta)e^{\mi m\varphi}
\text{ and }
\mathcal H_n :=\{Y_l^m,0\leq l\leq n,-l\leq m\leq l\}.
$$
\begin{lmm}
\label{lmm:dimFn}
Given $n\in\mathbb N$, the spaces of trigonometric functions and spherical harmonics are such that:
\begin{itemize}
\item $Span\widetilde{\mathcal F}_n = Span\mathcal F_n = Span\mathcal G_n = Span\mathcal H_n$,
\item  they are all of dimension $(n+1)^2$.
\end{itemize}
\end{lmm}
\begin{proof}
%we will proceed following two steps: first show that $Span\mathcal F_n = Span\mathcal G_n$, then show that $\dim \mathcal G_n=(n+1)^2$. 
{\bf Step 1.}
One clearly sees that $\widetilde{\mathcal F}_n \subset {\mathcal F}_n$, we show that every element of ${\mathcal F}_n$ belongs to $Span\widetilde{\mathcal F}_n$.
Considering any function $f_i$ in $\mathcal F_n$, according to $i_1$ being even or odd,  we will treat the two cases as follows.

%\begin{itemize}
%\item 
If $i_1$ is even, since $(\sin\varphi )^{i_1}=(1-\cos^2\varphi) ^{i_1/2}$, then we have:
$$
f_i(\theta,\varphi)
 = \sum_{k_\varphi = 0}^{i_1/2}
 \begin{pmatrix} i_1/2\\k_\varphi \end{pmatrix}
 (-1)^{k_\varphi}
  (\cos\varphi)^{2k_\varphi+i_2}(\sin\theta)^{i_1+i_2}(\cos\theta)^{i_3}/i!,
$$
and since $(\sin\theta)^{i_1+i_2} = (\sin\theta)^{2k_\varphi+i_2}(1-\cos^2\theta)^{i_1/2-k_\varphi}$, we can write:
$$
f_i(\theta,\varphi)
 = \sum_{k_\varphi = 0}^{i_1/2} \sum_{k_\theta = 0}^{i_1/2-k_\varphi}
 \begin{pmatrix} i_1/2\\k_\varphi \end{pmatrix}
 \begin{pmatrix} i_1/2-k_\varphi\\k_\theta \end{pmatrix}
 (-1)^{k_\varphi+k_\theta}
\sin^0\varphi  (\cos\varphi)^{2k_\varphi+i_2}(\sin\theta)^{2k_\varphi+i_2}(\cos\theta)^{2k_\theta+i_3}/i!.
$$
We then easily see that $0+(2k_\varphi+i_2) =2k_\varphi+i_2 $ while $(2k_\varphi+i_2)+(2k_\theta+i_3)\leq 2k_\varphi+i_2+i_1-2k_\varphi+i_3=|i|$, hence in particular $f_i\in Span\widetilde{\mathcal F}_n$.

%\item 
If $i_1$ is odd, since  $(\sin\varphi )^{i_1}=\sin\varphi(1-\cos^2\varphi) ^{(i_1-1)/2}$, then we have:
$$
f_i(\theta,\varphi)
 = \sum_{k_\varphi = 0}^{(i_1-1)/2}
 \begin{pmatrix} (i_1-1)/2\\k_\varphi \end{pmatrix}
 (-1)^{k_\varphi}
\sin\varphi  (\cos\varphi)^{2k_\varphi+i_2}(\sin\theta)^{i_1+i_2}(\cos\theta)^{i_3}/i!,
$$
and since $(\sin\theta)^{i_1+i_2} = (\sin\theta)^{2k_\varphi+1+i_2}(1-\cos^2\theta)^{(i_1-1)/2-k_\varphi}$, we can write:
$$
f_i(\theta,\varphi)
 = \sum_{k_\varphi = 0}^{\frac{i_1-1}2} \sum_{k_\theta = 0}^{\frac{i_1-1}2-k_\varphi}
 \begin{pmatrix} \frac{i_1-1}2\\k_\varphi \end{pmatrix}
 \begin{pmatrix} \frac{i_1-1}2-k_\varphi\\k_\theta \end{pmatrix}
 (-1)^{k_\varphi+k_\theta}
\sin\varphi  (\cos\varphi)^{2k_\varphi+i_2}(\sin\theta)^{2k_\varphi+1+i_2}(\cos\theta)^{2k_\theta+i_3}/i!.
$$
We then easily verify that $1+(2k_\varphi+i_2) =2k_\varphi+1+i_2 $ while $(2k_\varphi+1+i_2)+(2k_\theta+i_3)\leq 2k_\varphi+1+i_2+i_1-1-2k_\varphi+i_3=|i|$, hence again $f_i\in Span\widetilde{\mathcal F}_n$.
%\end{itemize}

This proves that $Span\ {\mathcal F}_n = Span\ \widetilde{\mathcal F}_n$.

{\bf Step 2.}
Let's start by considering any element $f_i$ of $\mathcal F_n$ to show that it belongs to $ Span\mathcal G_n$.
For all $i\in(\mathbb N_0)^3$ such that $|i|\leq n$, writing $\cos\varphi,\sin\varphi$ under their exponential form, and expanding their powers according to the binomial formula, we have:
%{\color{red} tu m'aideras a decider le bon niveau de d\'etail pour la calcul qui suit}
$$
f_i(\theta,\varphi)
= 
\frac1{i!(2\mi)^{i_1}2^{i_2}}(\sin\theta)^{i_1+i_2}\left(
\sum_{j_1=0}^{i_1}
\sum_{j_2=0}^{i_2}
(-1)^{j_1}
\begin{pmatrix} i_1\\j_1\end{pmatrix}
\begin{pmatrix} i_2\\j_2\end{pmatrix}
\exp\big( \mi (i_1-2j_1+i_2-2j_2)\varphi\big)
\right)
(\cos\theta)^{i_3}.
$$
In order to write each 
$(\sin\theta)^{i_1+i_2}
\exp \big(\mi (i_1-2j_1+i_2-2j_2)\varphi\big)
(\cos\theta)^{i_3}$ 
term in this sum as a linear combination of $Z_l^m$s, the power of $\sin\theta$ must be the absolute value of the power of $\exp \mi \varphi$.
Therefore, defining $M_{i,j_1,j_2}:=\min( i_1-j_1+i_2-j_2,j_1+j_2)$, we can write $i_1+i_2$ independently of the sign of $i_1-2j_1+i_2-2j_2$ as:
$$
i_1+i_2 = |i_1-2j_1+i_2-2j_2| +2 M_{i,j_1,j_2}.
$$
This leads to
\begin{equation*}
\begin{array}{l}
(\sin\theta)^{i_1+i_2}
\exp \big(\mi (i_1-2j_1+i_2-2j_2)\varphi\big)
(\cos\theta)^{i_3}\\\displaystyle \phantom{(\sin)}
=(\sin\theta)^{ |i_1-2j_1+i_2-2j_2|}
(1-\cos^2\theta)^{M_{i,j_1,j_2}}
\exp\big( \mi (i_1-2j_1+i_2-2j_2)\varphi\big)
(\cos\theta)^{i_3}
\\\displaystyle \phantom{(\sin)}
=
\sum_{j_3=0}^{M_{i,j_1,j_2}}
\begin{pmatrix} M_{i,j_1,j_2}\\j_3\end{pmatrix}
(-1)^{j_3}
(\sin\theta)^{|i_1-2j_1+i_2-2j_2|}
(\cos\theta)^{2j_3+i_3}
\exp \big(\mi (i_1-2j_1+i_2-2j_2)\varphi\big),
\end{array}
\end{equation*}
so
\begin{equation*}
f_i(\theta,\varphi)
= 
\frac1{i!(2\mi)^{i_1}2^{i_2}}
\sum_{j_1=0}^{i_1}
\sum_{j_2=0}^{i_2}
(-1)^{j_1}
\begin{pmatrix} i_1\\j_1\end{pmatrix}
\begin{pmatrix} i_2\\j_2\end{pmatrix}
\sum_{j_3=0}^{M_{i,j_1,j_2}}
\begin{pmatrix} M_{i,j_1,j_2}\\j_3\end{pmatrix}
(-1)^{j_3}
Z_{|i_1-2j_1+i_2-2j_2|+2j_3+i_3}^{|i_1-2j_1+i_2-2j_2|} (\theta,\varphi).
\end{equation*}
To verify that these $Z_l^m$s belong to the space $Span \mathcal G_n$, we must make sure that their indices satisfy  $0\leq l\leq n$ and $|m|\leq l$. In fact:
\begin{itemize}
\item if $i_1-2j_1+i_2-2j_2\geq 0$ then
$$
\begin{array}{rl}
j_3\leq j_1+j_2&\Rightarrow i_1-2j_1+i_2-2j_2+2j_3+i_3\leq i_1+i_2+i_3,\\
i_1-2j_1+i_2-2j_2\geq 0&\Rightarrow i_1-2j_1+i_2-2j_2+2j_3+i_3\geq 0,\\
 i_1-2j_1+i_2-2j_2= |i_1-2j_1+i_2-2j_2|&\Rightarrow |i_1-2j_1+i_2-2j_2|\leq  i_1-2j_1+i_2-2j_2+2j_3+i_3.
\end{array}
$$
\item if $i_1-2j_1+i_2-2j_2 < 0$ then
$$
\begin{array}{rl}
j_3\leq i_1-j_1+i_2-j_2 &\Rightarrow 2j_1-i_1+2j_2-i_2+2j_3+i_3\leq i_1+i_2+i_3,\\
2j_1-i_1+2j_2-i_2\geq 0&\Rightarrow 2j_1-i_1+2j_2-i_2+2j_3+i_3\geq 0,\\
2j_1-i_1+2j_2-i_2= |i_1-2j_1+i_2-2j_2|&\Rightarrow |2j_1-i_1+2j_2-i_2|\leq  2j_1-i_1+2j_2-i_2+2j_3+i_3.
\end{array}
$$
\end{itemize}
To summarize, since  $i_1+i_2+i_3=|i|$, the indices of each $Z_{|i_1-2j_1+i_2-2j_2|+2j_3+i_3}^{|i_1-2j_1+i_2-2j_2|} $ are such that:
$$\left\{
\begin{array}{l}
| i_1-2j_1+i_2-2j_2|+2j_3+i_3\leq |i|,\\
| i_1-2j_1+i_2-2j_2|+2j_3+i_3 \geq 0,\\
 |i_1-2j_1+i_2-2j_2|\leq  | i_1-2j_1+i_2-2j_2|+2j_3+i_3,
\end{array}\right.
$$
Hence, for $|i|\leq n$, $f_i$ can indeed be written as a linear combination of $Z_{l}^{m}$ with $0\leq l\leq n$ and $|m|\leq l$. In other words,  each element of $\mathcal F_n$ belongs to $Span\mathcal G_n$

Besides, any element $Z_l^{\pm n}$ of $\mathcal G_n$ can also be written as a linear combination of $f_i$s. In fact, expanding $e^{\mi m\varphi}$ via the binomial formula we immediately see that for all $(l,m)\in\mathbb Z^2$ such that $0\leq l \leq n$ and $-l\leq m\leq l$ we have:
$$
Z_l^{\pm m}(\theta,\varphi) = \sum_{m' = 0}^{|m|} \begin{pmatrix}{|m|}\\m'\end{pmatrix} (\cos\varphi)^{m'}(\pm \mi\sin\varphi)^{|m|-m'}(\sin\theta)^{|m|}(\cos\theta)^{l-|m|},
$$
or equivalently
$$
Z_l^{\pm m}(\theta,\varphi) 
= \sum_{m' = 0}^{|m|}  (\pm \mi )^{|m|-m'} 
(|m|)!(l-|m|)!f_{\left(|m|-m',m',l-|m|\right)},
%= \sum_{m' = 0}^{|m|} \begin{pmatrix}{|m|}\\m'\end{pmatrix} (\pm \mi )^{|m|-m'} 
%(|m|-m')!(m')!(l-|m|)!f_{\left(|m|-m',m',l-|m|\right)},
\text{ where } \left|\left(|m|-m',m',l-|m|\right)\right| = l\leq n,
$$
hence each element of $\mathcal G_n$ belongs to $Span\mathcal F_n$. 

As a consequence, $Span\mathcal G_n=Span\mathcal F_n$.

{\bf Step 3.}
Any function in $\mathcal H_n$ can be written  
as 
$$
\begin{array}{rl}
\text{if }m\geq0,\
Y_l^m(\theta,\phi) 
&\displaystyle
=C_l^m \frac{(-1)^m}{2^ll!} (\sin\theta)^m e^{\mi m\phi}
\sum_{l'=\left\lceil\frac{l+m}2\right\rceil}^l 
\begin{pmatrix} l\\l'\end{pmatrix} (-1)^{l-l'} \frac{(2l')!}{(2l'-l-m)!}(\cos\theta)^{2l'-l-m}
\\&\displaystyle
=C_l^m \frac{(-1)^m}{2^ll!}
\sum_{l'=\left\lceil\frac{l+m}2\right\rceil}^l 
\begin{pmatrix} l\\l'\end{pmatrix} (-1)^{l-l'} \frac{(2l')!}{(2l'-l-m)!} Z_{2l'-l}^m
\end{array}
$$
and 
$$
\begin{array}{rl}
\text{if }m<0,\
Y_l^m(\theta,\phi) 
&\displaystyle
=C_l^m  \frac{(l-|m|)!}{(l+|m|)!}\frac{1}{2^ll!} (\sin\theta)^{|m|} e^{\mi m\phi}
\sum_{l'=\left\lceil\frac{l+m}2\right\rceil}^l 
\begin{pmatrix} l\\l'\end{pmatrix} (-1)^{l-l'} \frac{(2l')!}{(2l'-l-m)!}(\cos\theta)^{2l'-l-|m|}
\\&\displaystyle
=C_l^m  \frac{(l-|m|)!}{(l+|m|)!}\frac{1}{2^ll!} 
\sum_{l'=\left\lceil\frac{l+m}2\right\rceil}^l 
\begin{pmatrix} l\\l'\end{pmatrix} (-1)^{l-l'} \frac{(2l')!}{(2l'-l-m)!}Z_{2l'-l}^m
\end{array}
$$
hence each element of $\mathcal H_n$ belongs to $Span\mathcal G_n$. 
So $Span\mathcal H_n\subset Span\mathcal G_n$.

Moreover, by property of the spherical harmonics, these are linearly independent, therefore
$$
\dim \mathcal H_n = (n+1)^2.
$$

We have then shown that
$$
(n+1)^2\leq \dim Span\mathcal G_n.
$$
But the space $\mathcal G_n$ has $(n+1)^2$ elements, so $\dim Span\mathcal G_n=(n+1)^2$.

{\bf Conclusion}
The result follows trivially from combining the previous steps.
\end{proof}
\begin{rmk}
The spaces $ \mathcal F_n $ and $ \mathcal H_n $ are the spaces of traces on the unit sphere respectively of homogeneous polynomials of degree at most equal to $n$ and  of harmonic homogeneous polynomials of degree at most equal to $n$.
The fact that they coincide follows from a more general result from harmonic polynomials theory,
see for instance \cite{waveletsref}.
\end{rmk}

As direct consequences of \Cref{lmm:dimFn,lmm:matsfuncsevals} as well as \Cref{cor:rkmatexp,cor:submatexp}, we get the following.
\begin{cor}
\label{cor:maxrk}
Given $(n,p)\in\mathbb N^2$,
 for any choice of $p$ directions $\{ \mathbf d_l\}_{1\leq l\leq p}$3,
the rank of $\mate^{n,p}$, and hence the rank of $\matr^{n,p}$, cannot be larger than $(n+1)^2$. 
%\end{cor}
%\begin{cor}
Moreover the submatrices of $\mate^{n,p}$ and $\matr^{n,p}$ corresponding to removing the rows $i$ such that $i_1>1$, denoted respectively $\overline{\mate}^{n,p} $ and $\overline{\matr}^{n,p} $, have the same rank as $\matr^{n,p}$ as well.
\end{cor}
The next natural question is that of existence of a set of $(n+1)^2$ directions $\{ (\theta_l,\varphi_l)\in\mathbb R^2\}_{1\leq l\leq (n+1)^2}$  such that the corresponding matrices are of rank $(n+1)^2$.
Hence it is natural to now fix $p=(n+1)^2$, denoting hereafter the corresponding matrices $\mate^{[n]}$ and $\matr^{[n]}$. If such a set of directions exists, then increasing the value of $p$  will not increase the rank according to \Cref{cor:maxrk}. 

\begin{theorem}
\label{thm:matmaxrk}
Given $n\in\mathbb N$, let $p=(n+1)^2$ directions on $\mathbb S^2$ be chosen as:
$$
\mathbf d_{l,m} = (\sin\theta_l\cos\varphi_{lm},\sin\theta_l\sin\varphi_{lm},\cos\theta_l)
$$
for all $l$ from $0$ to $n$ with $|m|\leq l$, where the $n+1$ colatitude angles $\{\theta_l\}_{0\leq l\leq n}\subset (0,\pi)$ are all different from each other, and the azimuths $\{\varphi_{lm}\}_{0\leq l\leq n,|m|\leq l}\subset [0,2\pi)$ satisfy $\varphi_{lm}\neq\varphi_{lm'}$ for every $m\neq m'$.
Then the two  $\frac{(n+1)(n+2)(n+3)}{6} \times (n+1)^2$ matrices $\matr^{[n]} $ and $\mate^{[n]} $  are of rank $(n+1)^2$.

Moreover, their square sub-matrices of $\matr^{[n]} $ and $\mate^{[n]} $ corresponding to removing the rows $i$ such that $i_1>1$, denoted respectively  $\overline{\matr}^{[n]} $ and $\overline{\mate}^{[n]} $, are also of rank $(n+1)^2$.
%
%Moreover, $\mate^{[n]} $ {\color{red} avec les m\^emes angles} is also of rank $(n+1)^2$.
\end{theorem}
Here again, properties of spherical harmonics will be at the center of the proof.
\begin{proof}
%start by showing that the rank of the ref matrix is that of the spherical harm mat
Since $\mathcal H_n$ is a generating set for $Span\ \mathcal F_n$, then according to \Cref{lmm:matsfuncsevals} the rank of $\matr^{[n]}$ is that of the $(n+1)^2\times(n+1)^2$ matrix $\mats_n$  defined element wise by evaluating the elements of the function space $\mathcal H_n$ at a set of $p=(n+1)^2$ points $\{(\theta_l,\varphi_{lm})\}_{0\leq l\leq n,|m|\leq l}\subset (0,\pi)\times[0,2\pi)$. 
As stated in Lemmas 3.4.1 and 3.4.2 from \cite{moiolaTH}, this matrix $\mats_n$ is invertible. As a conclusion  $\matr^{[n]} $ is indeed of rank $(n+1)^2$. As a direct consequence of \Cref{cor:rkmatexp}, $\mate^{[n]} $ is also of rank $(n+1)^2$.

Since $\overline{\mathcal F}_n$ is a generating set for $Span\ \mathcal F_n$ from \Cref{lmm:dimFn}, applying the second statement from \Cref{lmm:matsfuncsevals} then proves the second statement.
\end{proof}

\begin{rmk}
In terms of directions, the result is actually true for almost any set of points $\{ (\theta_l,\varphi_l)\in(0,\pi)\times[0,2\pi)\}_{0\leq l\leq n, |m|\leq l}$ as Lemma 3.4.1 from \cite{moiolaTH}  states that the rank of the matrix $\mats_n$ is $(n+1)^2$ for a dense open set of $(\mathbb S^2)^p$.
\end{rmk}

%%%%%%%%%%%%%%%%%%%%%%%%%%%%%%%%%%%
\subsection{Relating GPW and reference matrices}\label{ssec:Relate}

We are now interested in the rank of the GPW matrices, namely the $(n+1)(n+2)(n+3)/6\times p$ matrices $\mata$ and $\matp$.
In order to leverage the properties of the reference matrix to prove properties of these GPW matrices,
the missing link is then to understand their relation to the reference matrix.
%In other words, it is then natural to study the derivatives of each quasi-Trefftz function in terms of the corresponding $(\sin\varphi) ^{i_1} (\cos\varphi)^{i_2}(\sin\theta)^{i_1+i_2}(\cos\theta)^{i_3}$ term,

Because of their polynomial component, either in the phase, the amplitude, or the function itself, GPW functions have derivatives that share a common structure when evaluated at $\pG$ as long as their order is not higher than $q+1$. Even though higher order derivatives could also be studied, they would not share this common structure.
Hence we will start by studying the common properties of such derivatives of GPW functions.
Under the assumption that $q\geq n-1$, the entries in the GPW matrices are precisely such derivatives evaluated at a $\pG$. This will then allow us to establish a relation between the reference matrix and each of the GPW matrices.

Noticeably, neither the number of GPW functions chosen to construct each matrix, denoted $p$, nor the set of directions in the initialization come into play in this procedure: we establish relations between matrices independently of both. However, proving that the GPW  functions can be constructed to guarantee that the corresponding GPW matrices have maximal rank will rely on an appropriate choices for $p$ and the set of directions.

%\subsubsection{The GPW case}
To address the two GPW cases, it is natural to start from expressing the derivatives of GPW functions in terms of the  three non-zero free parameters in the initialization, leveraging \Cref{prop:Pblambdas,prop:Abmus}.
%\subsubsection{For the amplitude-based GPWs}
\begin{prop}
\label{prop:AbME}%Matrix Entries
Given $q\in\mathbb N$ and a point $\pG\in\mathbb R^3$, a set of complex-valued functions  $c = \{\funcci,i\in\mathbb N_0^3, |i|\leq 2\}$ is assumed to satisfy \Cref{hyp:PDop}.

Consider any amplitude-based GPW associated to partial differential operator $\Lc$, $J(\varx) := \polQ(\varx-\pG)\exp \Lambda \cdot\Big(\varx-\pG \Big)$ with $\displaystyle Q:=\sum_{i\in\mathbb N_0^3, |i|\leq q+1} \mu_i \mathbf X^i $, constructed via \Cref{Algo:AbS,Algo:SubS}, with the initialization introduced in \Cref{sec:norm} for $\cst\in\mathbb C^*$ and a unit vector $\dir\in\mathbb S^2$.
Then, for all $j\in(\mathbb N_0)^3$ such that $|j|\leq q+1$, the difference $\partial_x^j J(\pG)/j!-(\Lambda_{1})^{j_1}(\Lambda_{2})^{j_2}(\Lambda_{3})^{j_3}/j!$ can be expressed as a polynomial  in $\mathbb C[\Lambda_{1},\Lambda_{2},\Lambda_{3}]$, with degree smaller than $|j|$ and coefficients  depending on $\cst$ yet independent of $\dir$.
\end{prop}
\begin{proof}
Since $J(\varx) := \polQ(\varx-\pG)\exp \Lambda \cdot\Big(\varx-\pG \Big)$, we can express $\partial_\varx^j J$ thanks to Leibniz's rule as:
$$\forall x\in\mathbb R^3, 
\partial_\varx^j J(\varx) = \sum_{\widetilde j\in\mathbb N_0^3;\widetilde j\leq j} 
\begin{pmatrix}j\\\widetilde j\end{pmatrix}
\Lambda^{j-\widetilde j} \partial_\varx^{\widetilde j}\polQ(\varx-\pG),
$$
and therefore
$$
\partial_\varx^j J(\pG) = \sum_{\widetilde j\in\mathbb N_0^3;\widetilde j\leq j} 
\frac{j!}{(j-\widetilde j)!}
\Lambda^{j-\widetilde j} \mu_{\widetilde j}.
$$
Let's consider the terms in this linear combination, starting from considering the possible degree of each individual term as a polynomial in $\mathbb C[\Lambda_{1},\Lambda_{2},\Lambda_{3}]$:
\begin{enumerate}
\item for $|\widetilde j|=0$, then $\Lambda^{j-\widetilde j} \mu_{\widetilde j}=(\Lambda_{1})^{j_1}(\Lambda_{2})^{j_2}(\Lambda_{3})^{j_3}$ is a polynomial of degree $|j|$,
\item for $|\widetilde j|=1$, then $\Lambda^{j-\widetilde j} \mu_{\widetilde j}=0$,
\item for $|\widetilde j|>1$, then $\Lambda^{j-\widetilde j} \mu_{\widetilde j}$  can be expressed as a polynomial of degree at most equal to $|j-\widetilde j|+|\widetilde j|-1<|j|$
 according \Cref{prop:Abmus}.
\end{enumerate}
In the first case, the corresponding weight in the linear combination  is precisely $j!/j!=1$, while all other terms can be expressed as a polynomial in $\mathbb C[\Lambda_{1},\Lambda_{2},\Lambda_{3}]$ with total degree smaller than $|j|$. This concludes the proof.

These polynomials have coefficients independent of $\dir$ as a consequence of 
\Cref{prop:Abmus}.
\end{proof}

%\subsubsection{For the phase-based GPWs}
\begin{prop}
\label{prop:PbME}%Matrix Entries
Given $q\in\mathbb N$ and a point $\pG\in\mathbb R^3$, a set of complex-valued functions  $c = \{\funcci,i\in\mathbb N_0^3, |i|\leq 2\}$ is assumed to satisfy \Cref{hyp:PDop}.

Consider any phase-based GPW associated to differential operator $\Lc$, $G(\varx) :=\exp  \polP(\varx-\pG)$ with $\displaystyle P:=\sum_{i\in\mathbb N_0^3, |i|\leq q+1} \lambda_i \mathbf X^i $, constructed via \Cref{Algo:PbS,Algo:SubS}, with the initialization introduced in \Cref{sec:norm} for $\cst\in\mathbb C^*$ and a unit vector $\dir\in\mathbb S^2$.
Then, for all $j\in(\mathbb N_0)^3$ such that $|j|\leq q+1$, the difference $\partial_x^j G(\pG)/j!-(\lambda_{e_1})^{j_1}(\lambda_{e_2})^{j_2}(\lambda_{e_3})^{j_3}/j!$ can be expressed as a polynomial  in $\mathbb C[\lambda_{e_1},\lambda_{e_2},\lambda_{e_3}]$, with degree smaller than $|j|$ and coefficients depending on $\cst$ yet independent of $\dir$.
\end{prop}
\begin{proof}
Since $G(\varx)=\exp  \polP(\varx-\pG)$, we can express $\partial_\varx^j G$ thanks to Faa di Bruno's formula (see Appendix \ref{App:FdB}) as:
$$\forall x\in\mathbb R^3, 
\partial_\varx^j G(\varx) = 
\sum_{1\leq \widetilde m\leq |j|}
\exp(P(\varx-\pG)) \sum_{s=1}^{|j|}\sum_{p_s(j,\widetilde m)} j!\prod_{m=1}^s 
\frac{1}{k_m!}
\left(\frac{1}{l_m!} \partial_x^{l_m} P(\varx-\pG) \right)^{k_m},
$$
where sets $p_s$ as well as indices $l_m$ and $k_m$ are as defined in \Cref{App:FdB},
and therefore
$$
\partial_\varx^j G(\pG) = 
\sum_{1\leq \widetilde m\leq |j|}
 \sum_{s=1}^{|j|}\sum_{p_s(j,\widetilde m)} j!\prod_{m=1}^s 
\frac{1}{k_m!}
\left(\lambda_{l_m} \right)^{k_m}.
$$
%In this linear combination of products of $\lambda$s, the term $(\lambda_{e_1})^{j_1}(\lambda_{e_2})^{j_2}(\lambda_{e_3})^{j_3}$ corresponds to $s=3$, $k_m=j_m$ while $l_m=e_m$ for $m$ from $1$ to $3$, with a weight of $\frac{j!}{j_1!j_2!j_3!} = 1$. 
Let's consider the terms in this linear combination, starting from considering the possible degree of each individual $\lambda_{l_m}$ term as a polynomial in $\mathbb C[\lambda_{e_1},\lambda_{e_2},\lambda_{e_3}]$:
\begin{enumerate}
\item  if $|l_m|=1$, then  $\lambda_{l_m}$ is a polynomial of degree $1$,
\item  if $|l_m|>1$, then  $\lambda_{l_m}$ can be expressed as a polynomial of degree $\leq|l_m|-1$ according \Cref{prop:Pblambdas}.
\end{enumerate}
%starting from ...ing the partitions $p_s(i,\widetilde m)$.
Thus each $\prod_{m=1}^s \left(\lambda_{l_m} \right)^{k_m}$ can be expressed as a polynomial in $\mathbb C[\lambda_{e_1},\lambda_{e_2},\lambda_{e_3}]$ with total degree at most equal to:
$$
\sum_{|l_m| = 1}k_m +\sum_{|l_m|>1} k_m (|l_m|-1).
$$
Each partition of $j$ either does not or does contain any $l_m$ such that $|l_m|\neq 1$.
Accordingly,  each $\prod_{m=1}^s \left(\lambda_{l_m} \right)^{k_m}$ falls into one of the two following categories, it can be expressed as a polynomial with total degree:
\begin{enumerate}
\item either equal to $\sum_{m=1}^s k_m=\widetilde m$, when the partition contains only $l_m$ such that $|l_m|= 1$;
\item or smaller than $\sum_{m=1}^s k_m=\widetilde m\leq |j|$, when the partition contains at least one $l_m$ such that $|l_m|\neq 1$.
\end{enumerate}
In the first case,  each $l_m$ belongs to $\{e_1,e_2,e_3\}$, and, since the partition must satisfy $\sum_{m=1}^s k_m l_m =j$, it corresponds to $s=3$ with $j=j_1e_1+j_2e_2+j_3e_3$. Hence $(\lambda_{e_1})^{j_1}(\lambda_{e_2})^{j_2}(\lambda_{e_3})^{j_3}$ is precisely of degree equal to $|j|$, with a weight of $\frac{j!}{j_1!j_2!j_3!} = 1$, whereas all other $\prod_{m=1}^s \left(\lambda_{l_m} \right)^{k_m}$ terms can be expressed as a polynomial in $\mathbb C[\lambda_{e_1},\lambda_{e_2},\lambda_{e_3}]$ with total degree smaller than $|j|$. This concludes the proof.

These polynomials have coefficients independent of $\dir$ as a consequence of \Cref{prop:Pblambdas}.
\end{proof}

%\subsubsection{Common GPW property}
Since \Cref{prop:AbME,prop:PbME} state similar relations between the derivatives of GPW functions on the one hand and the three initialization parameters on the other hand, we can now prove the common property of both families of GPWs.
\begin{prop}
Given $(n,p)\in\mathbb N^2$ and a point $\pG\in\mathbb R^3$, a set of complex-valued functions  $c = \{\funcci,i\in\mathbb N_0^3, |i|\leq 2\}$ is assumed to satisfy \Cref{hyp:PDop}.

Consider any sets of phase-based and amplitude-based GPWs associated to partial differential operator $\Lc$, constructed via \Cref{Algo:PbS,Algo:AbS,Algo:SubS}, with $q=\max(n-1,1)$ and the initialization introduced in \Cref{sec:norm} for $\cst\in\mathbb C^*$ and any set of unit vectors $\{\dir_l\in\mathbb S^2,1\leq l\leq p\}$.
Then there exist square lower triangular matrices $\lata^n$ and $\latp^n$, whose diagonal coefficients are equal to $1$ and whose other non-zero coefficients depend only on (i) derivatives of the PDE coefficients $c$ evaluated at $\pG$ and (ii) the initialization parameter $\cst$, such that
\begin{equation*}
\mata^{n,p} = \lata^n \mate^{n,p} \text{ and } \matp^{n,p} = \latp^n \mate^{n,p}.
\end{equation*}
As a result, we have: 
%{\color{red} under the current form this statement is indep of the number of columns of these matrices $p$ - et c'est fait expr\`es}
\begin{equation*}
rk\left(\mata^{n,p} \right)= rk\left( \mate^{n,p}\right) \text{ and } rk\left( \matp^{n,p}\right) = rk\left( \mate^{n,p}\right).
\end{equation*}
\end{prop}
\begin{proof}
This is a direct consequence of \Cref{prop:AbME,prop:PbME} and the choice of numbering for the matrix entries. Indeed:
\begin{itemize}
\item as a consequence of \Cref{prop:AbME,prop:PbME}, for all $j$ such that $|j|\leq q$ there exist coefficients denoted $l^A_{ij} $ and $l^P_{ij}$ for $i\in\mathbb N_0^3$ with $|i|<|j|$ satisfying for any vector $\dir$ and the corresponding amplitude- or phase-based GPW function, denoted respectively  $J$ or $G$, for $|j|>0$:
$$
\partial_x^j J(\pG)=(\Lambda_{1})^{j_1}(\Lambda_{2})^{j_2}(\Lambda_{3})^{j_3}
+
\sum_{|i|< |j|}l^A_{ij} \frac{j!}{i!}(\Lambda_{1})^{i_3}(\Lambda_{2})^{i_2}(\Lambda_{3})^{i_3}
$$
$$
\partial_\varx^j G(\pG)= (\lambda_{e_1})^{j_1}(\lambda_{e_2})^{j_2}(\lambda_{e_3})^{j_3} 
+
\sum_{|i|< |j|}l^P_{ij} \frac{j!}{i!}(\lambda_{e_1})^{i_3}(\lambda_{e_2})^{i_2}(\lambda_{e_3})^{i_3}
$$
%$$
%\frac{\partial_\varx^j G(\pG)}{j!} = \frac{(\lambda_{e_1})^{j_1}(\lambda_{e_2})^{j_2}(\lambda_{e_3})^{j_3}}{j!} +
%\sum_{i\prec j}(L^P_n)_{\mathcal N(j),\mathcal N(i)}\frac{(\lambda_{e_1})^{i_3}(\lambda_{e_2})^{i_2}(\lambda_{e_3})^{i_3}}{i!}
%$$
while for $|j|=0$ both GPWs satisfy:
$$
J(\pG) = 1 \text{ and } G(\pG) = 1;
$$
\item we can then define two square matrices of size $(n+1)(n+2)(n+3)/6$, $\lata^n$ and $\latp^n$, by:
$$
 \left(\lata^n\right)_{\mathcal N(j),\mathcal N(i)}=
\left\{\begin{array}{l}
1 \text{ if } i=j,\\
 \frac{1}{i!}l^A_{ij} \text{ if } |i|<|j|,\\
0 \text{ otherwise},
\end{array}\right.
\text{ and similarly }
 \left(\latp^n\right)_{\mathcal N(j),\mathcal N(i)}=
\left\{\begin{array}{l}
1 \text{ if } i=j,\\
 \frac{1}{i!}l^P_{ij} \text{ if } |i|<|j| ,\\
0 \text{ otherwise},
\end{array}\right.
$$
both lower triangular matrices since $|i|<|j|$ implies that $\mathcal N(i)<\mathcal N(j)$,
\item
therefore, by definition of the matrices $\mata^{n,p}$, $\matp^{n,p}$ and $\matr^{n,p}$,  for $j\in\left(\mathbb N_0\right)^3\times \mathbb N$  with $|j|\leq n$ and $l\leq p$, we have the following relations:
$$
\mata^{n,p}_{\mathcal N(j),l } = \sum_{k=1}^{\frac{(n+1)(n+2)(n+3)}6} \left(\lata^n\right)_{\mathcal N(j),k}\matr^{n,p}_{ k,l}
\text{ and }
\matp^{n,p}_{\mathcal N(j),l} = \sum_{k=1}^{\frac{(n+1)(n+2)(n+3)}6} \left(\latp^n\right)_{\mathcal N(j),k}\matr^{n,p}_{k,l }
$$
\end{itemize}
\end{proof}
Fixing $p=(n+1)^2$, we denote hereafter the corresponding $(n+1)(n+2)(n+3)/6 \times (n+1)^2$ matrices $\mata^{[n]}:=\mata^{n,(n+1)^2/2}$ and $\matp^{[n]}:=\matp^{n,(n+1)^2/2}$.
As a direct consequence of \Cref{thm:matmaxrk}, we then obtain the following result.
\begin{cor}
\label{cor:GPWmats}
Given $n\in\mathbb N$ and a point $\pG\in\mathbb R^3$, a set of complex-valued functions  $c = \{\funcci,i\in\mathbb N_0^3, |i|\leq 2\}$ is assumed to satisfy \Cref{hyp:PDop}.

Let $p=(n+1)^2$ directions on $\mathbb S^2$ be chosen as:
$$
\mathbf d_{l,m} = (\sin\theta_l\cos\varphi_{lm},\sin\theta_l\sin\varphi_{lm},\cos\theta_l)
$$
for all $l$ from $0$ to $n$ with $|m|\leq l$, where the $n+1$ colatitude angles $\{\theta_l\}_{0\leq l\leq n}\subset (0,\pi)$ are all different from each other, and the azimuths $\{\varphi_{lm}\}_{0\leq l\leq n,|m|\leq l}\subset [0,2\pi)$ satisfy $\varphi_{lm}\neq\varphi_{lm'}$ for every $m\neq m'$.
Consider any sets of phase-based and amplitude-based GPWs associated to partial differential operator $\Lc$, constructed via \Cref{Algo:PbS,Algo:AbS,Algo:SubS}, with $q=\max(n-1,1)$ the initialization introduced in \Cref{sec:norm} for $\cst\in\mathbb C^*$ and directions $\{\dir_{l,m},0\leq l\leq n, |m|\leq l\}$.
Then the corresponding matrices $\mata^{[n]} $ and $\matp^{[n]} $, of size ${(n+1)(n+2)(n+3)}/{6} \times (n+1)^2$ are of rank $(n+1)^2$.
\end{cor}

\subsection{{The polynomial matrix}}
 
To address the polynomial case, we study the rank of the Taylor expansion matrix $\matq^{[n]}$.
 \begin{prop}\label{prop:polmat}
Given $(n,p)\in\mathbb N^2$ and a point $\pG\in\mathbb R^3$, a set of complex-valued functions  $c = \{\funcci,i\in\mathbb N_0^3, |i|\leq 2\}$ is assumed to satisfy \Cref{hyp:PDop}.

Consider the set of polynomial quasi-Trefftz functions associated to partial differential operator $\Lc$, constructed via \Cref{Algo:PolS,Algo:SubS}, with $q=\max(n-1,1)$ and  the initialization introduced in \Cref{sec:norm}. Then
$$
 rk \left(\matq^{[n]}\right) = (n+1)^2.
$$
\end{prop}
\begin{proof}
The matrix $\matq^{[n]}$ is of size ${(n+1)(n+2)(n+3)}/{6} \times (n+1)^2$, so in particular $ rk \left(\matq^{[n]}\right) \leq(n+1)^2$.

The set of polynomial quasi-Trefftz functions introduced in \Cref{sec:norm} are indexed by any $ j\in\left(\mathbb N_0\right)^3$, with $|j|\leq q+1$ and $j_1\in\{0,1\}$, so consider a numbering $\mathcal N_{pol}$ of the set of indices $\{ j\in\left(\mathbb N_0\right)^3,|j|\leq q+1,j_1\in\{0,1\}\}$. 
By definition of the matrix $\matq^{[n]}$ and according to the initialization, for all $(i_0,j_0)\in\{ j\in\left(\mathbb N_0\right)^3,|j|\leq q+1,j_1\in\{0,1\}\}^2$ we have:
$$
\left(\matq^{[n]}\right)_{\mathcal N(i_0),\mathcal N_{pol}(j_0)} 
%= \nu_{i_0}
=\delta(i_0-j_0).
$$
Hence the set of $(n+1)^2$ rows numbered $\mathcal N(i)$ for all $i\in\{ j\in\left(\mathbb N_0\right)^3,|j|\leq q+1,j_1\in\{0,1\}\}$ is clearly linearly independent, so $ rk \left(\matq^{[n]}\right) \geq(n+1)^2$. This concludes the proof.
\end{proof}

\subsection{Quasi-Trefftz bases Approximation properties}
%{\color{red} attention ou est defini cet espace - peut \^etre introduire l'id\'ee sans aucune notation math\'ematique pour \'eviter les r\'ep\'etitions}
We finally want to show that the three spaces spanned by the  quasi-Trefftz function sets introduced in \Cref{sec:QTfams}  %denoted here $\mathbb V_h^G$ 
satisfy the desired approximation property \eqref{IntPb}, 
%namely
%\begin{equation*}%\label{IntPb}
%\begin{array}{l}
%\forall u%\in \mathcal C^\infty 
%\text{ satisfying }\Lc u=0,
%\exists u_a\in\mathbb V_h^G \text{ s. t. }
%\forall \varx\in\mathbb R^2,
%|u(\varx)-u_a(\varx)|\leq C \| (\varx-\pG) \|^{n+1},
%\end{array}
%\end{equation*}
 via matching of $u_a$'s Taylor expansion to that of $u$ as announced in \eqref{eq:TEapp}.

\begin{theorem}\label{thm:approx}
Given $n\in\mathbb N$ and a point $\pG\in\mathbb R^3$, a set of complex-valued functions  $c = \{\funcci,i\in\mathbb N_0^3, |i|\leq 2\}$ is assumed to satisfy \Cref{hyp:PDop}. Let $\Omega$ be an open set in $\mathbb R^3$ such that $\pG\in\Omega$.

Consider the three quasi-Trefftz spaces  associated to partial differential operator $\Lc$, defined as the quasi-Trefftz spaces spanned by each of the three following sets:
\begin{itemize}
\item the set of amplitude-based GPWs, constructed via \Cref{Algo:AbS,Algo:SubS}, %with $q=n-1$, the initialization introduced in \Cref{sec:norm} for $\cst\in\mathbb C^*$ and $p=(n+1)^2$ directions $\{\dir_{l,m},0\leq l\leq n, |m|\leq l\}$ on $\mathbb S^2$,
\item the set of phase-based GPWs, constructed via \Cref{Algo:PbS,Algo:SubS}, %with $q=n-1$, the initialization introduced in \Cref{sec:norm} for $\cst\in\mathbb C^*$ and $p=(n+1)^2$ directions $\{\dir_{l,m},0\leq l\leq n, |m|\leq l\}$ on $\mathbb S^2$,
\item the set of polynomial functions, constructed via \Cref{Algo:PolS,Algo:SubS}, %with $q=n-1$, the initialization introduced in \Cref{sec:norm} for $\cst\in\mathbb C^*$ and $p=(n+1)^2$ directions $\{\dir_{l,m},0\leq l\leq n, |m|\leq l\}$ on $\mathbb S^2$. 
\end{itemize}
each of them constructed with $q=\max(n-1,1)$ and the initialization introduced in \Cref{sec:norm} for $\cst\in\mathbb C^*$ and $p=(n+1)^2$ directions %$\{\dir_{l,m},0\leq l\leq n, |m|\leq l\}$ 
on $\mathbb S^2$. 
As a side note, the polynomial space for these values of $q$ and $p$ is uniquely defined independently of the chosen initialization.
If the set of directions is chosen as:
$$
\forall(l,m)\in(\mathbb N_0)^2, l\leq n, |m|\leq l, 
\mathbf d_{l,m} = (\sin\theta_l\cos\varphi_{lm},\sin\theta_l\sin\varphi_{lm},\cos\theta_l),
$$
 where the $n+1$ colatitude angles $\{\theta_l\}_{0\leq l\leq n}\subset (0,\pi)$ are all different from each other, and the azimuths $\{\varphi_{lm}\}_{0\leq l\leq n,|m|\leq l}\subset [0,2\pi)$ satisfy $\varphi_{lm}\neq\varphi_{lm'}$ for every $m\neq m'$,
then any of these three spaces, denoted $\mathbb V^G_h$, satisfies the following approximation property:
\begin{equation}\label{eq:approx}
\begin{array}{l}
\forall u\in \mathcal C^{\max(2,n)}(\Omega)
\text{ satisfying }\Lc u=0,
\exists u_a\in\mathbb V_h^{G},\\
 \exists  C\in\mathbb R \text{ s. t. }
 \forall \varx\in\Omega,
\left\{ \begin{array}{l}
|u(\varx)-u_a(\varx)|\leq C \| \varx-\pG \|^{n+1},\\
\|\nabla u(\varx)-\nabla u_a(\varx)\|\leq C \| \varx-\pG \|^{n}.
\end{array}\right.
\end{array}
\end{equation}
The constant $C$ here depends on the desired order $n$, on the PDE solution $u$ in $\Omega$, as well as on the Taylor polynomials of the PDE coefficients $c$ evaluated at $\pG$.

It is to be noted that this
actually shows a convergence in the $H^1$ norm:
\begin{equation}
\forall u%\in \mathcal C^\infty 
\text{ satisfying }\mathcal Lu=0,
\exists u_a\in\mathbb V_h,\exists C\in\mathbb R, \text{ s. t. } \forall h \text{ small enough }
\|u-u_a\|_{H^1(B(\pG,h))}\leq C h^{n}.
\end{equation}
\end{theorem}
\begin{proof}
It is sufficient to prove that, if $\mathsf M^{[n]}$ is any of the three  $(n+1)(n+2)(n+3)/6 \times (n+1)^2$ quasi-Trefftz matrices, namely $\mata^{[n]}$, $\matp^{[n]}$ or $\matq^{[n]}$, then the linear system defined by:
$$
\mathsf M^{[n]} \mathsf x = \mathsf F_n
$$
has a solution for any $\mathsf F_n$ in the vector space $\Fspace$ defined as:
$$
\displaystyle
\Fspace :=\left\{\mathsf F\in\mathbb C^{(n+1)(n+2)(n+3)/6},\exists v\in\mathcal C^{\max(2,n)}(\Omega)\text{ s.t. } \Lc v=0  \text{ and }\forall i\in\left(\mathbb N_0\right)^3, |i|\leq n, 
\ \mathsf F_{\mathcal N(i)} = \partial_x^i v (\pG) /i!
\right\}.
$$
We first define a similar vector space by relaxing the Trefftz condition $\Lc v=0$ into a quasi-Trefftz condition:
$$
\begin{array}{l}
\displaystyle
\Kspace :=
\left\{
\mathsf K\in\mathbb C^{(n+1)(n+2)(n+3)/6},\exists v\in\mathcal C^{\max(2,n)}(\Omega)\text{ satisfying }\Lc v(\varx)=O\left(\|\varx-\pG\|^{n-1}\right)
\right.
\\\displaystyle
\phantom{
\mathfrak K :=
\mathfrak K :=
}
\left.
 \text{ s.t. }\forall i\in\left(\mathbb N_0\right)^3, |i|\leq n, \ \mathsf K_{\mathcal N(i)} = \partial_x^i v (\pG) /i! 
\right\}.
\end{array}
$$
It is then clear that $\Fspace\subset\Kspace$, but also that  the range of $\mathsf M^{[n]}$ is also included in $\Kspace$.

Next we want to show that the dimension of $\Kspace$ is $(n+1)^2$. For any function $f\in\mathcal C^{\max(2,n)}(\Omega)$ satisfying $\Lc f(\varx)=O(\|\varx-\pG\|^{n-1})$, by Leibniz rule we have:
$$
\forall i\in\left(\mathbb N_0\right)^3, |i|\leq n-2, 
\sum_{|j|\leq 2}\sum_{\tilde i\leq i} \begin{pmatrix}  i\\\tilde i \end{pmatrix}\partial_x^{i-\tilde i} c_j\partial_x^{j+\tilde i} f (\pG) = 0.
$$
This shows that $\Kspace$ is a subset of the kernel of an $(n-1)n(n+1)/6\times (n+1)(n+2)(n+3)/6$ matrix $\mathsf R$, with the following properties:
$$
\forall i\in\left(\mathbb N_0\right)^3, |i|\leq n-2, 
\left\{
\begin{array}{l}
\mathsf R_{\mathcal N(i)\mathcal N\left(i+2e_1\right)} = c_{2e_1}(\pG)% \text{ so } \mathsf R_{\mathcal N(i)\mathcal N\left(i+2e_1\right)} \neq 0 \text{ by \Cref{hyp:PDop}}
,\\
\mathsf R_{\mathcal N(i)\mathcal N(j)}=0 \text{ if } j>i+2e_1.
\end{array}
\right.
$$
Therefore, since $ c_{2e_1}(\pG) \neq 0$  by \Cref{hyp:PDop}, choosing a numbering scheme $\mathcal N$ such that if $|i|=|j|$ then $i\prec j$ implies $\mathcal N(i)\leq \mathcal N(j)$ highlights the echelon structure of $\mathsf R$.
As a result, the echelon structure of $\mathsf R$ guarantees that it has maximal rank, namely $(n-1)n(n+1)/6$, while its kernel is of dimension 
$$\frac{(n+1)(n+2)(n+3)}{6} -\frac{(n-1)n(n+1)}{6} = (n+1)^2. $$
Hence $\Kspace$ is a subset of a space of dimension $(n+1)^2$, but it also has a subset of the same dimension, namely the range of $\mathsf M^{[n]}$ according to \Cref{cor:GPWmats} and \Cref{prop:polmat}, so it is itself of dimension $(n+1)^2$.

This shows that the range of $\mathsf M^{[n]}$ and $\Kspace$ are the same space, and therefore any $\mathsf F\in\Fspace$ belongs to the range of $\mathsf M^{[n]}$. So this conclude the proof of the approximation of the function values.

Precisely because this proof relies on matching the Taylor expansions of $u$ and $u_a$, the result of approximation of the gradient can be obtained directly by taking derivatives of this Taylor expansion matching identity.
\end{proof}

\section{Numerical results}
We propose numerical experiments to illustrate the approximation properties presented in \Cref{thm:approx}, as well as the fact that the new polynomial quasi-Trefftz basis does not inherit the well-known conditioning issues of wave-like bases. These experiments include of course the implementation of the construction algorithms, but also the computation of the quasi-Trefftz approximation $u_a$ of the exact solution $u$ to a PDE.
Each test case is defined by 
a differential operator $\mathcal L$,
a domain $\Omega$,
an exact solution $u$, satisfying $\mathcal Lu=0$, to be approximated over $\Omega$.
The test cases are summarized below:
\begin{itemize}
    % ( (dx^2+dy^2+dz^2) + kappa^2*( 1-x ) ) u = 0
    % exact sol funcu=@(X) Ai( alpha x  ) * exp( i beta (y+z) )
                
    %Here we want (x,y,z)\in [-1,1]x[0,2 pi]x[-1,1]
\item (Test case 1) $\mathcal L = \Delta +\kappa^2$, 
	$\Omega =  [-1,1]\times[0,2 \pi]\times[-1,1]$,
	$u(x,y,z) = \exp i \kappa y$, and $\kappa = 3$;\\
	 for this constant-coefficient Helmholtz problem, plane waves are exact solution and performance of quasi-Trefftz bases can be compared to that of plane wave bases;
\item (Test case 2) $\mathcal L = \Delta +\kappa^2(1-x)$, 
	$\Omega = [-2,2]^3$, 
	$u(x,y,z) = Ai(\kappa^{2/3} x)\exp\left( i\kappa (y+z)/\sqrt{2}\right)$, and $\kappa = 2$;\\
	 for this variable-coefficient Helmholtz problem, the domain is chosen to include both a propagative zone ($\kappa^2(1-x)>0$) and an evanescent zone ($\kappa^2(1-x)<0$), as well as a smooth transition between them ($x=1$);
\item (Test case 3) $\mathcal L = \Delta - {M_0}^2\partial_z^2 +2i\kappa M_0 \partial_z +\kappa^2(1-x)$, 
	$\Omega =  [-2,2]^3$, \\
	$u(x,y,z) = Ai\left( \left( \frac{\kappa^2}{1-{M_0}^2} \right)^{1/3}z \right) \exp i\kappa \left( -\frac{M_0z}{1-{M_0}^2}  +\frac{x+y}{\sqrt{2-2{M_0}^2}}\right) $, $\kappa = 2$ and $M_0=0.2$;\\
	this convected Helmholtz problem introduces some anisotropy, both in the higher order terms and in the first order term of the governing differential operator.
\end{itemize}

Because the approximation properties are local, yet the quasi-Trefftz bases are aimed at discretizing a Galerkin formulation, instead of performing the tests at a single point $\pG$ we propose to do so at a set of random points in a given domain.
At each of $50$ random points in the domain, we follow the procedure described below.
\begin{enumerate}[label=\alph*)]
\item For each value of $n$ from $1$ to $8$, construct three quasi-Trefftz bases, two GPW bases and one polynomial basis, each of them with $p=(n+1)^2$ as the dimension of the basis and $q=\max(n-1,1)$ as the order of approximation of the Trefftz property.
For the initialization, the directions $\mathbf d_{l,m}$ are chosen with:
$$
\text{for } l \text{ from } 0 \text{ to } n, \theta_l= \frac{\pi}2+ \frac{(-1)^l \pi l}{2n+2} ,
\text{ and }
\text{for } m \text{ from } -l \text{ to } l, \varphi_{lm} = \frac{ 2m\pi}{2l+1}.
$$
\item For each basis we compute the linear combination of its elements $u_a$ as described in the proof of \Cref{thm:approx}, by solving the normal equation $\overline{\mathsf M^{[n]} }^T\mathsf M^{[n]} \mathsf x =\overline{\mathsf M^{[n]} }^T \mathsf F_n$. We are aware that this will raise to the square the condition number of the system and hence degrade the accuracy of the solutions, but this is not a concern to us because this study is mostly interested in validating the convergence orders, that will be observed {\it before} reaching the lowest errors.
\item For each quasi-Trefftz approximation $u_a$ we estimate the $L^\infty$ error between $u$ and $u_a$ over a ball centered at the random point and of radius $h$.
\end{enumerate}
We then compute for each value of $n$, for each basis, the worst error obtained at the $50$ random points.
Next we report the corresponding results, and clearly observe the convergence orders predicted by  \Cref{thm:approx}.

For the first test case, since the PDE is the constant-coefficient Helmholtz equation, plane wave functions with the appropriate wave number are exact solutions in this case. Hence we compare the results with a classical PW basis, which form here a set of quasi-Trefftz functions, for reference. 
\Cref{fig:HelmGPW,fig:HelmpqT} present respectively the results obtained from  the two GPW bases and the polynomial quasi-Trefftz basis. 
As anticipated, the GPW bases perform similarly to the PW basis: the expected orders of convergence are observed, they are the same, namely $n+1$, for a PW or GPW basis for a given number of elements, $p=2n+1$, and moreover the conditioning of the matrix  $\overline{\mathsf M^{[n]} }^T\mathsf M^{[n]} $ deteriorates as $n$ increases. This is a consequence of the well-known conditioning issues of wave-like bases.
As for the polynomial quasi-Trefftz basis,  the expected orders of convergence are observed as well, namely $n+1$, and match those for the PW basis. However, the conditioning of the polynomial quasi-Trefftz matrix $\overline{\mathsf M^{[n]} }^T\mathsf M^{[n]} $ increases much slower as $n$ increase. 
%{\color{red} check with Guillaume, what's the actual value of the conditioning in this case? can't we compute it?} 
This can be observed in the following table, indicating the approximate condition number of the matrices $\overline{\mathsf M^{[n]} }^T\mathsf M^{[n]} $ for the various bases depending on the value of $n$. 
\begin{center}
\begin{tabular}{|c||c|c|c|c|}
\hline
$n$ & plane wave basis & amplitude-based GPW basis & phase-based GPW basis &polynomial basis
\\\hline\hline
$n=1$ & $5.25\times 10^{1}$ & $5.25\times 10^{1}$ & $5.25\times 10^{1}$ & $1\times 10^0$\\
$n=2$ & $1.38\times 10^2$ & $1.38\times 10^2$ &  $1.38\times 10^2$ & $2.35\times 10^1$\\
$n=3$ & $5.51\times 10^3$ & $5.51\times 10^3$ & $5.51\times 10^3$ & $3.13\times 10^1$\\
$n=4$ & $4.28\times 10^5$ & $4.28\times 10^5$ & $4.28\times 10^5$ & $7.01\times 10^1$\\
$n=5$ & $8.18\times 10^7$ & $8.18\times 10^7$ & $8.18\times 10^7$ & $1.76\times 10^2$\\
$n=6$ & $3.10\times 10^{10}$ & $3.10\times 10^{10}$ & $3.10\times 10^{10}$ & $5.79\times 10^2$\\
$n=7$ & $2.47\times 10^{13}$ & $2.47\times 10^{13}$ & $2.47\times 10^{13}$ & $2.02\times 10^3$\\
$n=8$ & $1.14\times 10^{17}$ & $1.24\times 10^{17}$ & $2.24\times 10^{16}$ & $7.27\times 10^3$\\
\hline 
\end{tabular}
\end{center}
As a result, for the polynomial basis only, and even for increasing values of $n$, the linear combination of quasi-Trefftz functions $u_a$ can be computed up to machine precision, and so the error between $u$ and $u_a$ decreases until reaching machine precision. But for the other three bases $u_a$ cannot be computed up to machine precision.

For the second and third test cases, since the PDE has varying coefficients, plane wave functions are not exact solutions solutions anymore.  Hence we simply compare the performance of the three quasi-Trefftz bases. 
\Cref{fig:ConvHelm-1,fig:ConvHelm-2} present respectively the results obtained for the second and third test cases. 
Again as predicted by \Cref{thm:approx}, the expected orders of convergence are observed, namely $n+1$, for the three bases with $p=2n+1$ elements.
Similarly to the first test case,  the conditioning of the GPW matrices  $\overline{\mathsf M^{[n]} }^T\mathsf M^{[n]} $  deteriorates rapidly as $n$ increases, while that of the polynomial quasi-Trefftz matrix increases much slower. This can be observed in the two following tables, corresponding respectively to test cases 2 and 3.
\begin{center}
\begin{tabular}{|c||c|c|c|}
\hline
$n$ & amplitude-based GPW basis & phase-based GPW basis &polynomial basis
\\\hline\hline
$n=1$ & $9.55\times 10^4$ & $9.55\times 10^4$ & $1.00\times 10^0$ \\
$n=2$ & $1.12\times 10^{12}$ & $1.12\times 10^{12}$ & $3.54\times 10^1$\\
$n=3$ & $9.05\times 10^{15}$ & $9.05\times 10^{15}$ & $4.38\times 10^1$\\
$n=4$ & $5.78\times 10^{17}$ & $5.78\times 10^{17}$ & $9.68\times 10^1$\\
$n=5$ & $6.54\times 10^{17}$ & $6.54\times 10^{17}$ & $2.15\times 10^2$\\ 
$n=6$ & $9.14\times 10^{18}$ & $9.92\times 10^{18}$ & $6.63\times 10^2$\\
$n=7$ & $1.50\times 10^{19}$ & $6.62\times 10^{18}$ & $2.21\times 10^3$\\
$n=8$ & $4.08\times 10^{19}$ & $4.31\times 10^{19}$ & $7.77\times 10^3$\\
\hline 
\end{tabular}
\end{center}
\begin{center}
\begin{tabular}{|c||c|c|c|}
\hline
$n$ & amplitude-based GPW basis & phase-based GPW basis &polynomial basis
\\\hline\hline
$n=1$ & $2.64\times 10^2$ & $2.64\times 10^2$ & $1.00\times 10^0$\\
$n=2$ & $7.86\times 10^4$ & $7.86\times 10^4$ & $3.50\times 10^1$\\
$n=3$ & $5.25\times 10^8$ & $5.25\times 10^8$ & $5.37\times 10^1$\\
$n=4$ & $4.29\times 10^{12}$ & $4.29\times 10^{12}$ & $1.01\times 10^2$\\
$n=5$ & $3.51\times 10^{17}$ & $3.51\times 10^{17}$ & $2.29\times 10^2$\\
$n=6$ & $5.12\times 10^{17}$ & $7.66\times 10^{17}$ & $6.54\times 10^2$\\
$n=7$ & $9.85\times 10^{18}$ & $2.22\times 10^{18}$ & $2.20\times 10^3$\\
$n=8$ & $1.22\times 10^{19}$ & $1.38\times 10^{19}$ & $7.72\times 10^3$
\\\hline 
\end{tabular}
\end{center}
Here again, for the polynomial basis only, and even for increasing values of $n$, the linear combination of quasi-Trefftz functions $u_a$ can be computed up to machine precision, and so the error between $u$ and $u_a$ decreases until reaching machine precision. But for the other three bases $u_a$ cannot be computed up to machine precision.

\begin{figure}
\begin{center}
\begin{tikzpicture}
\begin{loglogaxis}[height=8cm,xlabel=$h$ (amplitude-based GPW bases), ylabel=max error on disks of radius h,
legend pos=outer north east,
xmin=10^(-7),xmax=20,ymin=10^(-16),ymax=10^(0),
xtick={1,0.01,.0001,.000001,.00000001},
very thick,cycle list name=twotimeseight,grid=major]
\addplot table [x=h, y=errAbGn1]{./BAE3D2OCaseexpikapy50pts.dat};
\addlegendentry{$n=1$}
\addplot table [x=h, y=errAbGn2]{./BAE3D2OCaseexpikapy50pts.dat};
\addlegendentry{$n=2$}
\addplot table [x=h, y=errAbGn3]{./BAE3D2OCaseexpikapy50pts.dat};
\addlegendentry{$n=3$}
\addplot table [x=h, y=errAbGn4]{./BAE3D2OCaseexpikapy50pts.dat};
\addlegendentry{$n=4$}
\addplot table [x=h, y=errAbGn5]{./BAE3D2OCaseexpikapy50pts.dat};
\addlegendentry{$n=5$}
\addplot table [x=h, y=errAbGn6]{./BAE3D2OCaseexpikapy50pts.dat};
\addlegendentry{$n=6$}
\addplot table [x=h, y=errAbGn7]{./BAE3D2OCaseexpikapy50pts.dat};
\addlegendentry{$n=7$}
\addplot table [x=h, y=errAbGn8]{./BAE3D2OCaseexpikapy50pts.dat};
\addlegendentry{$n=8$}
\addplot table [x=h, y=errPWfn1]{./BAE3D2OCaseexpikapy50pts.dat};
\addlegendentry{$n=1$ PW}
\addplot table [x=h, y=errPWfn2]{./BAE3D2OCaseexpikapy50pts.dat};
\addlegendentry{$n=2$ PW}
\addplot table [x=h, y=errPWfn3]{./BAE3D2OCaseexpikapy50pts.dat};
\addlegendentry{$n=3$ PW}
\addplot table [x=h, y=errPWfn4]{./BAE3D2OCaseexpikapy50pts.dat};
\addlegendentry{$n=4$ PW}
\addplot table [x=h, y=errPWfn5]{./BAE3D2OCaseexpikapy50pts.dat};
\addlegendentry{$n=5$ PW}
\addplot table [x=h, y=errPWfn6]{./BAE3D2OCaseexpikapy50pts.dat};
\addlegendentry{$n=6$ PW}
\addplot table [x=h, y=errPWfn7]{./BAE3D2OCaseexpikapy50pts.dat};
\addlegendentry{$n=7$ PW}
\addplot table [x=h, y=errPWfn8]{./BAE3D2OCaseexpikapy50pts.dat};
\addlegendentry{$n=8$ PW}
\addplot[dotted] coordinates {(2*10^-6, 10^-10) (2*10^-3, 10^-4)};
\addlegendentry{order 2}
\addplot[dashed] coordinates {(2*10^-1, 2*10^-10) (2*10^0, 2*10^-1)};
\addlegendentry{order 9}
\end{loglogaxis}
\end{tikzpicture}
\begin{tikzpicture}
\begin{loglogaxis}[height=8cm,xlabel=$h$ (phase-based GPW bases),ylabel=max error on disks of radius h,
legend pos=outer north east,
xmin=10^(-7),xmax=20,ymin=10^(-16),ymax=10^(0),
xtick={1,0.01,.0001,.000001,.00000001},
very thick,cycle list name=twotimeseight,grid=major]
\addplot table [x=h, y=errPbGn1]{./BAE3D2OCaseexpikapy50pts.dat};
\addlegendentry{$n=1$}
\addplot table [x=h, y=errPbGn2]{./BAE3D2OCaseexpikapy50pts.dat};
\addlegendentry{$n=2$}
\addplot table [x=h, y=errPbGn3]{./BAE3D2OCaseexpikapy50pts.dat};
\addlegendentry{$n=3$}
\addplot table [x=h, y=errPbGn4]{./BAE3D2OCaseexpikapy50pts.dat};
\addlegendentry{$n=4$}
\addplot table [x=h, y=errPbGn5]{./BAE3D2OCaseexpikapy50pts.dat};
\addlegendentry{$n=5$}
\addplot table [x=h, y=errPbGn6]{./BAE3D2OCaseexpikapy50pts.dat};
\addlegendentry{$n=6$}
\addplot table [x=h, y=errPbGn7]{./BAE3D2OCaseexpikapy50pts.dat};
\addlegendentry{$n=7$}
\addplot table [x=h, y=errPbGn8]{./BAE3D2OCaseexpikapy50pts.dat};
\addlegendentry{$n=8$}
\addplot table [x=h, y=errPWfn1]{./BAE3D2OCaseexpikapy50pts.dat};
\addlegendentry{$n=1$ PW}
\addplot table [x=h, y=errPWfn2]{./BAE3D2OCaseexpikapy50pts.dat};
\addlegendentry{$n=2$ PW}
\addplot table [x=h, y=errPWfn3]{./BAE3D2OCaseexpikapy50pts.dat};
\addlegendentry{$n=3$ PW}
\addplot table [x=h, y=errPWfn4]{./BAE3D2OCaseexpikapy50pts.dat};
\addlegendentry{$n=4$ PW}
\addplot table [x=h, y=errPWfn5]{./BAE3D2OCaseexpikapy50pts.dat};
\addlegendentry{$n=5$ PW}
\addplot table [x=h, y=errPWfn6]{./BAE3D2OCaseexpikapy50pts.dat};
\addlegendentry{$n=6$ PW}
\addplot table [x=h, y=errPWfn7]{./BAE3D2OCaseexpikapy50pts.dat};
\addlegendentry{$n=7$ PW}
\addplot table [x=h, y=errPWfn8]{./BAE3D2OCaseexpikapy50pts.dat};
\addlegendentry{$n=8$ PW}
\addplot[dotted] coordinates {(2*10^-6, 10^-10) (2*10^-3, 10^-4)};
\addlegendentry{order 2}
\addplot[dashed] coordinates {(2*10^-1, 2*10^-10) (2*10^0, 2*10^-1)};
\addlegendentry{order 9}
\end{loglogaxis}
\end{tikzpicture}
\end{center}
\caption{Local approximation of an exact solution from quasi-Trefftz bases: convergence results for the first test case, where PW functions are exact solutions, for $n$ from $1$ to $8$. For each value of $n$, the expected order of convergence, namely $n+1$, is observed and the error decreases until it reaches a threshold. Comparison of a PW and GPW bases using the same initialization, for amplitude-based (top) and phase-based (bottom) GPWs.}
\label{fig:HelmGPW}
\end{figure}
\begin{figure}
\begin{center}
\begin{tikzpicture}
\begin{loglogaxis}[height=8cm,xlabel=$h$, ylabel=max error on disks of radius h,
legend pos=outer north east,
xmin=10^(-7),xmax=20,ymin=10^(-16),ymax=10^(0),
xtick={1,0.01,.0001,.000001,.00000001},
very thick,cycle list name=twotimeseight,grid=major]
\addplot table [x=h, y=errPstn1]{./BAE3D2OCaseexpikapy50pts.dat};
\addlegendentry{$n=1$}
\addplot table [x=h, y=errPstn2]{./BAE3D2OCaseexpikapy50pts.dat};
\addlegendentry{$n=2$}
\addplot table [x=h, y=errPstn3]{./BAE3D2OCaseexpikapy50pts.dat};
\addlegendentry{$n=3$}
\addplot table [x=h, y=errPstn4]{./BAE3D2OCaseexpikapy50pts.dat};
\addlegendentry{$n=4$}
\addplot table [x=h, y=errPstn5]{./BAE3D2OCaseexpikapy50pts.dat};
\addlegendentry{$n=5$}
\addplot table [x=h, y=errPstn6]{./BAE3D2OCaseexpikapy50pts.dat};
\addlegendentry{$n=6$}
\addplot table [x=h, y=errPstn7]{./BAE3D2OCaseexpikapy50pts.dat};
\addlegendentry{$n=7$}
\addplot table [x=h, y=errPstn8]{./BAE3D2OCaseexpikapy50pts.dat};
\addlegendentry{$n=8$}
\addplot table [x=h, y=errPWfn1]{./BAE3D2OCaseexpikapy50pts.dat};
\addlegendentry{$n=1$ PW}
\addplot table [x=h, y=errPWfn2]{./BAE3D2OCaseexpikapy50pts.dat};
\addlegendentry{$n=2$ PW}
\addplot table [x=h, y=errPWfn3]{./BAE3D2OCaseexpikapy50pts.dat};
\addlegendentry{$n=3$ PW}
\addplot table [x=h, y=errPWfn4]{./BAE3D2OCaseexpikapy50pts.dat};
\addlegendentry{$n=4$ PW}
\addplot table [x=h, y=errPWfn5]{./BAE3D2OCaseexpikapy50pts.dat};
\addlegendentry{$n=5$ PW}
\addplot table [x=h, y=errPWfn6]{./BAE3D2OCaseexpikapy50pts.dat};
\addlegendentry{$n=6$ PW}
\addplot table [x=h, y=errPWfn7]{./BAE3D2OCaseexpikapy50pts.dat};
\addlegendentry{$n=7$ PW}
\addplot table [x=h, y=errPWfn8]{./BAE3D2OCaseexpikapy50pts.dat};
\addlegendentry{$n=8$ PW}
\addplot[dotted] coordinates {(2*10^-6, 10^-10) (2*10^-3, 10^-4)};
\addlegendentry{order 2}
\addplot[dashed] coordinates {(2*10^-1, 2*10^-10) (2*10^0, 2*10^-1)};
\addlegendentry{order 9}
\end{loglogaxis}
\end{tikzpicture}
\end{center}
\caption{Local approximation of an exact solution from quasi-Trefftz bases: convergence results for the first test case, where PW functions are exact solutions, for $n$ from $1$ to $8$. For each value of $n$, the expected order of convergence, namely $n+1$, is observed.  
Comparison of a PW and the polynomial quasi-Trefftz bases using the same initialization, evidencing the absence of conditioning problem with the polynomial basis.}
\label{fig:HelmpqT}
\end{figure}
\begin{figure}
\begin{center}
\begin{tikzpicture}
\begin{loglogaxis}[height=12cm,xlabel=$h$, ylabel=max error on disks of radius h,
legend pos=outer north east,
xmin=10^(-7),xmax=20,ymin=10^(-16),ymax=10^(0),
xtick={1,0.01,.0001,.000001,.00000001},
very thick,cycle list name=threetimeseight,grid=major]
\addplot table [x=h, y=errAbGn1]{./BAE3D2OCaseKAPPAairyxexpiypz50pts.dat};
\addlegendentry{$n=1$ Amp}
\addplot table [x=h, y=errAbGn2]{./BAE3D2OCaseKAPPAairyxexpiypz50pts.dat};
\addlegendentry{$n=2$ Amp}
\addplot table [x=h, y=errAbGn3]{./BAE3D2OCaseKAPPAairyxexpiypz50pts.dat};
\addlegendentry{$n=3$ Amp}
\addplot table [x=h, y=errAbGn4]{./BAE3D2OCaseKAPPAairyxexpiypz50pts.dat};
\addlegendentry{$n=4$ Amp}
\addplot table [x=h, y=errAbGn5]{./BAE3D2OCaseKAPPAairyxexpiypz50pts.dat};
\addlegendentry{$n=5$ Amp}
\addplot table [x=h, y=errAbGn6]{./BAE3D2OCaseKAPPAairyxexpiypz50pts.dat};
\addlegendentry{$n=6$ Amp}
\addplot table [x=h, y=errAbGn7]{./BAE3D2OCaseKAPPAairyxexpiypz50pts.dat};
\addlegendentry{$n=7$ Amp}
\addplot table [x=h, y=errAbGn8]{./BAE3D2OCaseKAPPAairyxexpiypz50pts.dat};
\addlegendentry{$n=8$ Amp}
\addplot table [x=h, y=errPbGn1]{./BAE3D2OCaseKAPPAairyxexpiypz50pts.dat};
\addlegendentry{$n=1$ Pha}
\addplot table [x=h, y=errPbGn2]{./BAE3D2OCaseKAPPAairyxexpiypz50pts.dat};
\addlegendentry{$n=2$ Pha}
\addplot table [x=h, y=errPbGn3]{./BAE3D2OCaseKAPPAairyxexpiypz50pts.dat};
\addlegendentry{$n=3$ Pha}
\addplot table [x=h, y=errPbGn4]{./BAE3D2OCaseKAPPAairyxexpiypz50pts.dat};
\addlegendentry{$n=4$ Pha}
\addplot table [x=h, y=errPbGn5]{./BAE3D2OCaseKAPPAairyxexpiypz50pts.dat};
\addlegendentry{$n=5$ Pha}
\addplot table [x=h, y=errPbGn6]{./BAE3D2OCaseKAPPAairyxexpiypz50pts.dat};
\addlegendentry{$n=6$ Pha}
\addplot table [x=h, y=errPbGn7]{./BAE3D2OCaseKAPPAairyxexpiypz50pts.dat};
\addlegendentry{$n=7$ Pha}
\addplot table [x=h, y=errPbGn8]{./BAE3D2OCaseKAPPAairyxexpiypz50pts.dat};
\addlegendentry{$n=8$ Pha}
\addplot table [x=h, y=errPstn1]{./BAE3D2OCaseKAPPAairyxexpiypz50pts.dat};
\addlegendentry{$n=1$ poly}
\addplot table [x=h, y=errPstn2]{./BAE3D2OCaseKAPPAairyxexpiypz50pts.dat};
\addlegendentry{$n=2$ poly}
\addplot table [x=h, y=errPstn3]{./BAE3D2OCaseKAPPAairyxexpiypz50pts.dat};
\addlegendentry{$n=3$ poly}
\addplot table [x=h, y=errPstn4]{./BAE3D2OCaseKAPPAairyxexpiypz50pts.dat};
\addlegendentry{$n=4$ poly}
\addplot table [x=h, y=errPstn5]{./BAE3D2OCaseKAPPAairyxexpiypz50pts.dat};
\addlegendentry{$n=5$ poly}
\addplot table [x=h, y=errPstn6]{./BAE3D2OCaseKAPPAairyxexpiypz50pts.dat};
\addlegendentry{$n=6$ poly}
\addplot table [x=h, y=errPstn7]{./BAE3D2OCaseKAPPAairyxexpiypz50pts.dat};
\addlegendentry{$n=7$ poly}
\addplot table [x=h, y=errPstn8]{./BAE3D2OCaseKAPPAairyxexpiypz50pts.dat};
\addlegendentry{$n=8$ poly}
\addplot[dotted] coordinates {(2*10^-6, 10^-10) (2*10^-3, 10^-4)};
\addlegendentry{order 2}
\addplot[dashed] coordinates {(2*10^-1, 2*10^-10) (2*10^0, 2*10^-1)};
\addlegendentry{order 9}
\end{loglogaxis}
\end{tikzpicture}
\end{center}
\caption{Local approximation of an exact solution from quasi-Trefftz bases: convergence results for the second test case, for $n$ from $1$ to $8$. For each value of $n$, the expected order of convergence, namely $n+1$, is observed.
Comparison of an amplitude-based GPW basis (Amp), a phase-based GPW basis (Pha) and a polyonomial quasi-Trefftz basis (poly), evidencing again the absence of conditioning problem with the polynomial basis. Only with the quasi-Trefftz polynomial approximation does the error decrease until it reaches machine precision.}
\label{fig:ConvHelm-1}
\end{figure}
\begin{figure}
\begin{center}
\begin{tikzpicture}
\begin{loglogaxis}[height=12cm,xlabel=$h$, ylabel=max error on disks of radius h,
legend pos=outer north east,
xmin=10^(-7),xmax=20,ymin=10^(-16),ymax=10^(0),
xtick={1,0.01,.0001,.000001,.00000001},
very thick,cycle list name=threetimeseight,grid=major]
\addplot table [x=h, y=errAbGn1]{./BAE3D2OCaseKAPPAConvHelm50pts.dat};
\addlegendentry{$n=1$ Amp}
\addplot table [x=h, y=errAbGn2]{./BAE3D2OCaseKAPPAConvHelm50pts.dat};
\addlegendentry{$n=2$ Amp}
\addplot table [x=h, y=errAbGn3]{./BAE3D2OCaseKAPPAConvHelm50pts.dat};
\addlegendentry{$n=3$ Amp}
\addplot table [x=h, y=errAbGn4]{./BAE3D2OCaseKAPPAConvHelm50pts.dat};
\addlegendentry{$n=4$ Amp}
\addplot table [x=h, y=errAbGn5]{./BAE3D2OCaseKAPPAConvHelm50pts.dat};
\addlegendentry{$n=5$ Amp}
\addplot table [x=h, y=errAbGn6]{./BAE3D2OCaseKAPPAConvHelm50pts.dat};
\addlegendentry{$n=6$ Amp}
\addplot table [x=h, y=errAbGn7]{./BAE3D2OCaseKAPPAConvHelm50pts.dat};
\addlegendentry{$n=7$ Amp}
\addplot table [x=h, y=errAbGn8]{./BAE3D2OCaseKAPPAConvHelm50pts.dat};
\addlegendentry{$n=8$ Amp}
\addplot table [x=h, y=errPbGn1]{./BAE3D2OCaseKAPPAConvHelm50pts.dat};
\addlegendentry{$n=1$ Pha}
\addplot table [x=h, y=errPbGn2]{./BAE3D2OCaseKAPPAConvHelm50pts.dat};
\addlegendentry{$n=2$ Pha}
\addplot table [x=h, y=errPbGn3]{./BAE3D2OCaseKAPPAConvHelm50pts.dat};
\addlegendentry{$n=3$ Pha}
\addplot table [x=h, y=errPbGn4]{./BAE3D2OCaseKAPPAConvHelm50pts.dat};
\addlegendentry{$n=4$ Pha}
\addplot table [x=h, y=errPbGn5]{./BAE3D2OCaseKAPPAConvHelm50pts.dat};
\addlegendentry{$n=5$ Pha}
\addplot table [x=h, y=errPbGn6]{./BAE3D2OCaseKAPPAConvHelm50pts.dat};
\addlegendentry{$n=6$ Pha}
\addplot table [x=h, y=errPbGn7]{./BAE3D2OCaseKAPPAConvHelm50pts.dat};
\addlegendentry{$n=7$ Pha}
\addplot table [x=h, y=errPbGn8]{./BAE3D2OCaseKAPPAConvHelm50pts.dat};
\addlegendentry{$n=8$ Pha}
\addplot table [x=h, y=errPstn1]{./BAE3D2OCaseKAPPAConvHelm50pts.dat};
\addlegendentry{$n=1$ poly}
\addplot table [x=h, y=errPstn2]{./BAE3D2OCaseKAPPAConvHelm50pts.dat};
\addlegendentry{$n=2$ poly}
\addplot table [x=h, y=errPstn3]{./BAE3D2OCaseKAPPAConvHelm50pts.dat};
\addlegendentry{$n=3$ poly}
\addplot table [x=h, y=errPstn4]{./BAE3D2OCaseKAPPAConvHelm50pts.dat};
\addlegendentry{$n=4$ poly}
\addplot table [x=h, y=errPstn5]{./BAE3D2OCaseKAPPAConvHelm50pts.dat};
\addlegendentry{$n=5$ poly}
\addplot table [x=h, y=errPstn6]{./BAE3D2OCaseKAPPAConvHelm50pts.dat};
\addlegendentry{$n=6$ poly}
\addplot table [x=h, y=errPstn7]{./BAE3D2OCaseKAPPAConvHelm50pts.dat};
\addlegendentry{$n=7$ poly}
\addplot table [x=h, y=errPstn8]{./BAE3D2OCaseKAPPAConvHelm50pts.dat};
\addlegendentry{$n=8$ poly}
\addplot[dotted] coordinates {(2*10^-6, 10^-10) (2*10^-3, 10^-4)};
\addlegendentry{order 2}
\addplot[dashed] coordinates {(2*10^-1, 2*10^-10) (2*10^0, 2*10^-1)};
\addlegendentry{order 9}
\end{loglogaxis}
\end{tikzpicture}
\end{center}
\caption{Local approximation of an exact solution from quasi-Trefftz bases: convergence results for the third test case, for $n$ from $1$ to $8$. For each value of $n$, the expected order of convergence, namely $n+1$, is observed.
Comparison of an amplitude-based GPW basis (Amp), a phase-based GPW basis (Pha) and a polyonomial quasi-Trefftz basis (poly), evidencing again the absence of conditioning problem with the polynomial basis. Only with the quasi-Trefftz polynomial approximation does the error decrease until it reaches machine precision.}
\label{fig:ConvHelm-2}
\end{figure}

\section{Conclusion}
Given a partial differential operator $\mathcal L$ and a parameter $q\in\mathbb N$, a quasi-Trefftz function $f$ is a function satisfying the following local property in the neighborhood of a given point $\pG$:
\begin{equation}
\label{eq:finalQTprop}
\forall \varx \text{ in a  neighborhood of }\pG,\ \mathcal L f(\varx) = O\left(|\varx-\pG|^q\right).
\end{equation}

The work presented here may be summarized as follows.
\begin{itemize}
\item We introduced three families of quasi-Trefftz functions for a class of 3D PDEs including the convected Helmholtz equation.
Two of these, the GPWs, are generalizations of a PW ansatz, and were first introduced for a class of 2D problems. 
The ansatz defining these two families are defined in a neighborhood of a given point $\pG$ under the following form: for all $\varx$,
\begin{equation*}
\left\{
\begin{array}{l}
\displaystyle
J(\varx) := \polQ(\varx-\pG)\exp \Lambda \cdot\Big(\varx-\pG \Big) 
\text{ for some polynomial } Q \text{ and some vector }\Lambda\in\mathbb C^3,\\
\displaystyle
G(\varx) := \exp \polP(\varx-\pG)
\text{ for some polynomial } P.
%,\\
%\displaystyle
%H(\varx) := \polR(\varx-\pG)
%\text{ for some polynomial } R,\\
\end{array}
\right.
\end{equation*}
The third one is fully polynomial, and this is the first introduction and study of polynomial quasi-Trefftz functions. The corresponding ansatz is defined in a neighborhood of a given point $\pG$ under the following forms: for all $\varx$,
\begin{equation*}
H(\varx) := \polR(\varx-\pG)
%\text{ with } R:=\sum_{i\in\mathbb N_0^3, |i|\leq d} \nu_i \mathbf X^i.
\text{ for some polynomial } R.
\end{equation*}
\item We provide explicit algorithms to construct  quasi-Trefftz functions belonging to each of these three families, by constructing the corresponding polynomial $P$, $Q$ or $R$. This is achieved by a careful study of the system obtained by setting to zero the degree-$(q-1)$ Taylor polynomial of the image of  each ansatz through the operator $\mathcal L$.
It is then possible to choose adequately the degree of the polynomial $P$, $Q$ or $R$ in order to split this system  into a hierarchy of linear triangular sub-system. Hence a solution to the system can be computed via an explicit formula by simply solving successively  the subsystems by substitution.
Interestingly, some of the polynomial coefficients of $P$, $Q$ and $R$ are free in the resulting algorithms, and thanks to these sets of linearly independent quasi-Trefftz functions can be constructed.
As a consequence, beyond the construction of individual quasi-Trefftz functions, we can construct spaces of quasi-Trefftz functions.
\item 
We prove that it is possible to construct quasi-Trefftz spaces $\mathbb V_h$ spanned by sets of such quasi-Trefftz functions enjoying high order approximation property for exact solutions of the PDE.  
More precisely, given a given point $\pG$, \Cref{thm:approx} states that in order to achieve a given order of accuracy $n+1$ of local approximation property in the following sense:
\begin{equation}
\label{eq:cclapp}
\forall u%\in \mathcal C^\infty 
\text{ satisfying the governing PDE, }%\mathcal Lu=0,
\exists u_a\in\mathbb V_h \text{ s. t. } \forall h \text{ small enough }
\|u-u_a\|_{L^\infty(B(\pG,h))}\leq C h^{n+1},
\end{equation}
where $B(\pG,h)$ denotes the sphere centered at $\pG$ of radius $h$ in $\mathbb R^3$, it is sufficient to construct a quasi-Trefftz space $\mathbb V_h$ of dimension $p=(n+1)^2$ with basis functions satisfying the quasi-Trefftz property \eqref{eq:finalQTprop} with $q=\max(n-1,1)$.
For reference, reaching the same order of approximation in  \eqref{eq:cclapp} using a standard polynomial space would yield a dimension $p=(n+1)(n+2)(n+3)/6$ (corresponding to the full space of polynomials of degrees at most equal to $n$).
For instance for $n=8$ the quasi-Trefftz space is of dimension $81$ while the polynomial space is of dimension $165$.
Besides, as noted in \Cref{thm:approx}, %precisely because it relies on matching the Taylor expansions of $u$ and $u_a$,
we actually show a convergence in $H^1$:
\begin{equation}
\forall u%\in \mathcal C^\infty 
\text{ satisfying the governing PDE, }%\mathcal Lu=0,
\exists u_a\in\mathbb V_h \text{ s. t. } \forall h \text{ small enough }
\|u-u_a\|_{H^1(B(\pG,h))}\leq C h^{n}.
\end{equation}
\item Most importantly, {the concept of polynomial quasi-Trefftz basis, explored here for the first time in the context of time-harmonic problems, represents a new avenue to leverage the benefits of Trefftz methods while avoiding the ill-conditioning problem inherent to wave-like basis. This problem has been a clear limitation to the further development of Trefftz methods. Moreover, since they do not rely on a wave-like ansatz, polynomial quasi-Trefftz bases do not only represent a promising way forward for wave propagation, but also can be applied beyond that to problems governed by other types of PDEs. 
However, the work presented here focuses on approximation properties at a single point, and does not tackle the $p$-regime, when $h$ is fixed and the number of basis functions $p$ is increased.}
\end{itemize}

Future plans include :
\begin{itemize}
\item 
investigating, beyond local properties, global best approximation error on a domain, by performing local approximation of each element of a mesh of the domain, and evaluating an $H^1$ error,
\item
comparing the performance of the three quasi-Trefftz bases with standard polynomial bases in terms of accuracy, computing time and stability, in particular in the high-frequency regime, as GPWs might retain an edge compared to polynomials on a sphere containing several wavelengths,
\item comparing the performance of quasi-Trefftz methods with these three quasi-Trefftz bases to other polynomial and wave-based methods - standardly used in the literature in aero-acoustics - on realistic industrial test cases, 
%including the impact of the choice of basis functions on the number and size of elements and the number of unknowns
%how does the choice of basis functions influence/impact on the 
including procedures to balance the size of mesh elements with the number of basis functions per element, and how it can be impacted by the type of basis functions.
\end{itemize}

%%%%%%%%%%%%%%%%%%%%%
\appendix
%%%%%%%%%%%%%%%%%%%%%

\section{A non-singular  matrix statement}
\label{app:eigval}
We will show here that under the assumption that $|\mathbf M(\pG)|<1$ then  the matrix
$$
\Cmat :=\begin{bmatrix}
\left( \mathbf M_1(\pG)\right)^2 -1 & \frac12 \mathbf M_1(\pG)\mathbf M_2(\pG)& \frac12\mathbf M_1(\pG)\mathbf M_3(\pG)\\
 \frac12\mathbf M_1(\pG)\mathbf M_2(\pG)&\left( \mathbf M_2(\pG)\right)^2-1& \frac12\mathbf M_2(\pG)\mathbf M_3(\pG)\\
  \frac12\mathbf M_1(\pG)\mathbf M_3(\pG) &\frac12 \mathbf M_2(\pG)\mathbf M_3(\pG)&\left( \mathbf M_3(\pG)\right)^2-1
\end{bmatrix}
$$
is not singular. For the sake of compactness we remove in this demonstration the dependecy of the entries of $\mathbf M$ on $\pG$.

First, we compute the determinant of $\Cmat$:
$$
\begin{array}{rl}
\det (\Cmat)
%%&= (\left( \mathbf M_1(\pG)\right)^2-1)\det
%%\begin{bmatrix}  
%%\left(\mathbf M_2(\pG)\right)^2-1 & \frac12\mathbf M_2(\pG)\mathbf M_3(\pG)
%%\\
%%\frac12 \mathbf M_2(\pG)\mathbf M_3(\pG) & \left( \mathbf M_3(\pG)\right)^2-1
%%\end{bmatrix}
%%\\&
%%\phantom{=}- \frac12\mathbf M_1(\pG)\mathbf M_2(\pG)\det
%%\begin{bmatrix}  
%% \frac12 \mathbf M_1(\pG)\mathbf M_2(\pG)& \frac12\mathbf M_1(\pG)\mathbf M_3(\pG)
%%\\
%%\frac12 \mathbf M_2(\pG)\mathbf M_3(\pG)&\left( \mathbf M_3(\pG)\right)^2-1
%%\end{bmatrix}
%%\\&
%%\phantom{=} +   \frac12\mathbf M_1(\pG)\mathbf M_3(\pG)   \det
%%\begin{bmatrix}  
%% \frac12 \mathbf M_1(\pG)\mathbf M_2(\pG)& \frac12\mathbf M_1(\pG)\mathbf M_3(\pG)
%%\\
%%\left( \mathbf M_2(\pG)\right)^2-1& \frac12\mathbf M_2(\pG)\mathbf M_3(\pG)
%%\end{bmatrix}
%%\\
%%&= (\left( \mathbf M_1(\pG)\right)^2-1)\left(1% for 1^2
%%-\left(\mathbf M_2(\pG)^2+\mathbf M_3(\pG)^2\right)% times 1
%%+\frac34 \mathbf M_2(\pG)^2\mathbf M_3(\pG)^2\right)
%%\\&\phantom { = }
%%- \frac14 \mathbf M_1(\pG)^2\mathbf M_2(\pG)^2\left( \frac12 \mathbf M_3(\pG)^2-1 \right)
%%+\frac14 \mathbf M_1(\pG)^2\mathbf M_3(\pG)^2\left(1- \frac12 \mathbf M_2(\pG)^2\right)
%%\\
%&=
%\frac12\left( 
%\mathbf M_1%(\pG)
%\mathbf M_2%(\pG)
%\mathbf M_3%(\pG)
% \right)^2
%+\left(\mathbf M_1%(\pG)
%^2+\mathbf M_2%(\pG)
%^2+\mathbf M_3%(\pG)
%^2\right)% times (lambda=1)^2
%-1% for (lambda=1)^3
%%\\&\phantom { = }
%-\frac34\left( \mathbf M_1%(\pG)
%^2\mathbf M_2%(\pG)
%^2+\mathbf M_1%(\pG)
%^2\mathbf M_3%(\pG)
%^2+\mathbf M_2%(\pG)
%^2\mathbf M_3%(\pG)
%^2 \right),%\times (lambda=1)
%\\
&=
\left[{\mathbf M_1}^2+{\mathbf M_2}^2+{\mathbf M_3}^2 - 1 \right]
+ \left[ \frac12\left( \mathbf M_1\mathbf M_2\mathbf M_3 \right)^2-\frac34\left( {\mathbf M_1}^2{\mathbf M_2}^2+{\mathbf M_1}^2{\mathbf M_3}^2+{\mathbf M_2}^2{\mathbf M_3}^2 \right) \right] .
\end{array}
$$
We will prove that $\det(\Cmat)<0$ by showing that both brackets in this last expression are negative, the first one strictly:
%Hence $\det (\mathcal C) = 0$ if and only if
%$$ 
%1
%+\frac34\left( \mathbf M_1^2\mathbf M_2^2+\mathbf M_1^2\mathbf M_3^2+\mathbf M_2^2\mathbf M_3^2 \right)
%=
%\frac12\left( \mathbf M_1\mathbf M_2\mathbf M_3\right)^2
%+\left(\mathbf M_1^2+\mathbf M_2^2+\mathbf M_3^2\right).
%$$
%%$$
%%\begin{array}{l}
%%\displaystyle 
%%1
%%+\frac34\left( \mathbf M_1(\pG)^2\mathbf M_2(\pG)^2+\mathbf M_1(\pG)^2\mathbf M_3(\pG)^2+\mathbf M_2(\pG)^2\mathbf M_3(\pG)^2 \right)
%%\\\displaystyle 
%%=
%%\frac12\left( \mathbf M_1(\pG)\mathbf M_2(\pG)\mathbf M_3(\pG) \right)^2
%%+\left(\mathbf M_1(\pG)^2+\mathbf M_2(\pG)^2+\mathbf M_3(\pG)^2\right).
%%\end{array}
%%$$
\begin{enumerate}
\item $1>|\mathbf M(\pG)|$ by assumption, so $\left[ {\mathbf M_1}^2+{\mathbf M_2}^2+{\mathbf M_3}^2-1\right]<0$;
\item 
defining 
$a_1 = {\mathbf M_1}^2{\mathbf M_2}^2$,
$a_2 = {\mathbf M_1}^2{\mathbf M_3}^2$ and
$a_3 = {\mathbf M_2}^2{\mathbf M_3}^2$, then the following MacLaurin's inequality:
$$
\frac{a_1+a_2+a_3}3 \geq \sqrt[3]{ a_1a_2a_3 }
$$
shows that:
$$
\frac{{\mathbf M_1}^2{\mathbf M_2}^2+{\mathbf M_1}^2{\mathbf M_3}^2+{\mathbf M_2}^2{\mathbf M_3}^2}3 
\geq
| \mathbf M_1\mathbf M_2\mathbf M_3|^{4/3},
$$
combined with:
 $$
 \left\{\begin{array}{l}
 \displaystyle \frac12\geq\frac13 \text{ and } 1\geq \frac13
  \\
  \displaystyle
|  \mathbf M_k|^{4/3}\geq|\mathbf M_k|^{2}\ \forall k\in\{1,2,3\} \text{ since } | \mathbf M_k|<1,
  \end{array}\right.
  $$
 in turn shows that:
  $$
  \frac{{\mathbf M_1}^2{\mathbf M_2}^2+{\mathbf M_1}^2{\mathbf M_3}^2+{\mathbf M_2}^2{\mathbf M_3}^2}2
\geq
\frac{| \mathbf M_1\mathbf M_2\mathbf M_3|^{2}}{3},
  $$
  and therefore the second bracket in $\det(\Cmat)$ is negative:
  $$
 \frac1{2}| \mathbf M_1\mathbf M_2\mathbf M_3|^{2}- \frac{3}{4}\left({\mathbf M_1}^2{\mathbf M_2}^2+{\mathbf M_1}^2{\mathbf M_3}^2+{\mathbf M_2}^2{\mathbf M_3}^2\right)
\leq
0.
  $$
\end{enumerate}
This actually shows that $\det (\Cmat)<0$, which indeed proves that the matrix $\Cmat$ is not singular. 

%%%%%%%%%%%%%%%%%%%%%
\section{Spherical harmonics reminder}
\label{sec:SHrem}%Spherical Harmonic reminder
Legendre polynomials are defined on $\mathbb R$ for $m\in\mathbb N_0$ as:
$$
P_l^{m}(x) = \frac{(-1)^m}{2^l l!} (1-x^2)^{m/2} \partial_x^{l+m} (x^2-1)^{l}, \forall x\in\mathbb R,
$$
$$
P_l^{-m} = (-1)^m\frac{(l-m)!}{(l+m)!}P_l^m.
$$
%So to express them explicitly, we note that: 
Moreover, for all $l$ and $m$ in $\mathbb N_0$ we note that:
$$
(x^2-1)^l = \sum_{l'=0}^l \begin{pmatrix} l\\l'\end{pmatrix} (-1)^{l-l'} x^{2l'},
$$
$$
\partial_x^{l+m}(x^2-1)^l = \sum_{l'=\left\lfloor\frac{l+m}2\right\rfloor+1}^l 
\begin{pmatrix} l\\l'\end{pmatrix} (-1)^{l-l'} \frac{(2l')!}{(2l'-l-m)!}x^{2l'-l-m}.
$$
As a result, from which we can write for all $\theta\in\mathbb R$:
$$
\text{if }m\geq0, P_l^m(\cos\theta) = \frac{(-1)^m}{2^ll!} (\sin\theta)^m
\sum_{l'=\left\lfloor\frac{l+m}2\right\rfloor+1}^l 
\begin{pmatrix} l\\l'\end{pmatrix} (-1)^{l-l'} \frac{(2l')!}{(2l'-l-m)!}(\cos\theta)^{2l'-l-m},
$$
$$
\text{if }m<0, P_l^m(\cos\theta) = 
  \frac{(l-|m|)!}{2^ll!(l+|m|)!} 
(\sin\theta)^m
\sum_{l'=\left\lfloor\frac{l+|m|}2\right\rfloor+1}^l 
\begin{pmatrix} l\\l'\end{pmatrix} (-1)^{l-l'} \frac{(2l')!}{(2l'-l-|m|)!}(\cos\theta)^{2l'-l-|m|}.
$$

%%%%%%%%%%%%%%%%%%%%%
\section{Faa di Bruno formula in 3D}
\label{App:FdB}
In dimension three, the Faa di Bruno formula presented in \cite{CS} reads: if $f$ is a function of one real variable and $g$ is a function defined on $\mathbb R^3$,
$\forall x\in\mathbb R^3$
$$
\partial_x^i f(g(x)) = 
\sum_{1\leq \widetilde m\leq |i|}
f^{(\widetilde m)}(g(x)) \sum_{s=1}^{|i|}\sum_{p_s(i,\widetilde m)} i!\prod_{m=1}^s 
\frac{1}{k_m!}
\left(\frac{1}{l_m!} \partial_x^{l_m} g(x) \right)^{k_m}
$$
where the linear order $\prec$ on $\mathbb N_0^3$ is defined in the introduction,
%by	
%$$
%\forall (\mu,\nu)\in\left(\mathbb N_0^3\right)^2,
% \mu\prec\nu 
% \Leftrightarrow
%\left\{
% \begin{array}{l}
% |\mu|<|\nu|\ , \text{ or }\\
% |\mu|=|\nu| \text{ and } \mu_1<\nu_1\ , \text{ or }\\
% |\mu|=|\nu|,\  \mu_1=\nu_1\text{ and } \mu_2<\nu_2\ ,
% \end{array}
% \right.
%$$
while the partition $p_s(i,\widetilde m)$ of multi-index $i\in\mathbb N^3$ is defined by:
$$
p_s(i,\widetilde m) = \left\{ (k_1,\dots,k_s;l_1,\dots,l_s) ; k_i>0, 0\prec l_1\prec\dots\prec l_s, \sum_{m=1}^s k_i = \widetilde m, \sum_{m=1}^s k_m l_m = i\right\}.
$$

\end{document}